\documentclass[12pt,bibliography=totoc,titlepage]{article}

\usepackage{graphics}
\usepackage{epsfig}
\usepackage{epstopdf}
\usepackage{slashed}
\usepackage{enumerate}
\usepackage{cyrillic}
\usepackage{scalerel}
\usepackage{epic}
\usepackage{indentfirst}

\usepackage{etoolbox}

\makeatletter
\pretocmd{\chapter}{\addtocontents{toc}{\protect\addvspace{15\p@}}}{}{}
\pretocmd{\section}{\addtocontents{toc}{\protect\addvspace{25\p@}}}{}{}
\pretocmd{\subsection}{\addtocontents{toc}{\protect\addvspace{15\p@}}}{}{}
\makeatother

\DeclareFontFamily{U}{mathb}{\hyphenchar\font45}
\DeclareFontShape{U}{mathb}{m}{n}{
      <5> <6> <7> <8> <9> <10> gen * mathb
      <10.95> mathb 10 <12> <14.4> <17.28> <20.74> <24.88> mathb12
      }{}
\DeclareSymbolFont{mathb}{U}{mathb}{m}{n}
\DeclareMathSymbol{\precneq}{3}{mathb}{"AC}

\usepackage[usenames,dvipsnames]{color}
\usepackage{float}
\floatstyle{plaintop}
\restylefloat{table}
\usepackage{mathtools}
\usepackage{afterpage}
\newcommand\blankpage{%
    \null
    \thispagestyle{empty}%
    \addtocounter{page}{-1}%
    \newpage}

\usepackage{amsthm}
\usepackage{amssymb}
\usepackage{amsmath}
\usepackage{verbatim}
\usepackage{makeidx}
\makeindex
\usepackage[totoc]{idxlayout}

\usepackage{calrsfs}
\DeclareMathAlphabet{\pazocal}{OMS}{zplm}{m}{n}

\DeclareMathAlphabet{\mathpzc}{OT1}{pzc}{m}{it}
\DeclareMathAlphabet{\mathcalligra}{T1}{calligra}{m}{n} 

\usepackage{txfonts}
\usepackage{lmodern}

\usepackage{hyperref}

\usepackage[T1]{fontenc}
\usepackage{calligra}
\newcommand{\mathscr}[1]{\text{\!\calligra{#1}\,}}
\usepackage{xcolor}
\makeatletter
\let\ps@plainorig\ps@plain
\makeatother

\newcommand\swabfamily{\usefont{U}{yswab}{m}{n}}
\DeclareTextFontCommand{\textswab}{\swabfamily}
\newcommand\frakfamily{\usefont{U}{yfrak}{m}{n}}
\DeclareTextFontCommand{\textfrak}{\frakfamily}
\newcommand\gothfamily{\usefont{U}{ygoth}{m}{n}}
\DeclareTextFontCommand{\textgoth}{\gothfamily}

\newcommand{\gllie}{\text{\swabfamily gl}}
\newcommand{\sllie}{\text{\swabfamily s:l}}
\newcommand{\solie}{\text{\swabfamily s:o}}
\newcommand{\splie}{\text{\swabfamily s:p}}

\newcommand{\Ulie}{\text{\swabfamily U}}
\newcommand{\Mlie}{\text{\swabfamily M}}
\newcommand{\Nlie}{\text{\swabfamily N}}
\newcommand{\Llie}{\text{\swabfamily L}}
\newcommand{\Plie}{\text{\swabfamily P}}

\newcommand{\Ilie}{\text{\swabfamily I}}
\newcommand{\Zlie}{\text{\swabfamily Z}}
\newcommand{\Klie}{\text{\swabfamily K}}
\newcommand{\Vlie}{\text{\swabfamily V}}
\newcommand{\tendsto}{\boldsymbol{\to}}

\newcommand{\gfrak}{\mathfrak{g}}
\newcommand{\hfrak}{\mathfrak{h}}
\newcommand{\kfrak}{\mathfrak{k}}
\newcommand{\bfrak}{\mathfrak{b}}
\newcommand{\nfrak}{\mathfrak{n}}
\newcommand{\pfrak}{\mathfrak{p}}

\newcommand{\cyrsl}{\fontencoding{OT2}\selectfont\textcyrsl}

\newcommand{\tcyr}{{\cyrsl{t}}}

\usepackage[refpage,intoc]{nomencl}
\usepackage{ifthen}
\usepackage{lipsum}

\makenomenclature

\renewcommand{\nomgroup}[1]
{%
\ifthenelse{\equal{#1}{A}}%
{\vspace{3\parsep}\item[]\hspace*{-\leftmargin}%
{\textbf{Abbreviations}}}%
{%
\ifthenelse{\equal{#1}{B}}%
{\item[]\hspace*{-\leftmargin}%
{\textbf{Mathematical Notations}}}%
}
}

\DeclareMathOperator{\Hom}{Hom}
\DeclareMathOperator{\Soc}{Soc}

\DeclareMathOperator{\End}{End}

\DeclareMathOperator{\im}{im}
\DeclareMathOperator{\ch}{ch}
\DeclareMathOperator{\Res}{Res}
\newcommand{\id}{\text{id}}

\newcommand{\dynkinradius}{.1cm}
\newcommand{\dynkinstep}{1cm}
\newcommand{\dynkindot}[2]{\fill (\dynkinstep*#1,\dynkinstep*#2) circle (\dynkinradius);}

\newcommand{\dynkinline}[4]{\draw[thin] (\dynkinstep*#1,\dynkinstep*#2) -- (\dynkinstep*#3,\dynkinstep*#4);}
\newcommand{\dynkindots}[4]{\draw[dotted] (\dynkinstep*#1,\dynkinstep*#2) -- (\dynkinstep*#3,\dynkinstep*#4);}
\newcommand{\dynkindoubleline}[4]{\draw[double,postaction={decorate}] (\dynkinstep*#1,\dynkinstep*#2) -- (\dynkinstep*#3,\dynkinstep*#4);}

\newenvironment{dynkin}{\begin{tikzpicture}[decoration={markings,mark=at position 0.7 with {\arrow{>}}}]}
{\end{tikzpicture}}

\DeclareSymbolFont{euletters}{U}{eur}{m}{n}
\SetSymbolFont{euletters}{bold}{U}{eur}{b}{n}
\DeclareMathSymbol{\varp}{\mathalpha}{euletters}{"7D}

\newcommand{\defeq}{\overset{\mathrm{def}}{=\!=}}
\newcommand{\amsbb}[1]{\mathbb{#1}}

\makeatletter
\newcommand*{\shifttext}[2]{%
  \settowidth{\@tempdima}{#2}%
  \makebox[\@tempdima]{\hspace*{#1}#2}%
}
\makeatother

\newcommand{\tensor}[1]{\underset{#1}{\otimes}}
\newcommand{\lsds}{{\mathrlap{\,\raisebox{0.8\depth}{\scaleobj{0.8}{\,\,$+$}}}\subset}}

\newcommand{\uniondot}{{\mathrlap{\raisebox{0.4pt}{\shifttext{5pt}{$\cdot$}}}{\cup}}}

\def\suchthat#1{%
\,\pmb{\left|\vphantom{#1}\right.}\,%
#1
}

\renewcommand{\qedsymbol}{{\hfill\raggedright \bf Q.E.D.}}
\newcommand{\qedcl}{{\hfill\raggedright  $\square\phantom{abc}$}}
\renewcommand{\qed}{\hfill \mbox{\raggedright \rule{.07in}{.1in}}}
\newcommand{\qedlem}{{\hfill\raggedright  $\#$}}

\makeatletter
\newenvironment{pf}[1][\proofname:]{\par
  \pushQED{\qedsymbol}%
  \normalfont \topsep6\p@\@plus6\p@\relax
  \list{}{\leftmargin=10pt
          \rightmargin=0pt
          \settowidth{\itemindent}{\bfseries#1}%
          \labelwidth=\itemindent
          \parsep=0pt \listparindent=\parindent 
          }

  \item[\hskip\labelsep
        \bfseries
    #1\@addpunct{:}]\ignorespaces
}{%
  \popQED\endlist\@endpefalse
}
\makeatother













\DeclareMathOperator{\Ext}{Ext}
\DeclareMathOperator{\Rad}{Rad}

\newcommand{\bggO}{\pazocal{O}}
\newcommand{\bbggO}{{\bar{\pazocal{O}}}}

\def\suchthat#1{%
\,\pmb{\left|\vphantom{#1}\right.}\,%
#1
}

\newtheoremstyle{mythm}
	{\topsep}
	{\topsep}
	{\upshape}
	{0pt}
	{\bfseries}
	{.}
	{5pt plus 1pt minus 1pt}
	{}
\theoremstyle{mythm}

\newtheorem{thrm}{Theorem}[section]

\usepackage{framed}

\colorlet{shadecolor}{blue!10}

\newenvironment{thm}
  {\colorlet{shadecolor}{orange!10}\begin{shaded}\begin{thrm}}
  {\end{thrm}\end{shaded}}

\newtheorem{propo}[thrm]{Proposition}

\newenvironment{prop}
  {\colorlet{shadecolor}{orange!15}\begin{shaded}\begin{propo}}
  {\end{propo}\end{shaded}}
  
\newtheorem{definit}[thrm]{Definition}

\newenvironment{define}
  {\colorlet{shadecolor}{yellow!15}\begin{shaded}\begin{definit}}
  {\end{definit}\end{shaded}}
  
\newtheorem{corol}[thrm]{Corollary}

\newenvironment{cor}
  {\colorlet{shadecolor}{orange!15}\begin{shaded}\begin{corol}}
  {\end{corol}\end{shaded}}

\newtheorem{lma}[thrm]{Lemma}

\newenvironment{lem}
  {\colorlet{shadecolor}{orange!15}\begin{shaded}\begin{lma}}
  {\end{lma}\end{shaded}}
  
\newtheorem{conjec}[thrm]{Conjecture}
\newtheorem{openquestion}[thrm]{Open Question}

\newenvironment{conj}
  {\colorlet{shadecolor}{red!15}\begin{shaded}\begin{conjec}}
  {\end{conjec}\end{shaded}}
  
  \newenvironment{openq}
  {\colorlet{shadecolor}{red!15}\begin{shaded}\begin{openquestion}}
  {\end{openquestion}\end{shaded}}

\newtheoremstyle{remark}
	{\topsep}
	{\topsep}
	{\upshape}
	{10pt}
	{\bfseries\itshape}
	{.}
	{5pt plus 1pt minus 1pt}
	{}
\theoremstyle{remark}

\everymath{\displaystyle}
\newtheorem{rmrk}[thrm]{Remark}

\newenvironment{rmk}
  {\colorlet{shadecolor}{gray!10}\begin{shaded}\begin{rmrk}}
  {\end{rmrk}\end{shaded}}

\newtheorem{exmp}[thrm]{Example}

\newenvironment{exm}
  {\colorlet{shadecolor}{blue!10}\begin{shaded}\begin{exmp}}
  {\end{exmp}\end{shaded}}

\usepackage{tikz,tikz-cd}
\usetikzlibrary{matrix,arrows,calc,decorations.pathreplacing,decorations.markings,intersections,shapes.geometric,through,fit,shapes.symbols,positioning,decorations.pathmorphing,snakes}

\makeatletter
\newcommand*\longcong[1]{%
    \setbox0=\hbox{\scriptsize#1}%
    \setbox2=\hbox{$\m@th{\simeq}$}%
    \stackrel{\copy0}{%
        \ifdim\wd2<\wd0
            \resizebox{\wd0}{\ht2}{$\m@th{\simeq}$}%
        \else
            \simeq
        \fi
    }%
}

\makeatletter
\newcommand*\longequal[1]{%
    \setbox0=\hbox{\scriptsize#1}%
    \setbox2=\hbox{$\m@th{=}$}%
    \stackrel{\copy0}{%
        \ifdim\wd2<\wd0
            \resizebox{\wd0}{\ht2}{$\m@th{=}$}%
        \else
            =
        \fi
    }%
}

\tikzset{cong/.style={draw=none,edge node={node [sloped, allow upside down, auto=false]{$\simeq$}}},
        Isom/.style={draw=none,every to/.append style={edge node={node [sloped, allow upside down, auto=false]{$\longcong{\phantom{abcdef}}$}}}},
       Equal/.style={draw=none,every to/.append style={edge node={node [sloped, allow upside down, auto=false]{$\longequal{\phantom{abcdef}}$}}}}
       }

\numberwithin{equation}{section}

\title{$\phantom{a}$\\[-1.7in]
\href{http://www.jacobs-university.de}{\centering\includegraphics[width=9cm]{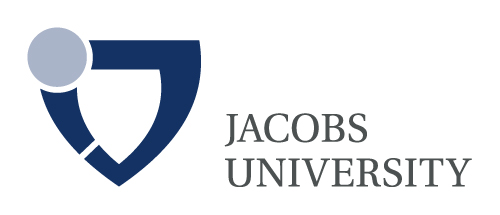}}
\\[0.2in]
{\bf Doctoral Thesis, Draft 04} \\[0.1in]
\bf\textit{Categories $\bggO$ for Dynkin Borel Subalgebras} \\[0.1in]
\bf\textit{of Root-Reductive Lie Algebras}\\[0.1in]
{\rm \normalsize by}\\[0.1in]}

\author{ {\bf Thanasin V. Nampaisarn}  \\[0.1in]
Department of Mathematics and Logistics \\[0.1in]
Focus Area Mobility \\[0.1in]
Jacobs University Bremen\\[0.1in]
{\small \href{mailto:t.nampaisarn@jacobs-university.de}{t.nampaisarn@jacobs-university.de}}\\[0.1in]
{\normalsize A Thesis submitted in partial fulfilment}\\[0.1in]
{\normalsize of the requirements for the degree of}\\[0.1in]
{\bf Doctor of Philosophy}\\[0.1in]
{\bf in Mathematics}
}
\date{\today}

\begin{document}
\pdfbookmark[0]{Front}{Front}
\clearpage

\pdfbookmark[1]{Cover}{Cover}
\begin{titlepage}
        \begin{flushright}\href{http://www.jacobs-university.de}{\centering\includegraphics[width=6cm]{Jacobs_University_Logo.jpg}}\end{flushright}
    
    \begin{center}
        \vspace{0.2in}
        {\bf\textit{\Large Categories $\bggO$ for Dynkin Borel Subalgebras}} \\[0.1in]
{\bf\textit{\Large of Root-Reductive Lie Algebras}}\\[0.2in]
{\rm \normalsize by}\\[0.2in]
        {\bf \href{mailto:t.nampaisarn@jacobs-university.de}{Thanasin V. Nampaisarn}}  \\[0.1in]
Department of Mathematics and Logistics \\[0.1in]
Focus Area Mobility \\[0.1in]
Jacobs University Bremen\\[0.1in]
        \vspace{0.2in}
        {\normalsize A Thesis Submitted in Partial Fulfilment}\\[0.1in]
{\normalsize of the Requirements for the Degree of}\\[0.2in]
{\bf Doctor of Philosophy}\\[0.1in]
{\bf in Mathematics}\\[0.3in]
    \end{center}
    
    \begin{flushleft}
   {\bf Approved, Thesis Committee}\\[0.2in]
    
    \href{mailto:i.penkov@jacobs-university.de}{Prof. Dr. Ivan Penkov}\\[-0.1in]
    
    \noindent\rule{4in}{0.4pt}
    
    Chair, Jacobs University Bremen\\[0.2in]
    
     \href{mailto:a.huckleberry@jacobs-university.de}{Prof. Dr. Alan Huckleberry}\\[-0.1in]
    
    \noindent\rule{4in}{0.4pt}
    
    Jacobs University Bremen / Ruhr-Universit\"at Bochum\\[0.2in]
    
     \href{mailto:serganov@math.berkeley.edu}{Prof. Dr. Vera Serganova}\\[-0.1in]
    
    \noindent\rule{4in}{0.4pt}
    
    University of California, Berkeley\\[0.3in]
    
    \end{flushleft}
    \vfill
    \begin{flushright}
    	Date of Defense: {\it Tuesday, June 06, 2017}
    \end{flushright}
\end{titlepage}
\thispagestyle{empty}

\pagenumbering{roman}

\pagebreak

$\phantom{a}$
\pagebreak

\pdfbookmark[1]{Abstract}{Abstract}
\begin{abstract}
\thispagestyle{plain}
\hspace{1cm} The purpose of my Ph.D. research is to define and study an analogue of the classical Bernstein-Gelfand-Gelfand (BGG) category $\bggO$ for the Lie algebra $\gfrak$, where $\gfrak$ is one of the finitary, infinite-dimensional Lie algebras $\gllie_\infty(\amsbb{K})$, $\sllie_\infty(\amsbb{K})$, $\solie_\infty(\amsbb{K})$, and $\splie_\infty(\amsbb{K})$.  Here, $\amsbb{K}$ is an algebraically closed field of characteristic $0$.  We call these categories \emph{extended categories $\bggO$} and use the notation $\bbggO$.  While the categories $\bbggO$ are defined for all splitting Borel subalgebras of $\gfrak$, this research focuses on the categories $\bbggO$ for very special Borel subalgebras of $\gfrak$ which we call Dynkin Borel subalgebras.  Some results concerning block decomposition and Kazhdan-Lusztig multiplicities carry over from usual categories $\bggO$ to our categories $\bbggO$.  There are differences which we shall explore in detail, such as the lack of some injective hulls.  In this connection, we study truncated categories $\bbggO$ and are able to establish an analogue of BGG reciprocity in the categories $\bbggO$.

\end{abstract}
\pagebreak

\setcounter{page}{2}
$\phantom{A}$

\pagebreak
\pdfbookmark[1]{Acknowledgement}{Acknowledgement}

\begin{center}{\bf Acknowledgement}\end{center}

This Ph.D. research is under the supervision of \href{http://math.jacobs-university.de/penkov/}{\bf Prof. Ivan Penkov}, Department of Mathematics and Logistics, Jacobs University Bremen.  I am greatly indebted to Prof. Penkov for introducing me to representation theory and to this research topic.  His guidance has been tremendous help towards the fulfillment of my study at Jacobs University Bremen.

The research topic has been partially suggested by \href{http://math.berkeley.edu/~serganov/}{\bf Prof. Vera Serganova} from the University of California, Berkeley.  She has given numerous significant insights to this dissertation.   \href{http://users.math.yale.edu/public_html/People/gjz2.html}{\bf Prof. Gregg J. Zuckerman} from Yale University has also provided many valuable suggestions.  I am very grateful for the assistance from both professors.   Additionally, \href{http://www2.math.uu.se/~mazor/}{\bf Prof. Volodymyr Mazorchuk} from Uppsala University has assisted me with finding some references.

I wish to acknowledge the DFG (Deutsche Forschungsgemeninschaft) for the financial support during my three-year doctoral study, without which this work would not have been possible.  I would also like to thank the whole Department of Mathematics and Logistics at Jacobs University Bremen for the opportunity to have nice discussions on various topics in mathematics.  In the end, I would like to emphasize that the constant support from my friends and my family, especially from my father, has been extremely crucial to the completion of my doctoral study.

\pagebreak

$\phantom{A}$

\pagebreak

\pdfbookmark[1]{Declaration}{Declaration}

\begin{center}
{\bf Declaration}

\end{center}

I hereby declare that this submission is my own work and that, to the best of my knowledge and belief, it contains no material previously published or written by another person nor material which to a substantial extent has been accepted for the award of any other degree or diploma of the university or other institute of higher learning, except where due acknowledgment has been made in the text.\\[0.5in]

\noindent Signature: \rule{1.5in}{0.4pt}$\phantom{aaa}$Name: Thanasin V. Nampaisarn$\phantom{aaa}$Date: \today

\pagebreak

\phantom{a}
\pagebreak

\bigskip

\pdfbookmark[1]{Table of Contents}{Table of Contents}
\tableofcontents

\pagebreak

\cleardoublepage
\pagenumbering{arabic}

\pdfbookmark[0]{Contents}{Contents}

\section{Introduction}

	 The classical Bernstein-Gelfand-Gelfand (BGG) category $\pazocal{O}$ for a finite-dimensional reductive Lie algebra $\gfrak$ has provided a rich ground for research topics for a few decades.   Over an algebraically closed field of characteristic $0$, the category $\bggO$ is well understood.  Amongst the well known results is the Kazhdan-Lusztig Conjecture, which is a connection between the category $\bggO$ and the Hecke algebra.
	 
	 Extensive studies of analogues of categories $\bggO$ have also been undertaken for affine Kac-Moody Lie algebras~\cite{FF,SVV}.  In a parallel fashion to the finite-dimensional theory, an affine Kac-Moody Lie algebra $\hat{\gfrak}$ has a Borel subalgebra $\hat{\bfrak}$ containing a Cartan subalgebra $\hat{\hfrak}$.  It is natural to define the category $\bggO$ for the Lie algebra $\hat\gfrak$ by requiring that $\hat{\hfrak}$ act semisimply and $\hat{\bfrak}$ act locally finitely on each module in $\bggO$.  
	 
	 In this work, we generalize the definition of the classical BGG categories $\bggO$ to a class of locally finite Lie algebras called ``root-reductive Lie algebras'' over an algebraically closed field of characteristic $0$.  Splitting maximal toral subalgebras play an essential role in this dissertation.  
	 
	 Cartan subalgebras and maximal toral subalgebras of a root-reductive Lie algebras have been studied in~\cite{DPS} and~\cite{NP}.  In addition, \cite{DPS} gives a rough classification of root-reductive Lie algebras as follows: if $\gfrak$ is a root-reductive Lie algebra, then $[\gfrak,\gfrak]$ is a direct sum of finite-dimensional semisimple Lie algebras and simple finitary Lie algebras $\sllie_\infty$, $\solie_\infty$, and $\splie_\infty$, each with at most countable multiplicity.
	 
	 To define a category $\bggO$ for a root-reductive Lie algebra $\gfrak$, we further need to understand the structure of Borel subalgebras of $\gfrak$.  Borel subalgebras for the simple finitary Lie algebras $\sllie_\infty$, $\solie_\infty$, and $\splie_\infty$ and for root-reductive Lie algebras are very well understood (see~\cite{BorelDC} and~\cite{DP1999}).  For a given Borel subalgebra $\bfrak$ of $\gfrak$ containing a splitting maximal toral subalgebra $\hfrak$, we define an extended category $\bggO$, denoted by $\bbggO$, for $\gfrak$ with respect to $\bfrak$ by demanding that the objects in the category be $\hfrak$-semisimple with locally finite $\bfrak$-action and with finite-dimensional $\hfrak$-weight spaces.
	 
	The condition that the objects in $\bbggO$ have finite-dimensional $\hfrak$-weight spaces allows us to decompose every module in $\bbggO$ into indecomposable direct summands.  However, without extra conditions on the Borel subalgebra $\bfrak$, the category $\bbggO$ lacks some important $\gfrak$-modules---the Verma modules.  To remedy this situation, we find a necessary and sufficient condition for the category $\bbggO$ to contain all Verma modules, namely, the Borel subalgebra $\bfrak$ must be a Dynkin Borel subalgebra.
	
	With this additional restriction on $\bfrak$, there are several important consequences.  First, every module in $\bbggO$ has an analogue of composition series, and therefore, composition factors and multiplicities are well defined.  Second, the block decomposition of $\bbggO$ resembles the block decomposition of a classical category $\bggO$.  In addition, Kazhdan-Lusztig theory may be applied to $\bbggO$, leading to a Kazhdan-Lusztig multiplicity formula.
	
	Nonetheless, $\bbggO$ also greatly differs from a classical category $\bggO$.  In the finite-dimensional case, $\bggO$ is a highest-weight category.  Therefore, BGG reciprocity comes naturally to a classical category $\bggO$.  On the other hand, $\bbggO$ does not have enough injectives, and therefore, is not a highest-weight category.  In an attempt to remedy the situation, a truncation method is applied to $\bbggO$.  This truncation technique yields several full subcategories of $\bggO$, which are called truncated subcategories.  While not all of these subcategories are highest-weight categories, they satisfy a version of BGG reciprocity.
	
	In this dissertation, Section~\ref{ch:intro} provides definitions of crucial terms such as root-reductive Lie algebras and Dynkin Borel subalgebras.  Explicit examples of root-reductive Lie algebras along with some of their important subalgebras are given in Section~\ref{ch:exm}.  Section~\ref{ch:keyresults} covers the finite-dimensional theory that is needed in subsequent sections.  The category $\bbggO$ for a root-reductive Lie algebra $\gfrak$ is defined in Section~\ref{ch:bbggo}, where we shall also prove key results such as direct sum decompositions and generalized composition series of objects in $\bbggO$. 
	
	Due to the importance of Verma modules in our study of the category $\bbggO$, they are discussed in detail in Section~\ref{ch:verma}.  In Section~\ref{sec:KLTheory}, the theory of Coxeter groups is generalized to the case of infinite generating sets.  With this generalization, we obtain results for the corresponding Hecke algebras and Kazhdan-Lusztig polynomials which are similar to the case of Coxeter groups with finite generating sets.  Section~\ref{ch:blocks} concerns the block structure of $\bbggO$, which is connected to the Kazhdan-Lusztig formula for the multiplicity of a composition factor of a Verma module.  Finally, the last section (Section~\ref{ch:trunc}) studies BGG reciprocity in $\bbggO$ by means of truncation.  

\pagebreak

\section{Preliminaries}
\label{ch:intro}

The base field is $\amsbb{K}$, which is assumed to be algebraically closed and of characteristic $0$.  All Lie algebras and vector spaces are defined over $\amsbb{K}$.   For a vector space $V$, $V^*$ denotes its algebraic dual space $\Hom_\amsbb{K}(V,\amsbb{K})$.

\subsection{Root-Reductive Lie Algebras}
\label{ch:rootreductive}

\begin{define}
	A Lie algebra $\mathfrak{g}$ is said to be \emph{locally finite}\index{local finiteness} if every finite subset of $\mathfrak{g}$ generates a finite-dimensional Lie subalgebra of $\mathfrak{g}$.\index{Lie subalgebra}\index{subalgebra of $\gfrak$|see {Lie subalgebra}}  A locally finite Lie algebra $\mathfrak{g}$ is \emph{locally solvable}\index{local solvability} if every finite-dimensional Lie subalgebra of $\mathfrak{g}$ is solvable.  Similarly, a locally finite Lie algebra $\mathfrak{g}$ is \emph{locally nilpotent}\index{local nilpotency} if every finite-dimensional Lie subalgebra of $\mathfrak{g}$ is nilpotent.\index{subalgebra of $\gfrak$|see {Lie subalgebra}}  
\end{define}

\begin{define}
	A Lie algebra $\mathfrak{g}$ is \emph{locally reductive} if it is the union $\bigcup_{n\in\amsbb{Z}_{> 0}}\,\gfrak_n$ of nested finite-dimensional reductive Lie subalgebras $\left(\gfrak_n\right)_{n\in\amsbb{Z}_{> 0}}$ such that $\gfrak_n$ is reductive in $\gfrak_{n+1}$.  We say that $\gfrak$ is \emph{locally semisimple} if each $\gfrak_n$ is semisimple.  Moreover, $\gfrak$ is \emph{locally simple} if each $\gfrak_n$ is simple.
\end{define}

\begin{define}
	An inclusion of finite-dimensional reductive Lie algebras $\mathfrak{l}\hookrightarrow\mathfrak{m}$ is a \emph{root inclusion} if, for some Cartan subalgebra $\mathfrak{c}$ of $\mathfrak{m}$, the subalgebra $\mathfrak{l}\cap \mathfrak{c}$ is a Cartan subalgebra of $\mathfrak{l}$ and each $(\mathfrak{l}\cap\mathfrak{c})$-root space  of $\mathfrak{l}$ is a root space of $\mathfrak{m}$.
	
	A \emph{root-reductive} Lie algebra $\gfrak$ is a locally reductive Lie algebra $\gfrak=\bigcup_{n\in\amsbb{Z}_{> 0}}\,\gfrak_n$, where $\left(\gfrak_n\right)_{n\in\amsbb{Z}_{>0}}$ is a nested system of finite-dimensional reductive subalgebras, with the property that there exist nested subalgebras $\kfrak_1\subseteq \kfrak_2\subseteq \ldots$, where $\kfrak_n\subseteq \gfrak_n$ is a Cartan subalgebra of $\gfrak_n$, such that each inclusion $\gfrak_n \hookrightarrow \gfrak_{n+1}$ is a root inclusion with respect to the Cartan subalgebra $\kfrak_{n+1}$ of $\gfrak_{n+1}$.
\end{define}

\begin{define}
	Let $\gfrak$ be a locally reductive Lie algebra.  A Lie subalgebra $\hfrak$ of $\gfrak$ is said to be a \emph{splitting maximal toral subalgebra}\index{Lie subalgebra!splitting maximal toral subalgebra} if there exists a directed system $\left(\gfrak_n\right)_{n\in \amsbb{Z}_{>0}}$ of finite-dimensional reductive Lie subalgebras of $\gfrak$ for which $\lim_{\underset{n}{\boldsymbol{\longrightarrow}}}\, \gfrak_n=\gfrak$, $\hfrak\cap\gfrak_n$ is a maximal toral subalgebra of $\gfrak_n$ for each $n \in \amsbb{Z}_{>0}$\nomenclature[B]{$\amsbb{Z}_{> 0}$}{the set of natural numbers, i.e., positive integers}\index{Lie subalgebra!Borel subalgebra}\nomenclature[B]{$\gfrak^\lambda$}{$\big\{x\in\mathfrak{g}\suchthat{[h,x]=\lambda(h)x}\big\}$ for $\lambda\in\hfrak^*$}, and $\gfrak$ has the following \emph{$\hfrak$-root space decomposition}\index{root space decomposition}:
		\begin{align}
		\mathfrak{g}=\bigoplus\limits_{\alpha\in\mathfrak{h}^*}\mathfrak{g}^\alpha=\mathfrak{h}\,\oplus\,\bigoplus\limits_{\alpha \in \Delta}\mathfrak{g}^\alpha\,,
		\end{align}
		where $\mathfrak{g}^\alpha$ is the eigenspace\index{root space} $\big\{x\in\mathfrak{g}\suchthat{[h,x]=\alpha(h)x}\text{ for all }h\in\hfrak\big\}$ for every $\alpha\in\mathfrak{h}^*$, and $\Delta$ is the set of \emph{roots}\index{root}, i.e.,  the nonzero linear functionals $\alpha \in \mathfrak{h}^*$ for which $\mathfrak{g}^\alpha \neq 0$\nomenclature[B]{$0$}{the additive identity of a ring or a ring module (sometimes, $0$ also denotes the set $\{0\}$)}.  For $\alpha \in \Delta$, $\gfrak^\alpha$ is known as the \emph{root space} associated to $\alpha$.\nomenclature[B]{$V^*$}{the dual vector space of a vector space $V$} \nomenclature[B]{$\oplus$}{direct sum} \nomenclature[B]{$\bigoplus$}{direct sum over an index}\nomenclature[B]{$[\cdot,\cdot]$}{the Lie bracket of a Lie algebra}\nomenclature[B]{$\hfrak$}{a splitting maximal toral subalgebra of $\gfrak$}  Note that $\gfrak^\alpha$ is always a one-dimensional vector space.
\end{define}

\begin{define}
	Let $\gfrak$ be a locally reductive Lie algebra.  A Lie subalgebra $\bfrak$ of $\gfrak$ is said to be a \emph{splitting Borel subalgebra}\index{Lie subalgebra!splitting Borel subalgebra} if $\bfrak$ is a maximal locally solvable subalgebra of $\gfrak$ containing a splitting maximal toral subalgebra of $\gfrak$.\nomenclature[B]{$\bfrak$}{a splitting Borel subalgebra of $\gfrak$}  
\end{define}
	
	In~\cite{dancohen},~\cite{BorelDC},~\cite{DPS} and~\cite{DP1999}, locally reductive Lie algebras and root-reductive Lie algebras are studied.  In the case of root-reductive Lie algebra, a (splitting) Borel subalgebra $\bfrak$ containing a splitting maximal toral subalgebra $\hfrak$ arises from a choice of a set of \emph{$\bfrak$-positive roots} $\Delta^+\subseteq \Delta$ subject to the requirements that both subsets $\Delta^+$ and $\Delta\smallsetminus\Delta^+$ are additive and that $\alpha\in\Delta^+$ if and only if $-\alpha \in \Delta\smallsetminus\Delta^+$.  We set $\Delta^-:=-\Delta^+=\Delta\smallsetminus\Delta^+$ and call $\Delta^-$ the set of \emph{$\bfrak$-negative roots}.  Then $\Delta$ is the disjoint union $\Delta^+\uniondot \Delta^-$, the locally nilpotent subalgebra $\nfrak=\nfrak^+:=[\bfrak,\bfrak]$ is the direct sum $\bigoplus_{\alpha \in \Delta^+}\,\gfrak^\alpha$, and the Borel subalgebra $\bfrak$ is given by $\bfrak=\bfrak^+=\hfrak\oplus\nfrak^+$ (this is a direct sum of vector spaces, not of Lie algebras).  The Lie algebra $\gfrak$ has the following \emph{triangular decomposition}
	\begin{align}
		\gfrak=\nfrak^-\oplus\hfrak\oplus\nfrak^+\,,
	\end{align}
	where $\nfrak^-$ is the opposite subalgebra to $\nfrak^+$, namely, $\nfrak^-=\bigoplus_{\alpha\in\Delta^-}\,\gfrak^\alpha$.  The Lie algebra $\hfrak\oplus\nfrak^-$ is denoted by $\bfrak^-$.  It is opposite to $\bfrak^+$ in the sense that $\bfrak^++\bfrak^-=\gfrak$ and $\bfrak^+\cap\bfrak^-=\hfrak$.  By the \emph{Poincar\'e-Birkhoff-Witt (PBW) Theorem}\nomenclature[A]{PBW}{Poinca\'e-Birkhoff-Witt}\index{Poincar\'e-Birkhoff-Witt Theorem}, we see that $\Ulie(\mathfrak{g})=\Ulie\left(\mathfrak{n}^-\right)\boldsymbol{\cdot} \Ulie\left(\mathfrak{h}\right)\boldsymbol{\cdot}\Ulie\left(\mathfrak{n}^+\right)$.  Here, $\Ulie(L)$ is the universal enveloping algebra\index{universal enveloping algebra} of a Lie algebra $L$.\nomenclature[B]{$\Ulie(L)$}{the universal enveloping algebra of a Lie algebra $L$} \nomenclature[B]{$\boldsymbol{\cdot}$}{the multiplication in a ring such as the universal enveloping algebras}
	
	Let the \emph{root lattice}\index{root lattice} $\Lambda$\nomenclature[B]{$\Lambda$}{the root lattice} be the $\amsbb{Z}$-span in $\hfrak^*$ of $\Delta$, and $\Lambda^+$ be the $\amsbb{Z}_{\geq 0}$-span in $\hfrak^*$ of $\Delta^+$.\nomenclature[B]{$\Lambda^-$}{the $\amsbb{Z}_{\geq 0}$-span of $\Delta^-$}  We equip $\mathfrak{h}^*$ with a partial order $\preceq$ given by\nomenclature[B]{$\amsbb{Z}$}{the set of integers} \nomenclature[B]{$\amsbb{Z}_{\geq 0}$}{the set of nonnegative integers} \nomenclature[B]{$\preceq$}{a weak partial order on $\hfrak^*$}\index{order!weak partial order} 
		\begin{align}
			\lambda \preceq \mu \text{ iff } \mu- \lambda \in \Lambda^+
		\end{align}  		
	 for $\lambda,\mu\in\mathfrak{h}^*$.  We also write $\Lambda^-:=-\Lambda^+$ for the $\amsbb{Z}_{\geq 0}$-span of $\Delta^-$.

	\begin{define}
		An element $\alpha \in \Delta^+$ is said to be a \emph{simple $\bfrak$-positive root}, or a \emph{simple root} with respect to $\bfrak$, if it cannot be decomposed as a (finite) sum of two or more $\bfrak$-positive roots.  We usually use the symbol $\Sigma^+$ or $\Sigma$ for the set of simple $\bfrak$-positive roots.\index{root!simple negative root}\nomenclature[B]{$\Sigma^-$}{the set of simple negative roots}
		
		Similarly, we say that $\alpha \in \Delta^-$ is a \emph{simple negative root} with respect to $\bfrak$ if it cannot be decomposed as a sum of two or more $\bfrak$-negative roots.  The symbol $\Sigma^-$ denotes the set of simple $\bfrak$-negative roots.\index{root!simple root}\index{root!simple positive root}\nomenclature[B]{$\Sigma=\Sigma^+$}{the set of simple (positive) roots}  Clearly, $\Sigma^-=-\Sigma^+$.
	\end{define}
	
	From now on, $\gfrak$ is a root-reductive Lie algebra with a splitting Borel subalgebra $\bfrak$ containing a splitting maximal toral subalgebra $\hfrak$.  We assume that $\gfrak$ is the union $\bigcup_{n\in\amsbb{Z}_{>0}}\,\gfrak_n$ of nested finite-dimensional reductive Lie subalgebras $\gfrak_n$ for which $\hfrak_n:=\hfrak\cap \gfrak_n$ is a maximal toral subalgebra of $\gfrak_n$, $\bfrak_n=\bfrak^+_n:=\bfrak\cap \gfrak_n$ is a Borel subalgebra of $\gfrak_n$, and $\nfrak_n=\nfrak^+_n:=\nfrak^+\cap\gfrak_n=\left[\bfrak^+_n,\bfrak^+_n\right]$ is a nilpotent subalgebra of $\gfrak_n$.   We also write $\bfrak_n^-:=\bfrak^-\cap\gfrak_n$ and $\nfrak^-_n:=\nfrak^-\cap\gfrak_n$.  In the case where $\gfrak$ is locally semisimple, we also assume that each $\gfrak_n$ is semisimple.  If $\gfrak$ is locally simple, each $\gfrak_n$ is taken to be simple.
	
	For each $n\in\amsbb{Z}_{>0}$, $W_n$ denotes the Weyl group of $\gfrak_n$.  Since the embedding $\gfrak_n\hookrightarrow \gfrak_{n+1}$ is a root inclusion, this induces an embedding $W_n\hookrightarrow W_{n+1}$.  The Weyl group $W$ of $\gfrak$ is simply the direct limit $\lim\limits_{\underset{n}{\boldsymbol{\longrightarrow}}}\,W_n$.
	
	\begin{thm}
	Let $\gfrak$ be a root-reductive Lie algebra with a splitting maximal toral subalgebra $\hfrak$.  Then, the image of $\hfrak$ under the quotient map $\gfrak\to\gfrak/[\gfrak,\gfrak]$ is the whole $\gfrak/[\gfrak,\gfrak]$.  
	In addition, $\hfrak\cap[\gfrak,\gfrak]$ is a splitting maximal toral subalgebra of $[\gfrak,\gfrak]$.
	\label{thm:splittingmaxtoral}
	
	\begin{pf}
		Suppose that $\hfrak$ does not surject onto $\gfrak/[\gfrak,\gfrak]$ under the canonical projection.  As $\hfrak$ is splitting, there exists an element $x\in\gfrak\smallsetminus\big(\hfrak+[\gfrak,\gfrak]\big)$ such that $x$ lies in a root space $\gfrak^\alpha$ for some $\alpha\in \Delta$.
		
		Since $\alpha\neq 0$, there exists $h\in\hfrak$ for which $\alpha(h)\neq 0$.  Then, $[h,x]=\alpha(h)\,x$ implies that $x\in [\gfrak,\gfrak]$.  This is a contradiction.  Therefore, $\hfrak$ indeed surjects onto $\gfrak/[\gfrak,\gfrak]$.
		
		
		To show that $\bar\hfrak:=\hfrak\cap[\gfrak,\gfrak]$ is a splitting maximal toral subalgebra of $\bar\gfrak:=[\gfrak,\gfrak]$, we first note that, for each positive integer $n$, 
		\begin{align}
			\left[\gfrak_n,\gfrak_n\right]=\bar\gfrak\cap\gfrak_n=:\bar\gfrak_n\,.
		\end{align}  
		As $\gfrak_n$ is a finite-dimensional reductive Lie algebra, 
		\begin{align}
			\gfrak_n=\bar\gfrak_n\oplus \mathfrak{z}_n\,,
		\end{align} where $\mathfrak{z}_n$ is the center of the Lie algebra $\gfrak_n$.  As $\hfrak_n=\hfrak\cap\gfrak_n$ is a maximal toral subalgebra of $\gfrak_n$, we conclude that 
		\begin{align}
			\hfrak_n=\bar\hfrak_n\oplus\mathfrak{z}_n
		\end{align} for some maximal toral subalgebra $\bar\hfrak_n$ of $\bar\gfrak_n$.  Ergo, $\bar\hfrak_n=\hfrak_n\cap\bar\gfrak_n$.  
		
		The image of the inclusion $\gfrak_n\hookrightarrow\gfrak_{n+1}$ restricted onto $\bar\hfrak_n$ lies in 
		\begin{align}\bar\gfrak_{n+1}\cap \hfrak_{n+1}=\bar\hfrak_{n+1}\end{align} as $\bar\hfrak_n\subseteq \bar\gfrak_n=\left[\gfrak_n,\gfrak_n\right]$ and $\bar\hfrak_n\subseteq\hfrak_n$.  Thus, the inclusion $\gfrak_n\hookrightarrow \gfrak_{n+1}$ induces an embedding $\bar\hfrak_n\hookrightarrow \bar\hfrak_{n+1}$.  This means $\bar\hfrak$ is the direct limit of the maximal toral subalgebras $\bar\hfrak_n\subseteq\bar\gfrak_n$, making $\bar\hfrak$ a maximal toral subalgebra of $\bar\gfrak$.  Since the inclusion $\gfrak_n\hookrightarrow\gfrak_{n+1}$ is a root inclusion with respect to $\hfrak_{n+1}$, the induced inclusion $\bar\gfrak_n\hookrightarrow\bar\gfrak_{n+1}$ is also a root inclusion with respect to $\bar\hfrak_{n+1}$.  Now, $\bar\gfrak_n$ has a $\bar\hfrak_n$-root decomposition, whence the $\bar\hfrak_n$-root spaces of $\bar\gfrak_n$ are $\bar\hfrak$-root spaces of $\bar\gfrak$.  Since $\bar\gfrak=\bigcup_{n\in\amsbb{Z}_{>0}}\,\bar\gfrak_n$, we deduce that $\bar\gfrak$ has a $\bar\hfrak$-root decomposition, and the proof is now complete.
	\end{pf}
\end{thm}

In general, a \emph{locally reductive Lie algebra} has a Jordan decomposition (see~\cite{DPS}).  Note that a root-reductive Lie algebra is also locally reductive.  Therefore,  the Jordan decomposition of $x\in\gfrak_n$ is the same as the Jordan decomposition of $x$ considered as an element of $\gfrak_{n+k}$ for every integer $k>0$.  Consequently, this defines a (unique) decomposition of $x$ as the sum of a semisimple part $x_\text{ss}$ and a nilpotent part $x_\text{nil}$.

\begin{define}
	For a root-reductive Lie algebra $\gfrak$, a subalgebra $\kfrak$ is a \emph{Cartan subalgebra} if $\kfrak$ is the centralizer of the subset $\kfrak_\text{ss}\defeq\big\{x_\text{ss}\,\boldsymbol{|}\,x\in\kfrak\big\}$.  A Cartan subalgebra $\kfrak\subseteq\gfrak$ is \emph{splitting} if $\gfrak$ is a $\kfrak_\text{ss}$-weight module (noting that $\kfrak_\text{ss}$ is in fact a toral subalgebra of $\gfrak$).
\end{define}

From Theorem~\ref{thm:splittingmaxtoral}, we obtain the corollary below.  In~\cite{DPS}, it is proven that, if $\kfrak$ is a Cartan subalgebra of $\gfrak$, then $\kfrak_{\text{ss}}$ is a maximal toral subalgebra.  This article also shows that, for the Lie algebras $\gllie_\infty$, $\sllie_\infty$, $\solie_\infty$, and $\splie_\infty$, a subalgebra is a splitting maximal toral algebra if and only if it is a splitting Cartan subalgebra.  

\begin{cor}
	If a subalgebra $\hfrak$ of a root-reductive Lie algebra $\gfrak$ is a splitting maximal toral subalgebra, then $\hfrak$ is a splitting Cartan subalgebra.
\end{cor}

In the rest of this dissertation, we work with splitting maximal toral subalgebras.
	
\subsection{Weight Modules}
\label{ch:wtmodules}

	 Let $M$ be a $\mathfrak{g}$-module.\index{module over $\gfrak$|see {Lie module}}  For each $\lambda \in \mathfrak{h}^*$, define $M^\lambda$ to be the subspace 
	 \begin{align}
	 	M^\lambda\defeq \big\{v \in M\suchthat{} h\cdot v=\lambda(h)\,v\text{ for all }h\in\hfrak\big\}\,.
	\end{align}
	If $M^\lambda \neq 0$, $\lambda$ is said to be a \emph{weight}\index{weight} of $M$, $M^\lambda$ is called a \emph{weight space}\index{weight space} of $M$, and elements of $M^\lambda$ are known as \emph{weight vectors}\index{vector!weight vector}.  If 
	\begin{align}
		M=\bigoplus\limits_{\lambda \in \mathfrak{h}^*}M^\lambda\,,
	\end{align}
	then $M$ is called an \emph{$\hfrak$-weight module}, or simply a \emph{weight module},\index{Lie module!weight module}\index{weight module|see {Lie module}} over $\mathfrak{g}$.  If $M$ is a weight $\gfrak$-module, then so are its submodules\index{Lie submodule}\index{submodule of a $\gfrak$-module|see {Lie submodule}} and quotients.  The proof of this fact is similar to Proposition 1.5 of \cite{kac}.  Note also that the direct sum of a collection of weight modules is also a weight module.
	
	A $\gfrak$-module $M$ is said to be a \emph{cyclic module}\index{cyclic module}\index{Lie module!cyclic module} over $\mathfrak{g}$ if it is generated as a $\Ulie(\mathfrak{g})$-module by a single nonzero element.   A cyclic $\gfrak$-module $M$ generated by $v \neq 0$ is said to be a \emph{highest-weight module}\index{Lie module!highest-weight module}\index{highest-weight module|see {Lie module}} with respect to the Borel subalgebra $\bfrak$ if
	\begin{align}
		\mathfrak{n}\cdot v=0
	\end{align}
	and there exists a weight $\lambda\in\mathfrak{h}^*$, known as the \emph{highest weight}\index{weight!highest weight} of $M$, such that 
	\begin{align}
		h\cdot v=\lambda(h)\,v
	\end{align} 
	for each $h \in \mathfrak{h}$.  The vector $v$ is then called a \emph{highest-weight vector} of $M$.\index{vector!highest-weight vector}  If $M$ is a highest-weight module over $\mathfrak{g}$ with highest weight $\lambda\in\mathfrak{h}^*$, then $M$ is a weight $\gfrak$-module.  Moreover, it follows immediately from the PBW Theorem that all possible weights of $M$ lie within $\lambda + \Lambda^-$. 	
	
	For a $\gfrak$-module $M$ (not necessarily a highest-weight module), a nonzero element $v\in M$ is called a \emph{singular vector} (with respect to $\bfrak$) if $\nfrak\cdot v=0$.  Furthermore, a weight $\lambda$ of $M$ is said to be a \emph{singular weight} if the weight space $M^\lambda$ contains a singular vector.  A subspace of $M^\lambda$ is a \emph{singular subspace} if it consists of singular vectors.  Note that a highest-weight vector of a highest-weight module is a singular vector, and the highest weight of this module is a singular weight.
	
	\begin{define}
	For every $\lambda\in\mathfrak{h}^*$, we define the \emph{Verma module}\index{Verma module|see {Lie module}}\index{Lie module!Verma module} over $\mathfrak{g}$ of highest weight $\lambda$ to be the left $\Ulie(\mathfrak{g})$-module
		\begin{align}
			\Mlie(\lambda;\mathfrak{g},\mathfrak{b},\mathfrak{h})\defeq \Ulie(\mathfrak{g})/I\,,
		\end{align}
		where $I$ is the left $\Ulie(\mathfrak{g})$-ideal generated by $\mathfrak{n}$ and $h-\lambda(h)\,1_{\Ulie(\gfrak)}$, for all $h \in \mathfrak{h}$.  If there is no confusion, we shall write $\Mlie(\lambda)$ for $\Mlie(\lambda;\mathfrak{g},\mathfrak{b,\mathfrak{h}})$.  \nomenclature[B]{$\Mlie(\lambda)=\Mlie(\lambda;\gfrak,\bfrak,\hfrak)$}{the Verma module over $\gfrak$ with $\bfrak$-highest weight $\lambda \in \hfrak^*$}
	\end{define}
	
	\begin{rmk}
	Note that $\Mlie(\lambda;\mathfrak{g},\mathfrak{b},\mathfrak{h})$ has a unique maximal proper $\gfrak$-submodule, which is denoted by $\Nlie(\lambda;\gfrak,\bfrak,\hfrak)$ (or by $\Nlie(\lambda)$ when there is no ambiguity).  The submodule $\Nlie(\lambda;\gfrak,\bfrak,\hfrak)$ is the sum of all proper submodules\index{Lie submodule!proper submodule} of $\Mlie(\lambda;\mathfrak{g},\mathfrak{b},\mathfrak{h})$.  Consequently,
	\begin{align}
		\Llie(\lambda;\mathfrak{g},\mathfrak{b},\mathfrak{h})\defeq\Mlie(\lambda;\mathfrak{g},\mathfrak{b},\mathfrak{h})/\Nlie(\lambda;\gfrak,\bfrak,\hfrak)\,,
	\end{align} 
	is the unique simple quotient of $\Mlie(\lambda;\mathfrak{g},\mathfrak{b,\mathfrak{h}})$ (also denoted by $\Llie(\lambda)$ if there is no confusion).\nomenclature[B]{$\Llie(\lambda)=\Llie(\lambda;\gfrak,\bfrak,\hfrak)$}{the simple quotient of $\Mlie(\lambda;\gfrak,\bfrak,\hfrak)$}  This has been pointed out in~\cite{penkov-dimitrov}.
	
	Let $v_\lambda$ be the image of $1_{\Ulie(\gfrak)}$ under the canonical projection $\Ulie(\mathfrak{g})\twoheadrightarrow\Ulie(\mathfrak{g})/I$.\nomenclature[B]{$X\twoheadrightarrow Y$}{an epimorphism from $X$ to $Y$}  Clearly, $v_\lambda$ is a highest-weight vector\index{vector!highest-weight vector} of $\Mlie(\lambda)$ with weight $\lambda$.  Thus, $\Mlie(\lambda)$, being generated by $v_\lambda$, is a highest-weight module over $\mathfrak{g}$ with highest weight $\lambda$.   It also follows that, for any highest-weight module $M$ over $\mathfrak{g}$ with highest weight $\lambda$, there is a unique epimorphism $\sigma:\Mlie(\lambda)\to M$ of $\Ulie(\mathfrak{g})$-modules, up to scalar multiples.  This is the \emph{universal property} of the Verma modules\index{universality property of Verma modules}.  Hence, $\Llie(\lambda)$ is (up to isomorphism) the unique simple highest-weight module over $\mathfrak{g}$ with highest weight $\lambda$. 
	
	Note also that $\Mlie(\lambda)$ is isomorphic to $\Ulie\left(\nfrak^-\right)$ as a left $\Ulie\left(\nfrak^-\right)$-module.  The homomorphism $\Ulie\left(\nfrak^-\right) \to \Mlie(\lambda)$ sending $1_{\Ulie\left(\nfrak^-\right)}$ to $v_\lambda$ is obviously an isomorphism of left $\Ulie\left(\nfrak^-\right)$-modules.
	\end{rmk}
	
\subsection{Dynkin Borel Subalgebras}	
\label{ch:dynkin}

	In this subsection, $\gfrak$ need not be finite-dimensional.  Furthermore, let $\rho_n$ denote the half sum of positive roots of $\gfrak_n$ with respect to $\bfrak_n$ (sometimes, the linear functionals $\rho_n$ are called the \emph{local half sum of positive roots}).
	
	\begin{define}
		We say that $\bfrak$ is a \emph{Dynkin Borel subalgebra}  of $\gfrak$ if it is generated by $\hfrak$ and the simple root spaces.  
	\end{define}
	
	\begin{define}
		A $\bfrak$-positive root $\alpha$ is \emph{of finite length}  (with respect to $\bfrak$) if there are only finitely many ways to express it as a sum of positive roots; otherwise, $\alpha$ is \emph{of infinite length}  (with respect to $\bfrak$).  A $\bfrak$-negative root $\alpha$ is said to be \emph{of finite length}  (with respect to $\bfrak$) if the positive root $-\alpha$ is of finite length; otherwise, $\alpha$ is \emph{of infinite length}  (with respect to $\bfrak$). 
	\end{define}
	
	It is an easy exercise to prove that $\bfrak$ is Dynkin if and only if every root is of finite length (with respect to $\bfrak$).  Furthermore, the partial order $\preceq$ on $\hfrak^*$ induced by $\bfrak$ is locally finite\footnote{A partially ordered set $(\pazocal{P},\leq)$ is said to be \emph{locally finite} if the set $\{x\in\pazocal{P}\,\boldsymbol{|}\,a\leq x\leq b\}$ is finite for every $a,b\in\pazocal{P}$.}  if and only if $\bfrak$ is Dynkin.  In other words, a Dynkin Borel subalgebra $\bfrak$ is the direct sum of $\hfrak$ and the root spaces corresponding to roots of finite length.  
	
	\begin{prop}
		Let $\bfrak$ be a Dynkin borel subalgebra of $\gfrak$.  Then,
		\begin{align}
			\left.\rho_{n+1}\right|_{\hfrak_n}=\rho_n
		\end{align}
		for every positive integer $n$.
		\label{prop:hspr}
		
		\begin{pf}
			Let $\left(\alpha_m\right)_{m\in\amsbb{Z}_{> 0}}$ be an ordering the simple positive roots of $\gfrak$ such that, for some positive integers $t_1<t_2<t_3<\ldots$, the $\amsbb{K}$-span of the coroots $h_{\alpha_{1}},h_{\alpha_2},\ldots,h_{\alpha_{t_n}}$ is precisely $\hfrak_n$.  Let $\omega_{n,i}\in\hfrak_n^*$ be the \emph{fundamental weight} associated to $\left.\alpha_i\right|_{\hfrak_n}$ (of $\gfrak_n$, with respect to the Cartan subalgebra $\hfrak_n$); in other words, 
			\begin{align}
				\omega_{n,i}\left(h_{\alpha_j}\right)=\deltaup_{i,j}
			\end{align} for $i,j\in\{1,2,\ldots,t_n\}$, where $\deltaup$ is the Kronecker delta.  We extend $\omega_{n,i}$ to $\omega_i\in\hfrak^*$ by setting 
			\begin{align}
				\omega_i\left(h_{\alpha_j}\right)=\left\{
				\begin{array}{ll}	
					\omega_{n,i}\left(h_{\alpha_j}\right)&\text{if }j\leq t_n\,,
					\\
					0&\text{if }j>t_n\,.
				\end{array}\right.
			\end{align} Note that $\omega_i$ does not depend on $n$.
			
			Now, we have 
			\begin{align}\rho_n=\sum_{i=1}^{t_n}\,\omega_{n,i}=\left(\sum_{i=1}^{t_n}\,\omega_i\right)\Biggr|_{\hfrak_n}\,.\end{align}  Similarly,  \begin{align}\rho_{n+1}=\left(\sum_{i=1}^{t_{n+1}}\,\omega_i\right)\Biggr|_{\hfrak_{n+1}}\end{align} so that $\left.\rho_{n+1}\right|_{\hfrak_n}=\left(\sum_{i=1}^{t_{n+1}}\,\omega_i\right)\Biggr|_{\hfrak_{n}}$.  Ergo,
			\begin{align}
				\left.\rho_{n+1}\right|_{\hfrak_n}-\rho_n=\left(\sum_{i=t_n+1}^{t_{n+1}}\,\omega_i\right)\Biggr|_{\hfrak_n}\,.
			\end{align}
			As $\omega_i\big|_{\hfrak_n}=0$ whenever $i>t_n$, the proposition follows.
		\end{pf}
	\end{prop}
	
	\begin{cor}
		There exists $\rho\in\hfrak^*$ such that $\rho|_{\hfrak_n}=\rho_n$ for every $n\in\amsbb{Z}_{> 0}$ if and only if $\bfrak$ is a Dynkin Borel subalgebra.  That is, a \emph{global half sum of $\bfrak$-positive roots} $\rho$ is well defined.  Furthermore, if $\gfrak$ is locally semisimple (i.e., each $\gfrak_n$ is semisimple) and $\bfrak$ is Dynkin, then $\rho$ is unique.
		
		\begin{pf}
			The uniqueness part when $\gfrak$ is locally semisimple is clear.  For the existence when $\bfrak$ is Dynkin, we consider $h=\sum_{i=1}^l \,k_l\,h_{\alpha_l}$ with $k_i\in\amsbb{K}$, where the $\alpha_i\in\hfrak^*$ and $h_{\alpha_i}\in\hfrak$ are as in the proof of Proposition~\ref{prop:hspr}.  Take $\rho\in\hfrak^*$ to be the map such that
			\begin{align}
				\rho(h)=\sum_{i=1}^l\,k_l\,.
			\end{align}
			Obviously, $\rho$ is well defined and $\rho|_{\hfrak_n}=\rho_n$ for every $n$.
			
			If $\bfrak$ is not Dynkin, then there is a positive root $\alpha$ of infinite length.  We can then write $\alpha$ as an arbitrary large sum of positive roots.  Suppose that 
			\begin{align}
				\alpha=\beta_1+\beta_2+\ldots+\beta_m
			\end{align} for some positive roots $\beta_1,\beta_2,\ldots,\beta_m$.  If $\rho$ exists, then we note that $\rho\left(h_\beta\right)$ is a positive integer for every positive root $\beta$.  Hence,
			\begin{align}
			\rho\left(h_\alpha\right)=\sum_{i=1}^m\,\rho\left(h_{\beta_i}\right)\geq m\,.
			\end{align}
			As $m$ can be arbitrarily large, $\rho\left(h_\alpha\right)$ is not defined, leading to a contradiction.
		\end{pf}
	\end{cor}

	From the corollary above, Dynkin Borel subalgebras play a distinguished role because the existence of the global half sum of positive roots $\rho$ allows us to define the \emph{dot action} of the Weyl group $W$ of $\gfrak$:
	\begin{align}
		w\cdot\lambda\defeq w(\lambda+\rho)-\rho
	\end{align}
	for all $w\in W$ and $\lambda\in\hfrak^*$.  
	
\subsection{The Lie Algebras $\gllie_\infty$, $\sllie_\infty$, $\solie_\infty$, and $\splie_\infty$}
\label{ch:glinfty}
	
	The material in this subsection is based on~\cite{penkov-styrkas}.  For every positive integer $k$, we shall write $\gllie_k$ for $\gllie_k(\amsbb{K})$.  The embeddings $\gllie_k\overset{\subseteq}{\longrightarrow}\gllie_{k+1}$\nomenclature[B]{$X\overset{\subseteq}{\longrightarrow}Y$}{the inclusion map from $X$ into $Y$, where $X\subseteq Y$}\nomenclature[B]{$\subseteq$}{is a subset of} sending 
	\begin{align}
		\textbf{X}\mapsto\left[
				\begin{array}{cc}
					\mathbf{X}		&	\boldsymbol{0}_{k\times1}	\\
					\boldsymbol{0}_{1\times k}	&	0
				\end{array}
				\right]\end{align}
for every $\textbf{X} \in \gllie_{k}$ make $\left(\gllie_k\right)_{k\in\amsbb{Z}_{>0}}$ a directed system.\footnote{Here, $\boldsymbol{0}_{m \times n}$ is the $m$-by-$n$ zero matrix for each $m,n\in\amsbb{Z}_{> 0}$}\nomenclature[B]{$\boldsymbol{0}_{m\times n}$}{the $m$-by-$n$ zero matrix for $m,n\in\amsbb{Z}_{> 0}$}  The direct limit $\lim_{\underset{k}{\boldsymbol{\longrightarrow}}} \,\gllie_k$ is denoted by $\gllie_\infty$.\index{direct limit}  For simplicity, we write $\gfrak$ for $\gllie_\infty$ and $\gfrak_n$ for the subalgebra $\gllie_n\subseteq\gfrak$, where $n\in\amsbb{Z}_{> 0}$. \nomenclature[B]{$\lim_{\boldsymbol{\longrightarrow}}$}{direct limit}
				
				Now, let $V$ and $V_*$ be countable-dimensional vector spaces along with a nondegenerate bilinear form $\langle\bullet,\bullet\rangle:V\times V_*\to \amsbb{K}$.\nomenclature[B]{$\langle\cdot,\cdot\rangle$}{a nondegenerate bilinear form between a natural $\gllie_\infty$-module and its restricted dual}  Then there exist ordered bases $\left(v^i\right)_{i \in \amsbb{Z}_{> 0}}$ of $V$ and $\left(v_*^j\right)_{j\in \amsbb{Z}_{> 0}}$ of $V_*$ such that 
				\begin{align}
					\left\langle v^i,v_*^j\right\rangle = \deltaup^{i,j}
				\end{align} for all $i,j\in\amsbb{Z}_{> 0}$, where $\deltaup$ is the Kronecker delta.\nomenclature[B]{$\deltaup$}{the Kronecker delta}  We can identify $V_*$ with a subspace of the dual space $V^*$ of $V$ via the identification $v_* \mapsto \left\langle \bullet,v_*\right\rangle$\nomenclature[B]{$\mapsto$}{is mapped to} for each $v_* \in V_*$.  We call $V_*$ the \emph{restricted dual} of $V$.\nomenclature[B]{$V_*$}{the restricted dual of a countable-dimensional vector space $V$} 
\nomenclature[B]{$\left(x_n\right)_{n \in \amsbb{Z}_{> 0}}$}{a sequence indexed by $\amsbb{Z}_{> 0}$}
			
			For every $n \in \amsbb{Z}_{> 0}$, we write $V^n$ for the $\amsbb{K}$-span of $\left\{v^1,v^2,\ldots,v^n\right\}$ and $V_*^n$ for the $\amsbb{K}$-span of $\left\{v_*^1,v_*^2,\ldots,v_*^n\right\}$.  Clearly, for each positive integer $n$, $V^n$ and $V_*^n$ have natural $\gfrak_n$-module structures and $V^n\tensor{} V_*^n$ is isomorphic as a Lie algebra to $\gfrak_n$ (the tensor product is defined over $\amsbb{K}$).  Trivially, 
			\begin{align}
				V=\lim\limits_{\underset{n}{\boldsymbol{\longrightarrow}}}\,V^n\text{ and }V_*=\lim\limits_{\underset{n}{\boldsymbol{\longrightarrow}}}\,V_*^n
			\end{align} (where the maps are the usual inclusions) are then $\gfrak$-modules, known as the \emph{natural $\gfrak$-module} and the \emph{conatural $\gfrak$-module}, respectively.\index{module over $\gllie_\infty$!natural ...}\index{module over $\gllie_\infty$!conatural ...}

			Note that $V\tensor{} V_*$ is an associative algebra under the multiplication defined by
			\begin{align}
				\left(u\otimes u_*\right) \boldsymbol{\cdot} \left(v\otimes v_*\right) \defeq  \left\langle v,u_*\right\rangle u\otimes v_*
			\,,
			\label{eq:multiplication}
			\end{align}
			for all $u,v\in V$ and $u_*,v_*\in V_*$.  Consequently, $\gfrak=\gllie_\infty$ is precisely $V\otimes V_*$ if one endows $V\otimes V_*$ with the Lie algebra structure associated to the multiplication (\ref{eq:multiplication}).
			
			Let $\boldsymbol{\prec}$ be a (strict) total order on the set $\amsbb{Z}_{> 0}$.\nomenclature[B]{$\boldsymbol{\prec}$}{a strict total order on $\amsbb{Z}_{> 0}$}\index{order!strict total order}  For $i,j\in\amsbb{Z}_{> 0}$, write $\textbf{E}_{i,j} \in \gfrak$ for the matrix with $1$ on the $(i,j)$-entry and $0$ everywhere else.\nomenclature[B]{$\textbf{E}_{i,j}$}{the element of $\gllie_\infty$ with $1$ on the $(i,j)$-entry and $0$ everywhere else} The subalgebra 
			\begin{align}\hfrak =\hfrak_\text{diag}:=\bigoplus\limits_{i \in \amsbb{Z}_{> 0}}\amsbb{K}\,\textbf{E}_{i,i}\end{align}\nomenclature[B]{$\hfrak_\text{diag}$}{the subalgebra of $\gllie_\infty$ consisting of diagonal elements} of diagonal elements is a splitting maximal toral subalgebra of $\gfrak$ (see~\cite{NP}).  The symbol $\epsilon_i$ represents the linear map in $\hfrak^*$ sending $h \in \hfrak$ to its $i$-th diagonal entry.\nomenclature[B]{$\epsilon_i$}{the map sending $x \in \gllie_\infty$ to its $i$-th diagonal entry}  Let $\bfrak$ be the splitting Borel subalgebra of $\gfrak$ with positive roots of the form $\epsilon_i - \epsilon_j$, where $i,j \in \amsbb{Z}_{> 0}$ such that $i \boldsymbol{\prec} j$.  Then we have a root space decomposition 
			\begin{align} \gfrak = \hfrak\oplus\bigoplus\limits_{\alpha \in \Delta} \gfrak^\alpha\,,
			\end{align} where $\Delta$ is the set of roots (with $\Delta^+$ representing the set of positive roots and $\Delta^-=-\Delta^+$ being the set of negative roots).  Hence, $\bfrak=\hfrak\oplus\nfrak$, where $\nfrak=\nfrak^+$ is the locally nilpotent subalgebra $\bigoplus\limits_{\alpha\in \Delta^+} \gfrak^\alpha$, whose opposite algebra\index{opposite Lie algebra} is 
			\begin{align}
				\nfrak^-=\bigoplus\limits_{\alpha\in\Delta^-}\gfrak^\alpha\end{align} (whence $\bfrak^-=\hfrak\oplus\nfrak^-$ is the opposite algebra of $\bfrak=\bfrak^+$).  
			
			Here are some Dynkin Borel subalgebras of $\gllie_\infty$ containing $\hfrak$.  The \emph{one-sided Dynkin Borel subalgebra}\index{Borel subalgebra of $\gllie_\infty$!one-sided standard ...} $\bfrak=\bfrak_\text{1st}=\bfrak_\text{1st}^+$\nomenclature[B]{$\bfrak_\text{1st}$}{the one-sided Dynkin Borel subalgebra of $\gllie_\infty$} is given by the  natural order on $\amsbb{Z}_{> 0}$: 
			\begin{align}
				1\boldsymbol{\prec} 2 \boldsymbol{\prec} 3 \boldsymbol{\prec} \ldots\,,
			\end{align} 
			whereas the \emph{two-sided Dynkin Borel subalgebra}\index{Borel subalgebra of $\gllie_\infty$!two-sided standard ...}\nomenclature[B]{$\bfrak_\text{2st}$}{the two-sided Dynkin Borel subalgebra of $\gllie_\infty$} $\bfrak=\bfrak_\text{2st}$ is given by the two-tailed ordering 
			\begin{align}
				\ldots \boldsymbol{\prec} 5\boldsymbol{\prec} 3\boldsymbol{\prec} 1 \boldsymbol{\prec} 2 \boldsymbol{\prec}4 \boldsymbol{\prec} 6 \boldsymbol{\prec} \ldots
			\,.
			\end{align}  
			The opposite subalgebra $\bfrak_\text{1st}^-$ to $\bfrak_\text{1st}$ is also a Dynkin Borel subalgebra and it is given by the ordering\index{Borel subalgebra of $\gllie_\infty$!standard ...} 
			\begin{align}
				\ldots \boldsymbol{\prec} 3 \boldsymbol{\prec} 2 \boldsymbol{\prec}1\,.
			\end{align}

The Lie algebras $\sllie_\infty$, $\solie_\infty$, and $\splie_\infty$ are subalgebras of $\gllie_\infty$ defined in the obvious manner.  That is, these Lie algebras are the direct limits of the finite-dimensional Lie algebras $\sllie_n$, $\solie_n$, and $\splie_{2n}$ with respect to the embeddings 
\begin{align}
\textbf{X}\mapsto \begin{bmatrix}\textbf{X}&\boldsymbol{0}_{n\times j}\\\boldsymbol{0}_{j\times n}&\boldsymbol{0}_{j\times j}\end{bmatrix}\,,\end{align}
where $j=1$ in the case of $\sllie_\infty$ and $\solie_\infty$, and $j=2$ in the case of $\splie_\infty$.    The following classification theorem of root-reductive Lie algebras is fully stated in~\cite{DPS} and proven in~\cite{DP1999}.

\begin{thm}
	Let $\gfrak$ be a root-reductive Lie algebra.  Then, $[\gfrak,\gfrak]$ is isomorphic to a direct sum of finite-dimensional simple Lie algebras, as well as copies of $\sllie_\infty$, $\solie_\infty$, and $\splie_\infty$, each with at most countable multiplicities.  Furthermore, $\gfrak$ is a semidirect sum $[\gfrak,\gfrak]\lsds\,\mathfrak{a}$ for some abelian Lie algebra $\mathfrak{a}$.  (The sign $\!\lsds\,$ denotes semidirect sum of Lie algebras, with the round side pointing towards the ideal.)
\end{thm}

For example, the root-reductive Lie algebra $\gllie_\infty$ is given by $\sllie_\infty\lsds\, \amsbb{K}$.  Note that this semidirect sum is not direct.
			
\pagebreak

\section{Some Examples}
\label{ch:exm}

\subsection{Countable-Dimensional Finitary Simple Lie Algebras}

	A result by Baranov and Strade~\cite{BS} gives a complete list of countable-dimensional finitary simple Lie algebras.  All such Lie algebras are isomorphic to $\sllie_\infty$, $\solie_\infty$, or $\splie_\infty$.  To clarify, $\sllie_\infty$ is the Lie algebra of traceless finitary matrices whose columns and rows are indexed by $\amsbb{Z}_{>0}$.  Then, $\solie_\infty$ is a subalgebra of $\sllie_\infty$  consisting of antisymmetric matrices.  Let $\textbf{J}$ denote the matrix
	\begin{align}
		\textbf{J}:=\left[\begin{array}{cc|cc|cc|c}
		0&1&0&0&0&0&\cdots\\
		-1&0&0&0&0&0&\cdots\\
		\hline
		0&0&0&1&0&0&\cdots\\
		0&0&-1&0&0&0&\cdots\\
		\hline
		0&0&0&0&0&1&\cdots\\
		0&0&0&0&-1&0&\cdots\\
		\hline
		\vdots&\vdots&\vdots&\vdots&\vdots&\vdots&\ddots
		\end{array}\right]\,.
	\end{align}
	Then $\splie_\infty$ consists of matrices $\textbf{X}\in\sllie_\infty$ such that $\textbf{X}^{\boldsymbol{\top}}\,\textbf{J}+\textbf{J}\,\textbf{X}=\boldsymbol{0}$.
	
	Let $\hfrak_\text{A}$ denote the subalgebra of $\sllie_\infty$ consisting of diagonal matrices.  Thus, $\hfrak_\text{A}$ is spanned over $\amsbb{K}$ by
	\begin{align}
		\textbf{E}_{k,k}-\textbf{E}_{k+1,k+1}\,,
	\end{align}
	for $k=1,2,3,\ldots$.  Up to automorphism of $\sllie_\infty$, $\hfrak_\text{A}$ is the unique splitting maximal toral subalgebra of $\sllie_\infty$.  As in Section~\ref{ch:glinfty}, let $\textbf{E}_{i,j}$ denote the matrix in $\gllie_\infty$ whose $(i,j)$-entry is $1$ and all other entries are $0$.   Write $\epsilon_k\in\hfrak_\text{A}^*$ to be the map sending $\textbf{E}_{k,k}$ to $1$ and $\textbf{E}_{j,j}$ to $0$ for all $j\neq k$.  
	
	For $\solie_\infty$, there are two distinct choices of splitting maximal toral subalgebras (up to automorphism of $\solie_\infty$).  Let $\hfrak_\text{B}$ denote the subalgebra of $\solie_\infty$ spanned over $\amsbb{K}$ by the matrices 
	\begin{align}
		\sqrt{-1}\left(\textbf{E}_{2k,2k+1}-\textbf{E}_{2k+1,2k}\right)
	\end{align} for $k=1,2,3,\ldots$.   The subalgebra $\hfrak_\text{D}$ is the $\amsbb{K}$-span of the matrices 
	\begin{align}
		\sqrt{-1}\left(\textbf{E}_{2k-1,2k}-\textbf{E}_{2k,2k-1}\right)
	\end{align} 
	for $k=1,2,3,\ldots$.  Up to automorphism of $\solie_\infty$, $\hfrak_\text{B}$ and $\hfrak_\text{D}$ are the only splitting maximal toral subalgebras of $\solie_\infty$.
	
	For $\splie_\infty$, up to automorphism of $\splie_\infty$, there is a unique splitting maximal toral subalgebra, which we shall denote by $\hfrak_\text{C}$.  This subalgebra $\hfrak_\text{C}$ is the $\amsbb{K}$-span of 
	\begin{align}
		\textbf{E}_{2k-1,2k-1}-\textbf{E}_{2k,2k}
	\end{align} 
	for $k=1,2,3,\ldots$.
	
	In the case of $\hfrak_\text{A}\subseteq \sllie_\infty$, we let $\epsilon_k\in\hfrak_\text{A}^*$ to be the map sending a diagonal matrix in $\hfrak$ to its $k$-th diagonal entry.  Let $\deltaup$ denote the Kronecker delta.  In the case of $\hfrak_\text{B}\subseteq\solie_\infty$, we define $\epsilon_k\in\hfrak_\text{B}^*$ to be the map sending 
	\begin{align}
		\sqrt{-1}\left(\textbf{E}_{2j,2j+1}-\textbf{E}_{2j+1,2j}\right)\mapsto\deltaup_{j,k}
	\end{align} for each $j=1,2,3,\ldots$.   In the case of $\hfrak_\text{D}\subseteq\solie_\infty$, we take $\epsilon_k\in\hfrak_\text{D}^*$ to be the map sending 
	\begin{align}
		\sqrt{-1}\left(\textbf{E}_{2j-1,2j}-\textbf{E}_{2j,2j-1}\right)\mapsto\deltaup_{j,k}
	\end{align}
	for all $j=1,2,3,\ldots$.  In the case of $\hfrak_\text{C}\subseteq \splie_\infty$, we set $\epsilon_k\in\hfrak_\text{C}^*$ to be the map sending 
	\begin{align}
		\textbf{E}_{2j-1,2j-1}-\textbf{E}_{2j,2j}\mapsto\deltaup_{j,k}
	\end{align}
	for every $j=1,2,3,\ldots$.

	Splitting Borel subalgebras $\bfrak$ of $\sllie_\infty$ containing $\hfrak_\text{A}$ are in a one-to-one correspondence with total orders $\prec$ on $\amsbb{Z}_{>0}$ via 
	\begin{align}	
		\bfrak :=\hfrak_\text{A}\oplus \bigoplus_{i\prec j}\,\amsbb{K}\,\textbf{E}_{i,j}\,.
	\end{align}
	Dynkin Borel subalgebras $\bfrak$ (defined in Section~\ref{ch:dynkin}) of $\sllie_\infty$ containing $\hfrak_\text{A}$  correspond to locally finite total orderings of $\amsbb{Z}_{>0}$.  Up to automorphism of $\sllie_\infty$, there are only two possible Dynkin Borel subalgebras.  One is $\bfrak_{\text{1A}}$ given by the natural order $1\prec 2\prec 3\prec\ldots$ on $\amsbb{Z}_{>0}$, and has the following Dynkin diagram $ \text{A}_\infty^{\text{1-sided}}$:
	\begin{align}
 \begin{dynkin}
    \dynkinline{1}{0}{2}{0};
    \dynkinline{2}{0}{3}{0};
    \dynkinline{3}{0}{4}{0};
    \dynkinline{4}{0}{5}{0};
    \dynkinline{5}{0}{6}{0};
    \dynkindots{6}{0}{7}{0};
    \foreach \x in {1,...,7}
    {
      \ifnum \x=7 {\fill (\dynkinstep*7,\dynkinstep*0) circle (0.01);}
       \else {
       \dynkindot{\x}{0}
       }
       \fi
    }
  \end{dynkin}
\end{align}
The other is $\bfrak_{\text{2A}}$ given by the two-sided ordering $\ldots \prec 5\prec 3 \prec 1 \prec 2\prec 4\prec6\prec\ldots$ on $\amsbb{Z}_{>0}$ and has the following Dynkin diagram $
  \text{A}^{\text{2-sided}}_\infty$: 
	\begin{align}\begin{dynkin}
    \dynkindots{1}{0}{2}{0};
    \dynkinline{2}{0}{3}{0};
    \dynkinline{3}{0}{4}{0};
    \dynkinline{4}{0}{5}{0};
    \dynkinline{5}{0}{6}{0};
    \dynkindots{6}{0}{7}{0};
    \foreach \x in {1,...,7}
    {
      \ifnum \ifnum \x=1 1 \else \ifnum \x=7 1 \else 0 \fi \fi =1 {\fill (\dynkinstep*\x,\dynkinstep*0) circle (0.01);}
       \else {
       \dynkindot{\x}{0}
       }
       \fi
    }
  \end{dynkin}
\end{align}

There are uncountably many non-Dynkin Borel subalgebras of $\sllie_\infty$ containing $\hfrak_\text{A}$, for instance, the Borel subalgebra associated to the following total ordering on $\amsbb{Z}_{>0}$:
	\begin{align}
		1\prec 3 \prec 5\prec \ldots \prec 6\prec 4 \prec 2\,.
	\end{align}
	There are also dense total orders\footnote{A partially ordered set $(\pazocal{P},\prec)$ is said to be \emph{dense} if, for any $a,b\in \pazocal{P}$ with $a\prec b$, there exists $c\in \pazocal{P}$ such that $a\prec c \prec b$.} on $\amsbb{Z}_{>0}$, making $\amsbb{Z}_{>0}$ isomorphic to $\amsbb{Q}$, $\amsbb{Q}\cup\{+\infty\}$, $\amsbb{Q}\cup\{-\infty\}$, and $\amsbb{Q}\cup\{-\infty,+\infty\}$ as an ordered set, and they correspond to ``highly non-Dynkin'' Borel subalgebras of $\sllie_\infty$ containing $\hfrak_\text{A}$.
	
For $\solie_\infty$, up to automorphism of $\solie_\infty$, there is one Dynkin Borel subalgebra $\bfrak_\text{B}$ containing $\hfrak_\text{B}$.  It is associated to the set of simple positive roots
\begin{align}
	\Sigma^+_\text{B}:=\left\{-\epsilon_1,\epsilon_1-\epsilon_2,\epsilon_2-\epsilon_3,\epsilon_3-\epsilon_4,\ldots\right\}
\end{align}
and has the following Dynkin diagram $\text{B}_\infty$:
	\begin{align}
  \begin{dynkin}
    \dynkindoubleline{2}{0}{1}{0};
    \dynkinline{2}{0}{3}{0};
    \dynkinline{3}{0}{4}{0};
    \dynkinline{4}{0}{5}{0};
    \dynkinline{5}{0}{6}{0};
    \dynkindots{6}{0}{7}{0};
    \foreach \x in {1,...,7}
    {
      \ifnum \x=7 {\fill (\dynkinstep*7,\dynkinstep*0) circle (0.01);}
       \else {
       \dynkindot{\x}{0}
       }
       \fi
    }
  \end{dynkin}
\end{align}

Up to automorphism of $\solie_\infty$, there is also one Dynkin Borel subalgebra $\bfrak_\text{D}$ containing $\hfrak_\text{D}$.  It is associated with the set of simple positive roots
\begin{align}
	\Sigma^+_\text{D}:=\left\{-\epsilon_1-\epsilon_2,\epsilon_1-\epsilon_2,\epsilon_2-\epsilon_3,\epsilon_3-\epsilon_4,\ldots\right\}
\end{align}
and has the following Dynkin diagram  $\text{D}_\infty$:
	\begin{align}
 \begin{dynkin}
    \dynkinline{1}{1}{2}{0};
    \dynkinline{1}{-1}{2}{0};
    \dynkinline{2}{0}{3}{0};
    \dynkinline{3}{0}{4}{0};
    \dynkinline{4}{0}{5}{0};
    \dynkinline{5}{0}{6}{0};
    \dynkindots{6}{0}{7}{0};
    \foreach \x in {2,...,7}
    {
        \ifnum \x=7 {\fill (\dynkinstep*7,\dynkinstep*0) circle (0.01);}
       \else {
       \dynkindot{\x}{0}
       }
       \fi
    }
    {\fill (\dynkinstep*1,\dynkinstep*-1) circle (\dynkinradius);}
     {\fill (\dynkinstep*1,\dynkinstep*1) circle (\dynkinradius);}
  \end{dynkin}
\end{align}

For $\splie_\infty$, there is only one Dynkin Borel subalgebra $\bfrak_\text{C}$ containing $\hfrak_\text{C}$ up to automorphism of $\splie_\infty$.  It is given by the set of simple positive roots
\begin{align}
	\Sigma^+_\text{C}:=\left\{-2\epsilon_1,\epsilon_1-\epsilon_2,\epsilon_2-\epsilon_3,\epsilon_3-\epsilon_4,\ldots\right\}
\end{align}
and has the following Dynkin diagram  $\text{C}_\infty$:
	\begin{align}
  \begin{dynkin}
    \dynkindoubleline{1}{0}{2}{0};
    \dynkinline{2}{0}{3}{0};
    \dynkinline{3}{0}{4}{0};
    \dynkinline{4}{0}{5}{0};
    \dynkinline{5}{0}{6}{0};
    \dynkindots{6}{0}{7}{0};
    \foreach \x in {1,...,7}
    {
      \ifnum \x=7 {\fill (\dynkinstep*7,\dynkinstep*0) circle (0.01);}
       \else {
       \dynkindot{\x}{0}
       }
       \fi
    }
  \end{dynkin}
  \end{align}
  
  Since the Lie algebras $\sllie_\infty$, $\solie_\infty$, and $\splie_\infty$ are locally simple, and the respective Borel subalgebras $\bfrak_\text{1A}$, $\bfrak_\text{2A}$, $\bfrak_\text{B}$, $\bfrak_\text{C}$, and $\bfrak_\text{D}$ are Dynkin, there is a unique global half sum $\rho$ of positive roots in each case.  In the case of $\bfrak_\text{1A}$, $\bfrak_\text{B}$, $\bfrak_\text{C}$, and $\bfrak_\text{D}$, 
  \begin{align}
  	\rho = \left\{
  		\begin{array}{ll}
  			-\sum_{k=1}^\infty\,(k-1)\,\epsilon_k\,,&\text{if }\bfrak=\bfrak_{\text{1A}}\text{ or }\bfrak=\bfrak_\text{D}\,,\\
  			-\sum_{k=1}^\infty\,(-1)^k\,\left\lfloor\frac{k}{2}\right\rfloor\,\epsilon_k\,,&\text{if }\bfrak=\bfrak_{\text{2A}}\,,\\
  			-\sum_{k=1}^\infty\,\frac{2k-1}{2}\,\epsilon_k\,,&\text{if }\bfrak=\bfrak_\text{B}\,,\\
  			-\sum_{k=1}^\infty\,k\,\epsilon_k\,,&\text{if }\bfrak=\bfrak_{\text{C}}\,,\\
  		\end{array}
  	\right.
  	\label{eq:rho}
\end{align}
If we represent $\lambda=\sum_{k=1}^\infty\,\lambda^k\,\epsilon_k\in\hfrak^*$ by the a sequence $\left(\lambda^1,\lambda^2,\ldots\right)$, then
  \begin{align}
  	\rho = \left\{
  		\begin{array}{ll}
  			(0,-1,-2,-3,\ldots)\,,&\text{if }\bfrak=\bfrak_{\text{1A}}\text{ or }\bfrak=\bfrak_\text{D}\,,\\
  			(0,1,-1,2,-2,3,-3,\ldots)\,,&\text{if }\bfrak=\bfrak_{\text{2A}}\,,\\
  			\left(-\frac{1}{2},-\frac{3}{2},-\frac{5}{2},\ldots\right)\,,&\text{if }\bfrak=\bfrak_\text{B}\,,\\
  			(-1,-2,-3,\ldots)\,,&\text{if }\bfrak=\bfrak_{\text{C}}\,,\\
  		\end{array}
  	\right.
\end{align}
  
  \subsection{$r$-Layered $\gllie_\infty$}
  
 Fix $r\in\amsbb{Z}_{>0}$.  For each $k\in \{2,\ldots,r\}$, let
 \begin{align}
 	\textbf{T}_k:=\sum_{j=1}^\infty\,\textbf{E}_{kj,kj}\,,
 \end{align}
 where $\textbf{E}_{i,j}$ is as defined in Section~\ref{ch:glinfty}.  The \emph{$r$-layered general linear Lie algebra}, denoted by $\gllie_{\infty}^{[r]}$, is given by
 \begin{align}
 	\gllie_{\infty}^{[r]}:=\left\{\begin{array}{ll}
 	\sllie_\infty\,,&\text{for }r=0\,,\\
 	\gllie_\infty\,,&\text{for }r=1\,,\\
 	\gllie_\infty\oplus\bigoplus_{k=2}^r\,\amsbb{K}\,\textbf{T}_k\,,&\text{for }r\geq 2\,,
 	\end{array}\right.
 	\label{eq:rlayered}
 \end{align}
 where the direct sum in the case $r\geq 2$ is a direct sum of vector spaces.  Since we have embeddings $\gllie_\infty^{[r]}\overset{\subseteq}{\longrightarrow}\gllie_\infty^{[r+1]}$ for all $r\in\amsbb{Z}_{>0}$, we set
 \begin{align}
 	\gllie_\infty^{[\omegaup]}:=\lim_{\underset{r}{\boldsymbol{\longrightarrow}}}\,\gllie^{[r]}_\infty = \bigcup_{r=1}^\infty\,\gllie^{[r]}_\infty\,,
 \end{align}
 where $\omegaup$ denotes the least infinite ordinal number.

 \begin{lem}
 	For $r\in\{1,2,\ldots\}\cup\{\omegaup\}$, there exists a splitting short exact sequence 
 \begin{align}
 	0\to\sllie_\infty\to\gllie^{[r]}_\infty\to\amsbb{K}^r\to0\,.
 	\label{eq:sesrgl}
 \end{align}
 	Thus, $\gllie^{[r]}$ is a semidirect sum $\sllie_\infty\lsds\,\amsbb{K}^r$.  This sum is, however, not direct.
 	
 	\begin{pf}
 		Equation (\ref{eq:rlayered}) implies that the short exact sequence (\ref{eq:sesrgl}) is splitting for finite $r$.  For $r=\omega$, we have $\gllie_\infty^{[\omegaup]}=\gllie_\infty\oplus\bigoplus_{k\geq 2}\,\amsbb{K}\textbf{T}_k$, which also implies that  (\ref{eq:sesrgl}) is splitting.
 		
 		Now, the semidirect sum $\gllie^{[r]}=\sllie_\infty\lsds\,\amsbb{K}^r$ is not direct because the only matrices commuting with $\sllie_\infty$ are scalar multiples of the identity matrix $\textbf{I}$.  Nonetheless, the matrix $\textbf{I}$ is not in $\gllie^{[r]}_\infty$.
 	\end{pf}
 \end{lem}

 \begin{rmk}
Let $r\in\{1,2,\ldots\}\cup\{\omegaup\}$.  we have the following filtration of $\gfrak=\gllie^{[r]}_\infty$ by ideals:
 \begin{align}
 	0\subsetneq \gllie_\infty^{[0]} \subsetneq \gllie_\infty^{[1]}\subsetneq \gllie_\infty^{[2]}\subsetneq \ldots \subsetneq \gllie_\infty^{[r]}\,,
 \end{align}
 Furthermore, for every nonnegative integer $k<r$,  $\gllie^{[k]}_\infty$ is not a direct summand of $\gllie_\infty^{[r]}$, and $\gllie^{[r]}_\infty=\gllie_\infty^{[k]}\lsds \,\amsbb{K}^{r-k}$. (Here, $\omegaup-k$ is equal to the ordinal number $\omegaup$ for every $k=0,1,2,\ldots$.)
\end{rmk}
 
 The subalgebra $\hfrak:=\hfrak_\text{diag}^{[r]}$ consisting of diagonal matrices is a splitting maximal toral subalgebra of $\gllie^{[r]}_\infty$.  Moreover, we have the following root space decomposition of $\gllie^{[r]}_\infty$:
 \begin{align}
\gllie^{[r]}_\infty=\hfrak\oplus\bigoplus_{i\neq j}\,\amsbb{K}\,\textbf{E}_{i,j}\,,
 \label{eq:rootglinfty}
 \end{align}
 This decomposition is identical to that of $\gllie_\infty$ (with respect to $\hfrak_\text{diag}$), or $\sllie_\infty$ (with respect to $\hfrak_\text{A}$)
 Hence, splitting Borel subalgebras $\bfrak$ of $\gllie^{[r]}_\infty$ are also described by total orders on $\amsbb{Z}_{>0}$.  Up to automorphism of $\gllie^{[r]}_\infty$, there are also two Dynkin Borel subalgebras containing $\hfrak$.
 
 \subsection{Twisted $\gllie_\infty$}
 
 For each $n\in\amsbb{Z}_{>0}$, let $\iota_n:\gllie_{2n}\to\gllie_{2(n+1)}$ be the embedding
 \begin{align}
 \textbf{X}\mapsto
 	\begin{bmatrix}
 		\textbf{X}&\boldsymbol{0}_{2n\times 1}&\boldsymbol{0}_{2n\times 1}
 		\\
 		\boldsymbol{0}_{1\times 2n}&\frac{1}{n}\,\text{Tr}(\textbf{X})&0
 		\\
 		\boldsymbol{0}_{1\times 2n}&0&0
 	\end{bmatrix}
 \end{align}
 for every $\textbf{X}\in\gllie_{2n}$.  The \emph{twisted general linear Lie algebra}, denoted by $\gllie^\#_\infty$, is the inductive limit of the inclusions $\iota_n$.  For simplicity, write $\gfrak_n$ for the image of $\gllie_{2n}$ in $\gllie^\#_\infty$, and $\gfrak$ for $\gllie^\#_\infty$ itself.
 
 It is evident that $[\gfrak,\gfrak]$ equals $\sllie_\infty$.  Therefore, we have the splitting short exact sequence of Lie algebras
 \begin{align}
 	0\to \sllie_\infty \to \gllie^\#_\infty \to \amsbb{K}\to 0\,,
 \end{align}
 where a section $\amsbb{K}\to\gllie^\#_\infty$ is given by $t\mapsto t\,\textbf{B}$ for every $t\in\amsbb{K}$, where
 \begin{align}
 	\textbf{B}:=
 	\left[
 		\begin{array}{cc|cc|cc|c}
 			1&0&0&0&0&0&\cdots\\
 			0&0&0&0&0&0&\cdots\\
 			\hline
 			0&0&1&0&0&0&\cdots\\
 			0&0&0&0&0&0&\cdots\\
 			\hline
 			0&0&0&0&1&0&\cdots\\
 			0&0&0&0&0&0&\cdots\\
 			\hline
 			\vdots&\vdots&\vdots&\vdots&\vdots&\vdots&\ddots
 		\end{array}
 	\right]
 \end{align}
 for every $t\in\amsbb{K}$.  Thus, $\gfrak=\gllie^\#_\infty$ is a semidirect sum $\sllie_\infty\lsds\,\amsbb{K}$.  However, $\gfrak\not\cong\gllie_\infty$ due to the following observation.  For $\gllie_\infty$, any maximal toral subalgebra surjects onto $\gllie_\infty/\left[\gllie_\infty,\gllie_\infty\right]$ via the quotient map (see~\cite{DPS}).  We shall illustrate a maximal toral subalgebra $\hfrak$ of $\gfrak$ which does not surject on to $\gfrak/[\gfrak,\gfrak]$.
 
 Let $\textbf{E}_{i,j}$ be defined as in Section~\ref{ch:glinfty}.  The maximal toral subalgebra $\hfrak$ is defined to be the $\amsbb{K}$-span of 
 \begin{align}
 \big\{\textbf{H}_{k}\suchthat{k\in\amsbb{Z}_{>0}}\big\}\cup \big\{\textbf{D}_k\,\suchthat{k\in\amsbb{Z}_{>0}}\big\}\,,
 \end{align}
 where 
 \begin{align}
 	\textbf{H}_k:=\textbf{E}_{2k-1,2k-1}+\textbf{E}_{2k,2k}-\textbf{E}_{2k+1,2k+1}-\textbf{E}_{2k+2,2k+2}
\end{align} and 
\begin{align}
	\textbf{D}_k:=\textbf{E}_{2k-1,2k}+\textbf{E}_{2k,2k-1}\,.
\end{align}  
Because $\hfrak$ lies entirely in $[\gfrak,\gfrak]$, we conclude that $\hfrak$ maps trivially under the quotient map $\gfrak\to\gfrak/[\gfrak,\gfrak]$.  Due to Theorem~\ref{thm:splittingmaxtoral}, $\hfrak$ is not a splitting maximal toral subalgebra.  The current development of our theory has not yet include the case where the maximal toral subalgebra is not splitting.

 For $\gllie_\infty^\#$, the subalgebra $\hfrak$ consisting of diagonal matrices is a splitting maximal toral subalgebra.  That is, we have the root space decomposition
 \begin{align}
 	\gllie_\infty^\#=\hfrak\oplus\bigoplus_{i\neq j}\,\amsbb{K}\,\textbf{E}_{i,j}\,.
 \end{align}

\pagebreak

\section{Finite-Dimensional Background}
\label{ch:keyresults}

	Suppose for now that $\gfrak$ is a finite-dimensional reductive Lie algebra, and $\hfrak$ is a Cartan subalgebra of $\gfrak$.  The subalgebras $\nfrak=\nfrak^+$ and $\nfrak^-$, as well as the sets $\Delta$, $\Delta^+$, and $\Delta^-$, are described as in Section~\ref{ch:rootreductive}. 

\subsection{Basics}
\label{ch:finitedimensional}
	
	The following three theorems are fundamental (see \cite{bggo} for more details).
	
	\begin{thm}
		Let $\lambda,\mu \in \hfrak^*$.
		\begin{enumerate}
			\item[(a)] There exists a unique simple submodule of $\Mlie(\lambda)$.\index{Lie submodule!simple submodule} \nomenclature[B]{$\defeq$}{is defined to be} This submodule is also a Verma module.
			\item[(b)] The dimension of $\Hom_{\Ulie(\gfrak)}\big(\Mlie(\lambda),\Mlie(\mu)\big)$ equals $0$ or $1$.\nomenclature[B]{$\Hom_R(M,N)$}{the set of $R$-module homomorphisms from $M$ to $N$}
			\item[(c)] Any nonzero $\Ulie(\gfrak)$-module homomorphism $\varphi:\Mlie(\lambda)\to\Mlie(\mu)$ is an embedding.  If such an embedding exists, then $\lambda \preceq \mu$.
		\end{enumerate}
		\label{thm:findim}
	\end{thm}

	\begin{thm}[\bf Verma's Theorem]
			Let $\lambda \in \hfrak^*$.   Given a positive root $\alpha$ such that 
			\begin{align}
				s_\alpha\cdot \lambda \preceq \lambda\,,
			\end{align}
			then there exists an embedding 
			\begin{align}
				\Mlie\left(s_\alpha\cdot\lambda\right)\overset{\subseteq}{\longrightarrow}\Mlie(\lambda)\,.
			\end{align}  (For $\alpha \in \Delta$, $s_\alpha$ is the reflection with respect to $\alpha$.)
			\index{Verma's Theorem}
			\label{thm:verma}
	\end{thm}
		
	\begin{thm}[\bf BGG Theorem]
		 For $\lambda,\mu\in\hfrak^*$, there exists a nontrivial $\gfrak$-module homomorphism from $\Mlie(\lambda)$ to $\Mlie(\mu)$ if and only if $\lambda$ is \emph{strongly linked}\index{strong linkage} to $\mu$, i.e., there exist positive roots $\alpha_1,\alpha_2,\ldots,\alpha_l$ such that
		\begin{align}
			\lambda = \left(s_{\alpha_l}\cdots s_{\alpha_2}s_{\alpha_1}\right)\cdot\mu \preceq \left(s_{\alpha_{l-1}}\cdots s_{\alpha_2}s_{\alpha_1}\right)\cdot\mu \preceq \ldots \preceq s_{\alpha_1}\cdot\mu \preceq \mu\,.
		\end{align}
		That is, every Verma submodule\index{Lie submodule!Verma submodule} of $\Mlie(\mu)$ is of the form $\Mlie(w\cdot \mu)$ for some element $w$ of the Weyl group.		
		\label{thm:stronglinkage}\index{BGG Theorem}
	\end{thm}
	
	For example, let $n\in\amsbb{Z}_{> 0}$ and $\gfrak:=\sllie_n$ (or, similarly, $\gllie_n$).  We can take $\hfrak$ to be the subalgebra of diagonal elements of $\gfrak$ and $\bfrak$ to be the subalgebra of upper triangular elements.   Let $\epsilon_i\in\hfrak^*$ be the map sending $h \in \hfrak$ to its $i$-th diagonal entry.   Then, each $\lambda\in\hfrak^*$ can be written as 
	\begin{align}
		\lambda=\sum_{i=1}^n \lambda^i \epsilon_i\,,
	\end{align} 
	where $\sum_{i=1}^n\,\lambda_i=0$ for $\gfrak=\sllie_n$.
	 We shall write $\lambda=\left(\lambda^1,\lambda^2,\ldots,\lambda^n\right)$ as a shorthand notation.   The Weyl group $W$ is the symmetric group on $n$ letters $\mathfrak{S}_n$ and acts by permuting $\epsilon_1,\epsilon_2,\ldots,\epsilon_n$.
	
\begin{exm}
Let $\gfrak$ be $\sllie_n$.  The Borel subalgebra $\bfrak$ is the subalgebra of upper triangular matrices, and the Cartan subalgebra $\hfrak$ is the subalgebra of diagonal matrices.

Take $n:=3$ and  $\lambda:=\rho=(1,0,-1)$. The Weyl group is the symmetric group on $\{1,2,3\}$.    We shall construct a filtration of the Verma module $\Mlie(\lambda)$ using three simple reflections $\varsigma_1$, $\varsigma_2$, and $\varsigma_3$, with the aim to find the (unique) simple submodule of $\Mlie(\lambda)$, i.e., the \emph{socle} of $\Mlie(\lambda)$.  

First, we may take $\varsigma_1 \in W$ to be the transposition $(1\;2)$\nomenclature[B]{$:=$}{is set to be}.   Thence,
	\begin{align}
		\lambda_1:=\varsigma_1\cdot\lambda=(-1,2,-1) \preceq \lambda\,.
	\end{align} 
	Next, with $\varsigma_2:=(2\;3)$, we have 
	\begin{align}
		\lambda_2:=\varsigma_2\cdot\lambda_1 =(-1,-2,3)\preceq \lambda_1
		\,.
	\end{align}  
	Finally, with $\varsigma_3:=\varsigma_1=(1\;2)$, we obtain
	\begin{align}
		\lambda_3:=\varsigma_3\cdot\lambda_2 = (-3,0,3) \preceq \lambda_2\,,
	\end{align} 
	which is antidominant\index{weight!antidominant weight}.  That is, $\Mlie\left(\lambda_3\right)$ is the unique simple submodule of $\Mlie(\lambda)$ by Theorem~\ref{thm:findim}(a).  Furthermore, due to Verma's Theorem (Theorem \ref{thm:verma}), we get the filtration\nomenclature[B]{$\subsetneq$}{is a proper subset of}
	\begin{align}
			0 \subsetneq \Mlie\left(\lambda_3\right) \subsetneq \Mlie\left(\lambda_2\right) \subsetneq \Mlie\left(\lambda_1\right) \subsetneq \Mlie\left(\lambda\right)
	\end{align} 
	of $\Mlie(\lambda)$ by Verma submodules.\footnote{We would get a different Verma filtration $0 \subsetneq \Mlie\left(\lambda_3\right)\subsetneq \Mlie\left(\lambda_2'\right) \subsetneq \Mlie\left(\lambda_1'\right) \subseteq \Mlie\left(\lambda\right)$ if we instead took $\varsigma_1:=(2\;\;\;3)$, $\varsigma_2:=(1\;\;\;2)$, and $\varsigma_3:=(2\;\;\;3)$.}    By refining this \emph{Verma filtration}, we obtain the following filtration of $\Mlie(\lambda)$:
	\begin{align}
		0 \subsetneq \Mlie\left(\lambda_3\right) \subseteq N_2 \subsetneq \Mlie\left(\lambda_2\right) \subsetneq N_1\subsetneq \Mlie\left(\lambda_1\right) \subsetneq N \subsetneq \Mlie\left(\lambda\right)\,,
		\label{eq:compseriesex}
	\end{align}
	where $N$, $N_1$, and $N_2$ are maximal proper submodules of $\Mlie(\lambda)$, $\Mlie\left(\lambda_1\right)$, $\Mlie\left(\lambda_2\right)$, respectively.  Hence, we conclude that $\Mlie(\lambda)$ is of length at least $6$.

	In fact, it can be easily shown that the filtration (\ref{eq:compseriesex}) is a composition series of $\Mlie(\lambda)$, i.e., that module $\Mlie(\lambda)$ has length $6$. .  Firstly, $\Mlie(\lambda)/N = \Llie(\lambda)$ and 
		\begin{align}
		N=\Mlie\left(\lambda_1\right)+\Mlie\left(\lambda_1'\right)\,,
		\end{align} 
		where $\lambda_1':=(1,-2,1)$.  This means $N/\Mlie\left(\lambda_1\right)=\Llie\left(\lambda'_1\right)$.  Secondly, $\Mlie\left(\lambda_1\right)/N=\Llie\left(\lambda_1\right)$ and 
		\begin{align}
		N_1=\Mlie\left(\lambda_2\right)+\Mlie\left(\lambda_2'\right)\,,
		\end{align} where $\lambda_2':=(-3,2,1)$.  Consequently, $N_1/\Mlie\left(\lambda_2\right)=\Llie\left(\lambda_2'\right)$.  Finally, $\Mlie\left(\lambda_2\right)/N_2=\Llie\left(\lambda_2\right)$ and  
		\begin{align}
			N_2=\Mlie\left(\lambda_3\right)=\Llie\left(\lambda_3\right)\,.
		\end{align}  
		Thus, all composition factors of $\Mlie(\lambda)$ are $\Llie\left(\lambda\right)$, $\Llie\left(\lambda_1\right)$, $\Llie\left(\lambda_1'\right)$, $\Llie\left(\lambda_2\right)$, $\Llie\left(\lambda_2'\right)$, and $\Llie\left(\lambda_3\right)$, each occurring with multiplicity $1$.\index{composition factor}\index{composition factor multiplicity}  That is, $\Mlie(\lambda)$ is indeed of length $6$.
		\label{exm:someexm}
\end{exm}

	\begin{figure}[h!]
  \caption{Geometry of the Singular Weights of $\Mlie\big((1,0,-1)\big)$ in Example~\ref{exm:someexm}.}
  \centering
    \includegraphics[width=\textwidth]{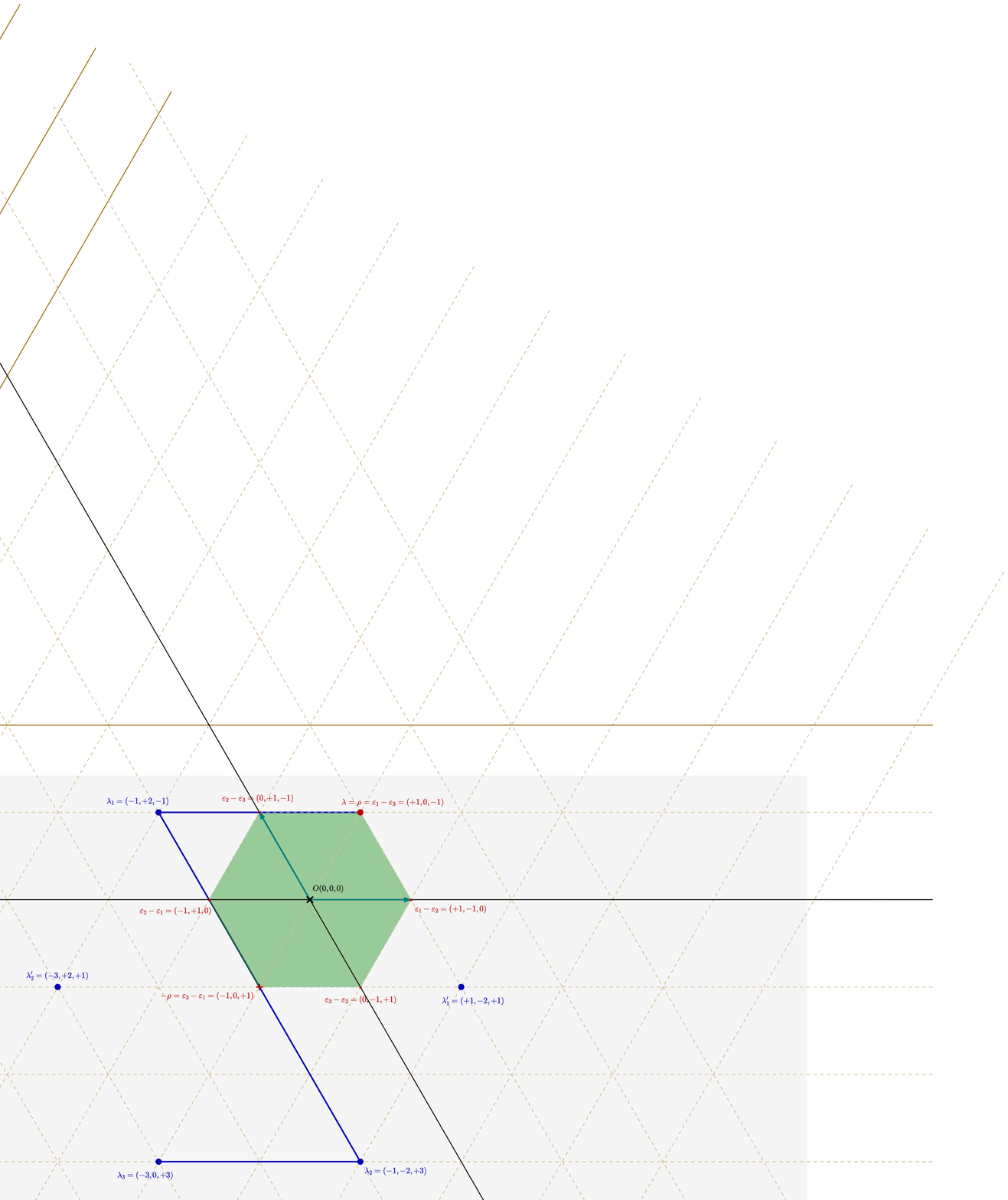}
\end{figure}
	
	
%
	
\subsection{The Bernstein-Gelfand-Gelfand Category $\pazocal{O}$}
	\label{ch:bggo}
	
	Consider the full subcategory $\pazocal{O}^\gfrak_\bfrak$ of the category of $\gfrak$-modules satisfying the conditions:\nomenclature[B]{$\pazocal{O}=\mathpzc{O}(\gfrak,\bfrak,\hfrak)$}{the BGG category $\pazocal{O}$ of a finite-dimensional semisimple Lie algebra $\gfrak$ with respect to a Borel subalgebra $\bfrak$ and a Cartan subalgebra $\hfrak\subseteq \bfrak$}
	\begin{itemize}
		\item[O1.]  Every $M \in \pazocal{O}^\gfrak_\bfrak$ is a finitely generated $\Ulie(\gfrak)$-module;
		\item[O2.]  Each $M \in \pazocal{O}^\gfrak_\bfrak$ is an $\hfrak$-weight $\gfrak$-module: $M=\bigoplus_{\lambda \in \hfrak^*} M^\lambda$; 	
		\item[O3.]  For each $M \in \pazocal{O}^\gfrak_\bfrak$  and $v \in M$, the subspace $\Ulie(\nfrak)\cdot v$ of $M$ is finite-dimensional (i.e., $M$ is \emph{locally $\nfrak$-finite}).
	\end{itemize}
	This category is known as the \emph{Bernstein-Gelfand-Gelfand (BGG) category $\pazocal{O}$} (see~\cite{bggo}).\index{BGG category $\pazocal{O}$}  
	In the remaining part of this subsection, we write $\bggO$ instead of $\bggO^\gfrak_\bfrak$.
	
	\begin{thm}
		The category $\pazocal{O}$ is both artinian and noetherian.  That is, every module in $\pazocal{O}$ are both artinian and noetherian over $\Ulie(\gfrak)$.\index{artinian module}\index{noetherian module}  In other words, each object in $\bggO$ has finite length.
	\end{thm}
	
	Many finite-dimensional results are based on the artinian and noetherian properties of the category $\pazocal{O}$.  These properties  do not hold in the infinite-dimensional case.  
	
	The category $\bggO$ has a direct sum decomposition into blocks.  The term ``block'' is defined in the same way as in~\cite{Penkov-Chirvasitu}, as given below.
	
	\begin{define}
		Let $\text{Indec}(\pazocal{C})$ denote the class of indecomposable objects of a given abelian category $\pazocal{C}$.  Suppose that $\text{Indec}(\pazocal{C})$ is a set.  The \emph{blocks} of $\pazocal{C}$ are equivalence classes of the finest equivalence relation $\approx$ on $\text{Indec}(\pazocal{C})$, requiring that, two indecomposable objects $X$ and $Y$ are equivalent when 
		\begin{align}
			\Hom_\pazocal{C}(X,Y)\neq 0\,.
		\end{align}
	\end{define}
	
	\begin{rmk}
		By abuse of language, if $\pazocal{B}$ is a full abelian subcategory of an abelian category $\pazocal{C}$ generated by a block (as described in the previous definition) of $\pazocal{C}$, then we also say that $\pazocal{B}$ is a \emph{block} of $\pazocal{C}$.   For example, the category of abelian groups consists of a single block---the category itself.
	\end{rmk}
	
	Let $\Zlie(\gfrak)$ be the \emph{center}\index{center} of $\Ulie(\gfrak)$.\nomenclature[B]{$\Zlie(L)$}{the center of $\Ulie(L)$}  An algebra homomorphism $\chi:\Zlie(\gfrak)\to \amsbb{K}$ is called a \emph{central character}.\index{central character}  We know from the standard theory (see~\cite{bggo}) that the blocks of $\pazocal{O}$ are of the form $\pazocal{O}\llbracket\lambda\rrbracket$\nomenclature[B]{$\pazocal{O}\llbracket\lambda\rrbracket$}{the block of $\pazocal{O}$ corresponding to a central character $\chi_\lambda$}, where $\lambda\in\hfrak^*$.  The subcategory $\pazocal{O}\llbracket\lambda\rrbracket$ consists of modules $M \in \pazocal{O}$ for which each weight of $M$ is in the shifted root lattice $\lambda+\Lambda$, and
	\begin{align}
		M=\Big\{v\in M\suchthat{}\forall z\in\Zlie(\gfrak), \exists n\in\amsbb{Z}_{> 0}: \big(z-\chi_\lambda(z)\big)^n\cdot v=0\Big\}\,,
	\end{align}
	where $\chi_\lambda$ is the central character corresponding to the weight $\lambda\in\hfrak^*$.  The \emph{principal block} is the block $\pazocal{O}\llbracket0\rrbracket$\index{block!principal block}.\nomenclature[B]{$\pazocal{O}_0$}{the principal block of $\pazocal{O}$}  Note that 
	\begin{align}
		\pazocal{O}=\bigoplus_{\llbracket\lambda\rrbracket\in \Omega}\, \pazocal{O}\llbracket\lambda\rrbracket\,,
	\end{align}
	where $\Omega$ is the set of equivalence classes $\llbracket\lambda\rrbracket$ under the equivalence relation $\sim$ on $\hfrak^*$ defined by $\lambda \sim \mu$ if $\lambda-\mu \in \Lambda$ and $\mu\in W\cdot\lambda$.  
	
	Unfortunately, for the Lie algebras of our interest, the enveloping algebras have trivial centers.  We shall need to devise a new method to study the blocks of (extended) categories $\bggO$ of such Lie algebras.
	
\begin{rmk}
	Let $\gfrak$ be a reductive Lie algebra of the form 
	\begin{align}
		\gfrak=\mathfrak{s}\oplus\mathfrak{a}\,,
	\end{align}
	 where  the semisimple ideal $\mathfrak{s}:=[\gfrak,\gfrak]$ of $\gfrak$ is finite dimensional and the center $\mathfrak{a}:=\text{\swabfamily z}(\gfrak)$ of $\gfrak$ is not necessarily finite dimensional (note that $\mathfrak{a}$ can be uncountable dimensional).   Then, a splitting maximal toral subalgebra of $\gfrak$ (and in fact, any Cartan subalgebra of $\gfrak$) is of the form 
	 \begin{align}
	 	\hfrak=\kfrak\oplus\mathfrak{a}\,,
	\end{align} where $\kfrak$ is a Cartan subalgebra of $\mathfrak{s}$.  A splitting Borel subalgebra containing $\hfrak=\kfrak\oplus\mathfrak{a}$ is of the form 
	\begin{align}
		\bfrak=\mathfrak{l}\oplus\mathfrak{a}\,,
	\end{align} 
	where $\mathfrak{l}$ is a Borel subalgebra of $\mathfrak{s}$ containing $\kfrak$.  
	
	The category $\bggO$ for $\gfrak$ with respect to the Borel subalgebra $\bfrak$ is defined in the same way as in the finite-dimensional case, and is denoted also by $\bggO^{\gfrak}_{\bfrak}$.  Most results from the finite-dimensional case carries over trivially to $\bggO^{\gfrak}_{\bfrak}$.
	\label{rmk:reductive}
\end{rmk}

\subsection{Kazhdan-Lusztig Theory}
	\label{ch:kltheory}\index{Kazhdan-Lusztig Theory}
	
	The Lie algebra $\gfrak$ is still a finite-dimensional reductive Lie algebra.  Statements in this subsection are based on~\cite{Coxeter} and~\cite{bggo}.
	
	\begin{thm}
		Let $(G,S)$ be a \emph{Coxeter system}\index{Coxeter system}, i.e., a group $G$ with identity $1_G$\nomenclature[B]{$1_G$}{the multiplicative identity of a group or a ring $G$} generated by a finite set $S$ with a presentation of the form 
		\begin{align} 
			\left\langle S\suchthat{} (st)^{m_{s,t}}=1_G\right\rangle\,,
		\end{align} where $m_{s,t}=m_{t,s}$ is a positive integer or $\infty$ for each $s,t\in S$, and $m_{s,s}=1$ for all $s \in S$. Then, there exists a unique partial order on $G$ (known as the \emph{Bruhat ordering}\index{order!weak partial order} )\index{Bruhat ordering}, denoted by $\preccurlyeq$, such that 
		\begin{itemize}
			\item[(i)] $1_G$ is the minimum element of $(G,\preccurlyeq)$ (i.e., $1_G \preccurlyeq g$ for every $g \in G$), and
			\item[(ii)] for each $g \in G$ and $s \in S$, if 
			\begin{align}
				\left[1_G,gs\right]=\left\{x\in G\suchthat{}1_G \preccurlyeq x \preccurlyeq gs\right\}
			\end{align} 
			has fewer elements than  
			\begin{align}
				\left[1_G,g\right]=\left\{x\in G\suchthat{}1_G \preccurlyeq x \preccurlyeq g\right\}\,,
			\end{align} then $\left[1_G,g\right]$ is the union of $\left[1_G,gs\right]$ and $\left[1_G,gs\right]s$.\nomenclature[B]{$[a,b]$}{the subset $\big\{x \in P\suchthat{} a\preccurlyeq x \preccurlyeq b\big\}$, where $(P,\preccurlyeq)$ is a (weakly ordered) poset}
		\end{itemize}
		The \emph{Bruhat length}\index{length!Bruhat length} of an element $g\in G$, denoted by $\ell(g)$, is defined to be the size of the interval $\left(1_G,g\right]:=\left\{x\in G\suchthat{}1_G \prec x \preccurlyeq g\right\}$.\nomenclature[B]{$\preccurlyeq$}{the Bruhat ordering on a Coxeter system} \nomenclature[B]{$\ell(g)$}{the Bruhat length of an element $g$ of a Coxeter group}\nomenclature[B]{$(a,b]$}{the subset $\big\{x \in P\suchthat{} a\prec x \preccurlyeq b\big\}$, where $(P,\preccurlyeq)$ is a (weakly ordered) poset} \nomenclature[B]{$\prec$}{the corresponding strict order of a weak order $\preccurlyeq$} \index{Coxeter group}
	\end{thm}
	
	\begin{define}
		Let $q$ be an indeterminate.   The \emph{Hecke algebra}\index{Hecke algebra} $\pazocal{H}=\mathpzc{H}(G,S)$ is a free module  generated by $\left\{T_g\right\}_{g\in G}$ over the ring $\amsbb{Z}\left[q^{-\frac{1}{2}},q^{+\frac{1}{2}}\right]$\nomenclature[B]{$R[S]$}{the ring of polynomials over a ring $R$ with indeterminates in the set $S$} of Laurent polynomials in $q^{\frac{1}{2}}$.  The multiplicative identity of $\pazocal{H}$ is $1_{\pazocal{H}}=T_{1_G}$ and the following multiplicative relations are satisfied: \nomenclature[B]{$\pazocal{H}=\mathpzc{H}(G,S)$}{the Hecke algebra of a Coxeter system $(G,S)$} 
		\begin{align}
			T_s^2 = (q-1)\,T_s + q\,T_{1_G}\,,
		\end{align}
		\begin{align}
			T_g\,T_s = T_{gs} \text{ if }g \prec gs\,,
		\end{align}
		and
		\begin{align}
			T_s\,T_g = T_{sg} \text{ if }g \prec sg\,,
		\end{align}
		for each $s \in S$ and $g \in G$.  Observe that $\pazocal{H}$ is an associative algebra.  Note that each $T_g$ is invertible, as 
		\begin{align}
			T_s^{-1}= q^{-1}\,T_s+\left(q^{-1}-1\right)\,T_{1_G}\,,
		\end{align}  
		for all $s \in S$.
		 There exists an \emph{involution}\index{involution} $\iota:\pazocal{H}\to\pazocal{H}$ sending $q^{+\frac{1}{2}}\mapsto q^{-\frac{1}{2}}$ and $T_g \mapsto T_{g^{-1}}^{-1}$ for all $g \in G$.  It is customary to write $\bar{X}$ for $\iota(X)$, where $X\in \pazocal{H}$.\nomenclature[B]{$T_g$}{the generator of $\pazocal{H}$ associated to the group element $g$} \nomenclature[B]{$\bar{X}$}{the involution of $X$ in $\pazocal{H}$}
	\end{define}
	
	\begin{thm}[\bf Kazhdan-Lusztig Theorem]
		There are unique elements $C_g$ with $g$ running over all $g\in G$ such that $C_g$ is fixed by the involution on $\pazocal{H}$ and 
		\begin{align}
			C_g = (-1)^{\ell(g)}\,q^{\frac{\ell(w)}{2}}\,\sum_{x \preccurlyeq g} \,(-1)^{\ell(x)}q^{-\ell(x)}\,\overline{P_{x,g}(q)}\,T_x\,,
		\end{align}
		where  $P_{x,g}(q) \in \amsbb{Z}[q]$ for every $x \in G$ with $x\preccurlyeq g$, $P_{g,g}(q)=1$, and 
		\begin{align}
			\deg\left(P_{x,g}(q)\right) \leq \frac{\ell(g)-\ell(x)-1}{2}
		\end{align} 
		whenever $x \prec g$.  The polynomials $P_{x,g}$ are known as the \emph{Kazhdan-Lusztig (KL) polynomials}\index{Kazhdan-Lusztig polynomials}.\nomenclature[A]{KL}{Kazhdan-Lusztig} \nomenclature[B]{$P_{x,g}$}{the Kazhdan-Lusztig polynomial for the group elements $x,g$ with $x \preccurlyeq g$}
	\end{thm}

	\begin{define}
		Let $A$ be a finite-dimensional algebra over an artinian ring $R$\index{artinian ring} (often, $R$ is a field).  The \emph{Grothendieck group} $\Klie(A)$ is the quotient of the free abelian group generated by the isomorphism classes $[X]$ of $A$-modules by the subgroup generated by elements of the form 
		\begin{align}[X]-[Y]+[Z]\,,
		\end{align}
		where
		\begin{align}
		0 \to X \to Y \to Z \to 0
		\end{align}\nomenclature[B]{$0_A$}{the additive identity of an abelian group or a ring $A$} is any exact sequence of $A$-modules $X$, $Y$, and $Z$.\index{Grothendieck group} \nomenclature[B]{$\Klie(A)$}{the Grothendieck group of an algebra $A$ over an artinian ring} \nomenclature[B]{[M]}{the element of the Grothendieck group associated to the module $M$}\index{exact sequence}
	\end{define}
	
	Let $\gfrak$ be a finite-dimensional semisimple Lie algebra, and $\hfrak$ a Cartan subalgebra contained in a Borel subalgebra $\bfrak$.\index{Lie subalgebra!Cartan subalgebra}   Write $\pazocal{O}$ for $\bggO^\gfrak_\bfrak$.  If $\lambda$ is a regular antidominant integral weight, then $\Mlie(\lambda)$ is simple (since $\lambda$ is antidominant\index{weight!antidominant weight}).  Write $M_w$ for $\Mlie(w\cdot\lambda)$ and $L_w$ for $\Llie(w\cdot\lambda)$.  The group $G$ is now replaced by the Weyl group $W$ of $\gfrak$.   The KL polynomials of the Weyl group play an important role in the computation of the multiplicities of the composition factors of Verma modules.  
	
	The theorem below was conjectured by Kazhdan and Lusztig in \cite{KL}.  It was proven independently in~\cite{BB} and in~\cite{BK}.
	
	\begin{thm}[{\bf Kazhdan-Lusztig Conjecture}]
		The composition factor multiplicities\index{composition factor multiplicity} of Verma modules are determined in $\Klie\big(\pazocal{O}\llbracket\lambda\rrbracket\big)$ by
		\begin{align}
			\left[L_w\right]=\sum_{x \preccurlyeq w}\,(-1)^{\ell(w)-\ell(x)}\,P_{x,w}(1)\,\left[M_x\right]\,,
		\end{align}
		where $w \in W$.  The inversion formula is
		\begin{align}
			\left[M_w\right]=\sum_{x\preccurlyeq w}\,P_{w^0w,w^0x}(1)\,\left[L_x\right]\,,
		\end{align}
		for all $w\in W$.  Here, $w^0$ denotes the longest element of $W$.
		\label{thm:KLconj}
	\end{thm}
	

\pagebreak

\section{The Extended Categories $\bggO$}
\label{ch:bbggo}

In this section, $\gfrak$ is a root-reductive algebra.  The subalgebra $\bfrak$ is any splitting Borel subalgebra of $\gfrak$ containing a splitting maximal toral subalgebra $\hfrak$.

\subsection{The Definition}
\label{sec:extbggo}

\begin{define}
	The \emph{extended category $\bggO$}, denoted by $\bar\bggO^\gfrak_\bfrak$, is the full subcategory of the category of $\gfrak$-modules satisfying the following two conditions:
	\begin{enumerate}[(i)]
		\item  Every $M\in\bar\bggO^\gfrak_\bfrak$ is an $\hfrak$-weight $\gfrak$-module with finite-dimensional $\hfrak$-weight spaces;
		\item Every $M\in\bar\bggO^\gfrak_\bfrak$ is locally $\nfrak$-finite (that is, $\Ulie(\nfrak)\cdot v$ is finite dimensional for every $v\in M$).
	\end{enumerate}
\end{define}

When this cannot cause confusion, we shall write $\bar\bggO$ for $\bar\bggO^{\gfrak}_\bfrak$.  Let us focus for now on the case $\gfrak$ is a finite-dimensional reductive Lie algebra and write $\bggO^\gfrak_\bfrak$ (or simply $\bggO$) for the classical BGG category $\bggO$ of $\gfrak$ with respect to the Borel subalgebra $\bfrak$.  As in Section~\ref{ch:bggo}, let $\Omega$ be the set of equivalence classes $\llbracket\lambda\rrbracket$, where $\lambda\in\hfrak^*$.  Then, it can easily be seen that each object $M$ of $\bar\bggO$ can be uniquely decomposed into a (possibly infinite) direct sum
\begin{align}
	M=\bigoplus_{\llbracket \lambda\rrbracket \in \Omega}\,M\llbracket\lambda\rrbracket\,,
\end{align}
where each $M\llbracket\lambda\rrbracket\in\bggO\llbracket\lambda\rrbracket$.  Note that the submodule $M\llbracket\lambda\rrbracket$ of $M$ is given by
\begin{align}
	M\llbracket\lambda\rrbracket := \Big\{v\in M\suchthat{\forall z\in\Zlie(\gfrak),\exists n\in\amsbb{Z}_{>0}:\,\big(z-\chi_\lambda(z)\big)^k\cdot v=0}\Big\}\,,
\end{align} 
where $\chi_\lambda$ is the central character corresponding to $\lambda\in\hfrak^*$ (for more details, see~\cite{bggo}).  If $M\in\bggO$, then such a direct sum must be finite.  Due to the proposition below, we can see that the blocks of $\bar\bggO$ are precisely the subcategories $\bggO\llbracket\lambda\rrbracket$.

\begin{prop}
	When $\gfrak$ is finite dimensional, each object of $\bbggO$ is a direct sum of objects in $\bggO$.
	\label{prop:findimdirsum}
	
	\begin{pf}	
		Let $M\in\bbggO$.  Then, $M$ can be written as a direct sum
		\begin{align}
			M=\bigoplus_{[\lambda]\in\hfrak^*/\Lambda}\,M^{[\lambda]}\,,
		\end{align}
		where the submodule $M^{[\lambda]}$ is given by
		\begin{align}
			M^{[\lambda]}:=\sum_{\mu\in[\lambda]}\,M^\mu
		\end{align}
		for each $[\lambda]\in\hfrak^*/\Lambda$.  Hence, it suffices to prove the proposition under the assumption that $M=M^{[\lambda]}$ for some $[\lambda]\in\hfrak^*/\Lambda$.
		
		Note that $M$ is countable dimensional.  It is generated by countably many weight vectors $v_1,v_2,\ldots$.  Write $M_n$ for the submodule of $M$ generated by $v_1,v_2,\ldots,v_n$.  Then, $M_n\in\bggO$.  We then can write 
		\begin{align}
			M_n=\bigoplus_{\left\llbracket \nu\right\rrbracket\in\Omega}\,M_n\left\llbracket\nu\right\rrbracket\,.
		\end{align}
		It can be easily seen that we have the inclusion $M_n\llbracket\nu\rrbracket\subseteq M_{n+1}\llbracket\nu\rrbracket$.  Write $M\left\llbracket\nu\right\rrbracket$ for the direct limit $\lim_{\underset{n}{\boldsymbol{\longrightarrow}}}\,M_n\llbracket\nu\rrbracket$.   Clearly, $M=\bigoplus_{\llbracket\nu\rrbracket\in\Omega}\,M\llbracket\nu\rrbracket$, whence it is sufficient to further assume that $M=M\llbracket\nu\rrbracket$ for some $\nu\in\hfrak^*$.
		
		We claim that $M=M_n$ for some sufficiently large $n$.  If this is not true, then there are infinitely many positive integers $n$ such that $M_n$ is a proper submodule of $M_{n+1}$, whence $M_{n+1}$ introduces at least one composition factor $L_n$ into $M_n$.  For such $n$, we have a simple subquotient $L_n$.  As there are only finitely many simple objects in $\bggO\llbracket\nu\rrbracket$ (up to isomorphism), there exists a simple object $L$ that appears as $L_n$ for infinitely many $n$.  Therefore, if $\xi$ is the highest weight of $L$, then the dimension of $M^\xi$ must be infinite, which is absurd.  Thus, the claim holds.
	\end{pf}
\end{prop}

We still assume that $\gfrak$ is a finite-dimensional reductive Lie algebra.  By Theorem~\ref{thm:directsuminjectives} below, $\bar\bggO$ has enough injectives.  Consequently, $\bar\bggO$, just like $\bggO$, is also a highest-weight category in the sense of Definition~\ref{def:hwcat} with respect to the partially ordered set $\hfrak^*$, the simple objects $\Llie(\lambda)=\Llie(\lambda;\gfrak,\bfrak,\hfrak)$, the co-standard objects $\Vlie(\lambda)=\Vlie(\lambda;\gfrak,\bfrak,\hfrak)$, and the injective objects $\Ilie(\lambda)=\Ilie(\lambda;\gfrak,\bfrak,\hfrak)$.  Here, $\Vlie(\lambda)$ is the co-Verma module with socle $\Llie(\lambda)$ and $\Ilie(\lambda)$ is the injective hull of $\Llie(\lambda)$ for each $\lambda\in\hfrak^*$.  The standard objects in the sense of Theorem~\ref{thm:brauer-humphreys} are the Verma modules $\Mlie(\lambda)=\Mlie(\lambda;\gfrak,\bfrak,\hfrak)$.   See also~\cite{hwcat}.

\begin{thm}
	For a finite-dimensional reductive Lie algebra $\gfrak$, every object $M$ in $\bbggO$ has a unique direct sum decomposition $M=\bigoplus_{j \in J}\,M_j$ where each $M_j$ is an indecomposable object in $\bggO$.  If $I_j$ is the injective hull of $M_j$ in $\bggO$, then the injective hull of $M$ in $\bbggO$ is $I:=\bigoplus_{j\in J}\,I_j$.
	\label{thm:directsuminjectives}
	
	\begin{pf}
		Since each block of $\bbggO$ is a block of $\bggO$, we have a unique direct sum decomposition $M=\bigoplus_{\llbracket\lambda\rrbracket\in\Omega}\,M\llbracket\lambda\rrbracket$ with $M\llbracket\lambda\rrbracket\in\bggO\llbracket\lambda\rrbracket$.  Now, since the objects in $\bggO\llbracket\lambda\rrbracket$ have finite length, $M\llbracket\lambda\rrbracket$ has a unique decomposition as a direct sum of indecomposable objects.  Thus, $M$ has a unique direct sum decomposition with indecomposable direct summands.
		
		To prove the second part of the theorem, let $0\to N\to L$ be an arbitrary exact sequence in $\bbggO$ along with a homomorphism $N\to I$.   Fix $\llbracket\lambda\rrbracket\in\Omega$.  By restricting the objects and morphisms to the block $\bbggO\llbracket\lambda\rrbracket=\bggO\llbracket\lambda\rrbracket$, we have the following diagram
		\begin{equation}
		\begin{tikzcd}
		0\arrow{r}& N\llbracket\lambda\rrbracket \arrow{r}{} \arrow{d}{} & L\llbracket\lambda\rrbracket\\
		&I\llbracket\lambda\rrbracket \,,&
		\end{tikzcd}
		\label{eq:findiminj}
		\end{equation}
		where the row is exact.  Note that $I\llbracket\lambda\rrbracket=\bigoplus_{j\in J}\,I_j\llbracket\lambda\rrbracket$.  For each $j\in J$, either $M_j\llbracket\lambda\rrbracket=0$ or $M_j\llbracket\lambda\rrbracket=M_j$ (as $M_j$ is indecomposable).  Because the weight spaces of $M$ are finite dimensional and there are only finitely many simple objects (up to isomorphism) in $\bggO\llbracket\lambda\rrbracket$, there are only finitely many $j\in J$ such that $M_j\llbracket\lambda\rrbracket \neq 0$.  Hence, $I_j\llbracket\lambda\rrbracket=0$ for all but finitely $j\in J$, and for $j\in J$ which $M_j\llbracket\lambda\rrbracket=M_j$, we see that $I_j\llbracket\lambda\rrbracket=I_j$ is injective in $\bggO$.  Therefore, $I\llbracket\lambda\rrbracket$ is a finite direct sum of some of the injective modules $I_j$.  Consequently, $I\llbracket\lambda\rrbracket$ is an injective object of $\bggO$.  Thus, there exists a map $L\llbracket\lambda\rrbracket\to I\llbracket\lambda\rrbracket$ such that the diagram below commutes:
		\begin{equation}
		\begin{tikzcd}
		0\arrow{r}& N\llbracket\lambda\rrbracket \arrow{r}{} \arrow{d}{} & L\llbracket\lambda\rrbracket\arrow{dl}{}\\
		&I\llbracket\lambda\rrbracket \,.&
		\end{tikzcd}
		\label{eq:findiminj2}
		\end{equation}
		
		Since $\llbracket\lambda\rrbracket$ is arbitrary, we take the direct sum of diagrams (\ref{eq:findiminj2}) over $\llbracket\lambda\rrbracket\in\Omega$ to get the following commutative diagram
		\begin{equation}
		\begin{tikzcd}
		0\arrow{r}& N \arrow{r}{} \arrow{d}{} & L\arrow{dl}{}\\
		&I \,.&
		\end{tikzcd}
		\label{eq:findiminj3}
		\end{equation}
		Hence, $I$ is indeed an injective object of $\bbggO$.
	\end{pf}
\end{thm}

\begin{define}
	A subcategory $\pazocal{C}$ of the category of $\amsbb{K}$-vector spaces is a \emph{highest-weight category} if the following conditions are met:
	\begin{enumerate}[(i)]
		\item $\pazocal{C}$ is \emph{locally artinian} in the sense that it admits arbitrary direct limits of subobjects and every object is a union of its subobjects of finite lengths;
		\item $\pazocal{C}$ has enough injectives;
		\item For every family of subobjects $\left\{A_\alpha\right\}_{\alpha \in J}$ and a subobject $B$ of $X\in\pazocal{C}$, we have
		\begin{align}
			B\cap\left(\bigcup_{\alpha \in J}\,A_\alpha\right)=\bigcup_{\alpha \in J}\,\left(B\cap A_\alpha\right)\,;
		\end{align}
		\item There exists a locally finite partially ordered set $(\mathcal{P},\preceq)$ which indexes an exhaustive family $\big\{S(\lambda)\big\}_{\lambda \in \mathcal{P}}$ of nonisomorphic simple objects in $\pazocal{C}$;
		\item The partially ordered  set $\mathcal{P}$ also indexes the family $\big\{A(\lambda)\big\}_{\lambda\in\mathcal{P}}$ of \emph{co-standard objects} of $\pazocal{C}$ such that there exist embeddings $S(\lambda)\to A(\lambda)$ and all composition factors $S(\mu)$ of $A(\lambda)/S(\lambda)$ satisfies $\mu\prec\lambda$;
		\item For all $\mu,\lambda\in \mathcal{P}$, the dimension of $\Hom_\pazocal{C}\big(A(\lambda),A(\mu)\big)$ is finite and the multiplicity $\big[A(\lambda):S(\mu)\big]$ is also finite; and
		\item Each $S(\lambda)$ has an injective hull $I(\lambda)$ in $\pazocal{C}$ equipped with an increasing filtration, possibly of infinite length,  called a \emph{co-standard filtration},
		\begin{align}
			0=F_0(\lambda)\subseteq F_1(\lambda)\subseteq F_2(\lambda)\subseteq \ldots\subseteq I(\lambda)
		\label{eq:goodfiltration}
		\end{align}
		such that
		\begin{enumerate}[1.]
			\item $F_1(\lambda)=A(\lambda)$,
			\item for $n>1$, $F_{n}(\lambda)/F_{n-1}(\lambda)\cong A\big(\mu(n)\big)$ for some $\mu(n)>\lambda$,
			\item for each $\nu\in\mathcal{P}$, $\mu(n)=\nu$ for only finitely many $n$, 
			\item $\bigcup_{n\geq 0}\,F_n(\lambda)=I(\lambda)$.
		\end{enumerate}
	\end{enumerate}	
	\label{def:hwcat}
\end{define}

\begin{thm}[Brauer-Humphreys Reciprocity,~\cite{hwcat}]
	Let $\pazocal{C}$ be a highest-weight category.  Then, $\pazocal{C}$ also contains enough projectives.  With notations as in Definition~\ref{def:hwcat}, we have the equality
	\begin{align}
		\big[V(\mu):S(\lambda)\big]=\big\{I(\lambda):A(\mu)\big\}\,,
	\end{align}
	where the multiplicity $\big\{I(\lambda):A(\mu)\big\}$ is the number of $n$ such that $F_n(\lambda)/F_{n-1}(\lambda)\cong A(\mu)$, with $\big(F_n(\lambda)\big)_{n\geq 0}$ as defined in (\ref{eq:goodfiltration}).    Here, for every $\lambda\in\Lambda$, $V(\lambda)$ denotes the largest quotient module of the projective cover of $S(\lambda)$ in $\pazocal{C}$, and it is known as the \emph{standard object} with respect to $\lambda$.
	\label{thm:brauer-humphreys}
\end{thm}

	Let us return to the case where $\gfrak$ may be infinite dimensional.  We can define the duality functor of the category $\bbggO$ in the same manner as the standard duality of the category $\bggO$.  More precisely, we have the following definition.

\begin{define}
	Let $M\in\bar\bggO$.  For $\lambda\in\hfrak^*$, let $M^{\vee,\lambda}$ denote the set of $f\in M^*$ such that $f$ vanishes on $M^\mu$ for every $\mu\in\hfrak^*\smallsetminus\{\lambda\}$.  The \emph{dual} of $M$ is defined to be $M^\vee\defeq \bigoplus_{\lambda\in\hfrak^*}\,M^{\vee,\lambda}$.
\end{define}

Now, if $\left\{x_{\pm\alpha}\,\boldsymbol{|}\,\alpha\in\Delta^+\right\}\cup \left\{h_\beta\,\boldsymbol{|}\,\beta\in\Sigma^+\right\}$ is a \emph{Chevalley basis}~\cite{humphreys,bggo} of $\gfrak$, then the \emph{transpose map} $\tau:\gfrak\to\gfrak$ is the linear map given by $\tau(h):=h$ for all $h\in\hfrak$, and $\tau\left(x_{\pm\alpha}\right):=x_{\mp\alpha}$ for all positive roots $\alpha$.  Note that $\big[\tau(x),\tau(y)\big]=-\tau\big([x,y]\big)$ for all $x,y\in\gfrak$.

\begin{thm}
	For every $M\in\bar\bggO$, $M^\vee$ is a $\gfrak$-module with respect to the twisted $\gfrak$-action $(g\cdot f)(v)\defeq f\big(\tau(g)\cdot v\big)$ for all $g\in\gfrak$, $v\in M$, and $f\in M^\vee$.  Furthermore, $M^\vee \in \bar\bggO$.
	
	\begin{pf}
		For $g_1,g_2\in\gfrak$, $f\in M^\vee$, and $v\in M$, we have
		\begin{align}
			\big(\left[g_1,g_2\right]\cdot f\big)(v)&=f\Big(\tau\big(\left[g_1,g_2\right]\big)\cdot v\Big)=f\Big(-\big[\tau\left(g_1\right),\tau\left(g_2\right)\big]\cdot v\Big)
			\nonumber\\
			&=f\Big(-\tau\left(g_1\right)\cdot\big(\tau\left(g_2\right)\cdot v\big)+\tau\left(g_2\right)\cdot\big(\tau\left(g_1\right)\cdot v\big)\Big)
		\end{align}
		That is,
		\begin{align}
			\big(\left[g_1,g_2\right]\cdot f\big)(v)&=f\Big(\tau\left(g_2\right)\cdot\big(\tau\left(g_1\right)\cdot v\big)\Big)-f\Big(\tau\left(g_1\right)\cdot\big(\tau\left(g_2\right)\cdot v\big)\Big)
			\nonumber\\
			&=\left(g_2\cdot f\right)\big(\tau\left(g_1\right)\cdot v\big)-\left(g_1\cdot f\right)\big(\tau\left(g_2\right)\cdot v\big)
			\nonumber\\
			&=\big(g_1\cdot\left(g_2\cdot f\right)\big)(v)-\big(g_2\cdot\left(g_1\cdot f\right)\big)(v)\,.
		\end{align}
		Consequently, $\left[g_1,g_2\right]\cdot f= g_1\cdot\left(g_2\cdot f\right)-g_2\cdot\left(g_1\cdot f\right)$, as desired.
		
		Assume now that $f\in M^{\vee,\lambda}$ for some weight $\lambda\in\hfrak^*$.  Then, for each $h\in\hfrak$, we have $h\cdot f=\lambda(h)\cdot f$, and it follows that $M^{\vee,\lambda}=\left(M^\vee\right)^{\lambda}$.  Thus, from the definition $M^\vee = \bigoplus_{\lambda\in\hfrak^*}\,M^{\vee,\lambda}$, $M^\vee$ is $\hfrak$-semisimple.  Furthermore, if $g\in\gfrak^{\alpha}$ for some positive root $\alpha$, then $g\cdot f\in M^{\vee,\lambda+\alpha}$.  Since $M$ is locally $\nfrak$-finite, we conclude that $M^{\vee}$ is also locally $\nfrak$-finite.  Finally, we clearly have $\dim_\amsbb{K}\left(M^{\vee,\lambda}\right)=\dim_\amsbb{K}\left(M^\lambda\right)<\infty$.  That is, $M^\vee$ has finite-dimensional $\hfrak$-weight spaces.  Ergo, $M^\vee\in\bar\bggO$.
	\end{pf}
\end{thm}
	
	Finally, we consider the categories $\bggO$ and $\bbggO$ for a finite-dimensional reductive Lie algebra $\gfrak$. With $\bar\bggO$ and $\bggO$ being highest-weight categories, we automatically have BGG reciprocity (which is a special case of Brauer-Humphreys reciprocity stated in Theorem~\ref{thm:brauer-humphreys}):
	\begin{align}
		\big[\Vlie(\mu):\Llie(\lambda)\big]=\big[\Mlie(\mu):\Llie(\lambda)\big]=\big\{\Ilie(\lambda):\Vlie(\mu)\big\}=\big\{\Plie(\lambda):\Mlie(\mu)\big\}\,,
	\end{align}
	where $\Vlie(\mu)$ is the dual Verma module $\big(\Mlie(\mu)\big)^\vee$, $\Ilie(\lambda)$ is the injective hull of $\Llie(\lambda)$, and $\Plie(\lambda)=\big(\Ilie(\lambda)\big)^\vee$ is the projective cover of $\Ilie(\lambda)$.  However, as we shall later prove, the category $\bar\bggO^\gfrak_\bfrak$ is not a highest-weight category if $\gfrak$ is infinite dimensional.

\subsection{Direct Sum Decompositions}

Let $\gfrak$ be a root-reductive Lie algebra.  The objective of this section is to verify that every object in $\bar\bggO$ has a decomposition into a direct sum of indecomposable objects.  Furthermore, this decomposition is unique up to isomorphism.  That is, if $M\in\bar\bggO$ can be written as $\bigoplus_{j\in J}\,M_j=M=\bigoplus_{j'\in J'}\,M_{j'}'$, where $J$ and $J'$ are index sets and $M_j,M'_{j'}\in\bar\bggO$ are indecomposable for all $j\in J$ and $j'\in J'$, then there exists a bijection $\psi:J\to J'$ such that $M_j\cong M'_{\psi(j)}$ for every $j\in J$.

First, we need the lemma below.  This lemma is in fact equivalent to Zorn's Lemma (which is equivalent to the Axiom of Choice).  See~\cite{settheory} for more information.

\begin{lem}[\bf Hausdorff Maximal Principle]
	Let $\mathcal{F}$ be a nonempty family of subsets of a given set $S$ such that, for every chain $\mathcal{C}\subseteq\mathcal{F}$, the union $\bigcup_{C\in\mathcal{C}}\,C$ belongs to $\mathcal{F}$.  Then, $\mathcal{F}$ has a maximal element with respect to inclusion.
\end{lem}

\begin{cor}
	For a nonempty partially ordered set $\mathcal{P}$, every chain $\mathcal{C}$ in $\mathcal{P}$ is contained in a maximal chain.
	\label{cor:hausdorff}
\end{cor}

\begin{thm}
	Every $M\in\bar\bggO$ is a direct sum of indecomposable objects.
	\label{thm:directsumdecomposition}
	
	\begin{pf}
		For an object $M\in\bar\bggO$, we say that $\xi\in\hfrak^*$ is a \emph{decomposable weight} of $M$ if there exist submodules $N$ and $L$ of $M$ such that $M=N\oplus L$ with $\dim_\amsbb{K}\left(N^\xi\right)>0$ and $\dim_\amsbb{K}\left(L^\xi\right)>0$.  Otherwise, $\xi$ is an \emph{indecomposable weight}.  (Note that, by abuse of language, an indecomposable weight of a $\gfrak$-module $X$ need not be a weight of $X$.  In other words, if $X^\xi=0$, then $\xi$ is an indecomposable weight of $X$, despite not actually being a weight of $X$.)

		For a semisimple $\hfrak$-module $X$, the \emph{support} $\text{supp}(X)$ of $X$ is the set of the $\hfrak$-weights of $X$.  For a subset $S\subseteq\hfrak^*$, we say that $S$ is an \emph{indecomposable weight set} of $M$ if every weight in $S$ is an indecomposable weight of $M$ and if $M$ cannot be written as a direct sum $M=N\oplus L$ such that $\text{supp}(N)\cap S$ and $\text{supp}(L)\cap S$ are both nonempty.
	
		Assume that $M\in\bar\bggO$ is nonzero.  First, let $\triangleleft$ be a well order on the set $\text{supp}(M)$.  The order $\triangleleft$ exists by the Well-Ordering Principle.  
		
		We shall prove by transfinite induction that, for each $\mu\in\text{supp}(M)$, there exists an index set $J_\mu$ such that $M$ has a direct sum decomposition
		\begin{align}
			M=\bigoplus_{j\in J_\mu}\,D_\mu(j)
		\end{align}
		such that the set $Q_\mu:=\big\{\nu\in\text{supp}(M)\suchthat{}\nu\trianglelefteq \mu\big\}$ is an indecomposable weight set of $D_\mu(j)$.  We further require that the decomposition above (with respect to the weight $\mu$) be compatible with the decomposition with respect to every weight $\upsilon$ satisfying $\upsilon\trianglelefteq \mu$ in the sense that, for any $j\in J_\mu$, there exists (uniquely) $i\in J_\upsilon$ such that
		\begin{align}
			D_\mu(j)\subseteq D_\upsilon(i)\,.
		\end{align}
		
		For the base case, let $\lambda$ be the least element of $\text{supp}(M)$ with respect to $\triangleleft$.  If $\lambda$ is already an indecomposable weight, then $M=M$ is a required decomposition.  If $\lambda$ is a decomposable weight, then there exists a submodule $D_\lambda(1)$ of $M$ such that $D_\lambda(1)$ is a direct summand of $M$.  We may chose $D_\lambda(1)$ so that $\dim_\amsbb{K}\Big(\big(D_\lambda(1)\big)^\lambda\Big)>0$ is minimal.  Then, $\lambda$ is an indecomposable weight of $D_\lambda(1)$.  Let $D'_\lambda(1)$ denote the complementary submodule of $M$ to $D_\lambda(1)$.  We proceed further by decomposing $D'_\lambda(1)$ as a direct sum of submodules.  As $\dim_\amsbb{K}\left(M^\lambda\right)<\infty$, the process will lead after finitely many steps to a decomposition
		\begin{align}
			M=D_\lambda(1)\oplus D_\lambda(2)\oplus\ldots\oplus D_\lambda(k)
		\end{align}
		such that $\lambda$ is an indecomposable weight of each $D_\lambda(i)$.  Thus, $Q_\lambda=\{\lambda\}$ is an indecomposable weight set of each $D_\lambda(i)$.
		
		Now, let $\mu\triangleright\lambda$ be an arbitrary element of $\text{supp}(M)$ such that we have a decomposition of $M$ as desired for every $\nu\triangleleft\mu$.  First, we handle the case where $\mu$ is the immediate successor of a weight $\upsilon$.  Then, by the induction hypothesis, we have a direct sum decomposition
		\begin{align}
			M=\bigoplus_{j\in J_\upsilon}\,D_\upsilon(j)
		\end{align}
		such that $Q_\upsilon$ is an indecomposable weight set of $D_\upsilon(j)$ for all $j\in J_\upsilon$.  We proceed to decompose $D_\upsilon(j)$ with respect to the weight $\mu$ instead of $\lambda$ in the same manner as the base case.  That is, 
		\begin{align}
			D_\upsilon(j)=\bigoplus_{i\in J_{\upsilon}^j}\,D_\upsilon^j(i)
		\end{align}
		for some index set $J_\upsilon^j$ and for some submodules $D_\upsilon^j(i)$ such that $\mu$ is an indecomposable weight.  Therefore,
		\begin{align}
			M=\bigoplus_{j\in J_\upsilon}\,\bigoplus_{i\in J_\upsilon^j}\,D_\upsilon^j(i)
		\end{align}
		is a decomposition in which $Q_\upsilon$ is an indecomposable weight set of each $D_\upsilon^j(i)$ and $\mu$ is an indecomposable weight of each $D_\upsilon^j(i)$.  Note that $Q_\mu$ may not be an indecomposable weight set of each $D_\upsilon^j(i)$, but when that is the case, we can further decompose $D_\upsilon^j(i)$ as follows:
		\begin{align}
			D_\upsilon^j(i)= \bar{D}_\upsilon^j(i)\oplus \tilde{D}_\upsilon^j(i)\,,
		\end{align}
		where $\left(\bar{D}_\upsilon^j(i)\right)^\mu=0$ and $\left(\tilde{D}_\upsilon^j(i)\right)^\nu=0$ for all $\nu\in Q_\upsilon$.  Let $\tilde{J}_\upsilon^i$ denote the subset of $J^i_\upsilon$ consisting of $i\in J_\upsilon^i$ such that $Q_\mu$ is not an indecomposable weight set of $D_\upsilon^j(i)$.  Then,
		\begin{align}
			M=\bigoplus_{j\in J_\upsilon}\,\left(\left(\bigoplus_{i\in J_\upsilon^j\smallsetminus\tilde{J}_\upsilon^j}\,D_\upsilon^j(i) \right)\oplus\left(\bigoplus_{i\in \tilde{J}_\upsilon^j}\,\left(\bar{D}_\upsilon^j(i)\oplus\tilde{D}_\upsilon^j(i)\right)\right)_{\vphantom{1_1}}^{\vphantom{1^1}}\right)
			\label{eq:hugedirectsum}
		\end{align}
		is a direct sum decomposition of $M$ with respect to weight $\mu$ and with the required properties.
		
		Now, suppose that $\mu$ is a limit element of $\text{supp}(M)$.  Then, let $\pazocal{P}_\mu$ be the set of all pairs of the form $\big(\nu,D_\nu(j)\big)$, where $\nu\triangleleft\mu$ and $j\in J_\nu$.  We equip $\pazocal{P}_\mu$ with a partial order $\leq $ by setting
		\begin{align}
			\big(\nu,D_\nu(j)\big) \leq \big(\tilde{\nu},D_{\tilde{\nu}}(\tilde{j})\big)\text{ if and only  if } \nu\trianglelefteq \tilde{\nu}\text{ and }D_\nu(j)\supseteq D_{\tilde{\nu}}(\tilde{j})\,.
			\label{eq:pairspartial}
		\end{align}  
		For each maximal chain $\pazocal{C}\subseteq \pazocal{P}_\mu$, write $D(\pazocal{C})$ for the intersection $\bigcap_{\big(\nu,D_{\nu}\left(j_\nu\right)\big)\in\pazocal{C}}\,D_\nu\left(j_\nu\right)$.  (The existence of maximal chains is guaranteed by Corollary~\ref{cor:hausdorff}.)  Let $\pazocal{M}_\mu$ be the set of all maximal chains of $\pazocal{P}_\mu$.
		
		We need to show that the sum $\sum_{\pazocal{C}\in\pazocal{M}_\mu}\,D(\pazocal{C})$ is direct and equals $M$.
		Let $\xi$ be an arbitrary weight of $M$.  We consider the following decomposition of $M^\xi$ as
		\begin{align}
			M^\xi=\bigoplus_{j\in J_\nu} M^\xi(\nu,j)\,,
			\label{eq:wspdec}
		\end{align}
		where $M^\xi(\nu,j) := M^\xi\cap D_\nu(j)$ for each $\nu \triangleleft \mu$ and $j\in J_\nu$.  Because $M^\xi$ is finite-dimensional, there are finitely many nonzero weight spaces $M^\xi(\nu,j)$ for each $\nu$.  Furthermore, each $M^\xi(\nu,j)$ is a subspace of some $M^\xi\left({\tilde{\nu}},\tilde{j}\right)$ for every $\tilde{\nu}\triangleleft \nu$.  Thus, for some $\upsilon \triangleleft \mu$ depending on $\nu$, the decomposition (\ref{eq:wspdec}) of the weight space $M^\xi$ stabilizes at every level $\nu$ with $\upsilon\trianglelefteq \nu\triangleleft\mu$.  This  implies that
		\begin{align}
			M^\xi =\bigoplus_{\pazocal{C}\in\pazocal{M}_\mu}\,\big(D(\pazocal{C})\big)^\xi
		\end{align}
		for every $\xi\in\text{supp}(M)$, whence we have
		\begin{align}
			M=\bigoplus_{\pazocal{C}\in\pazocal{M}_\mu}\,D(\pazocal{C})\,.
		\end{align}
		It may be the case that $\big(D(\pazocal{C})\big)^\nu=0$ for some $\nu\triangleleft\mu$ (in particular, $D(\pazocal{C})$ may be a zero module), but in any case, every $\nu\triangleleft\mu$ is an indecomposable weight of $D(\pazocal{C})$.  Furthermore, since each $Q_\nu$ is an indecomposable weight set of $D_\nu\left(j_\nu\right)$, we conclude that 
		\begin{align}
			Q_\mu':=\bigcup_{\nu\triangleleft\mu}\,Q_\nu=\big\{\nu\in\text{supp}(M)\suchthat{}\nu\triangleleft\mu\big\}
		\end{align} is an indecomposable weight set of each $D(\pazocal{C})$.
		
		We apply the same argument from the base case to each $D(\pazocal{C})$, using $\mu$ in place of $\lambda$ and write
		\begin{align}
			D(\pazocal{C})=\bigoplus_{j\in J(\pazocal{C})}\,D(\pazocal{C},j)
		\end{align}
		for some index set $J(\pazocal{C})$ and for some direct summands $D(\pazocal{C},j)$ fo $D(\pazocal{C})$ such that $\mu$ is an indecomposable weight of $D(\pazocal{C},j)$.  Thus, we have a required direct sum decomposition of $M$ with respect to the weight $\mu$:
		\begin{align}
			M=\bigoplus_{\pazocal{C}}\,\bigoplus_{j\in J(\pazocal{C})}\,D(\pazocal{C},j)\,.
		\end{align}
		We then utilize the same argument that leads to (\ref{eq:hugedirectsum}) in the immediate successor case to show that there exists a subset $\tilde{J}(\pazocal{C})$ of $J(\pazocal{C})$ such that
		\begin{enumerate}[(i)] 
		\item $D(\pazocal{C},j)=\bar{D}(\pazocal{C},j)\oplus\tilde{D}(\pazocal{C},j)$, where $\big(\bar{D}(\pazocal{C},j)\big)^\mu=0$ and $\big(\tilde{D}(\pazocal{C},j)\big)^\nu=0$ for all $\nu\in Q_\mu'$,
		\item if $j\in J(\pazocal{C})\smallsetminus \tilde{J}(\pazocal{C})$, then $Q_\mu=\big\{\nu\in\text{supp}(M)\suchthat{}\nu\trianglelefteq\mu\big\}$ is an indecomposable weight set of $D(\pazocal{C},j)$.
		\end{enumerate}Then, we have achieved a direct sum decomposition 
		\begin{align}
			M=\bigoplus_{\pazocal{C}}\,\left(\left(\bigoplus_{j\in J(\pazocal{C})\smallsetminus\tilde{J}(\pazocal{C})}\,D(\pazocal{C},j) \right)\oplus\left( \bigoplus_{j\in \tilde{J}(\pazocal{C})}\,\big(\bar{D}(\pazocal{C},j)\oplus \tilde{D}(\pazocal{C},j)\big)\right)_{\vphantom{1_1}}^{\vphantom{1^1}}\right)
			\label{eq:hugedirect}
		\end{align}
		of $M$ satisfying the condition that the set $Q_\mu$ is an indecomposable weight set of each $D(\pazocal{C},j)$, and from the construction of (\ref{eq:hugedirect}), we know that, for a given direct summand  $D(\pazocal{C},j)$, $\bar{D}(\pazocal{C},j)$, or $\tilde{D}(\pazocal{C},j)$ of (\ref{eq:hugedirect}) and  for every $\nu\triangleleft\mu$, there exists a unique $i\in J_\nu$ such that this direct summand is a submodule of $D_\nu(i)$.
		
		To complete the proof, we define the partially ordered set $\mathscr{P}$ to be the set of all pairs $\big(\mu,D_\mu(j)\big)$ where $\mu\in\text{supp}(M)$ and $j\in J_\mu$, equipped with the same partial order $\leq$ defined by (\ref{eq:pairspartial}).  We write $\mathscr{M}$ for the set of maximal chains in $\mathscr{P}$.   We apply a similar argument as in the previous paragraph, using the finite-dimensionality assumption on the weight spaces of $M$, to establish that
		\begin{align}
			M=\bigoplus_{\mathscr{C}\,\in\,\mathscr{M}}\,\mathcal{D}(\mathcal{C})\,,
		\end{align}
		and 
		\begin{align}
			\mathcal{D}(\mathcal{C}):=\bigcap_{\big(\mu,D_\mu\left(j_\mu\right)\big)\in\mathcal{C}}\,D_\mu\left(j_\mu\right)
		\end{align} 
		for every maximal chain $\mathcal{C}$ in $\mathscr{P}$.  For each $\mu\in\text{supp}(M)$, the set $Q_\mu$ is an indecomposable weight set of each $D_\mu\left(j_\mu\right)$.  Then, $\text{supp}(M)=\bigcup_{\mu\in\text{supp}(M)}\,Q_\mu$ is an indecomposable weight set of every $\mathcal{D}(\mathcal{C})$.  Consequently, $\mathcal{D}(\mathcal{C})$ is an indecomposable module.  (Note that a direct summand $\mathcal{D}(\mathcal{C})$ for some $\mathcal{C}$ may be the zero module, but this does not effect the proof or the statement of this theorem.)
	\end{pf}
\end{thm}

\begin{rmk}
	The proof of Theorem~\ref{thm:directsumdecomposition} does not explicitly use the fact that $M$ is a weight module.  In general, let $M$ be a module over a (not necessarily associative) $\amsbb{K}$-algebra $\pazocal{A}$.  Suppose that $M$ has a vector space decomposition 
	\begin{align}
		M=\bigoplus_{i\in I}\,M^i\,,
	\end{align}
	where $M^i$ is a finite-dimensional $\amsbb{K}$-vector subspace of $M$ and $I$ is an index set, and this decomposition has the property that, for every $\pazocal{A}$-submodule $N$ of $M$, we have
	\begin{align}	
		N=\bigoplus_{i \in I}\,\big(N\cap M^i\big)
	\end{align}
	as a direct sum of vector spaces.  Then, $M$ can be written as a direct sum of indecomposable $\pazocal{A}$-submodules.
\end{rmk}

\begin{prop}
	Let $M\in\bar\bggO$ be indecomposable.  Then, every $\varphi\in \End_{\bar\bggO}(M)$ is either an automorphism or is \emph{locally nilpotent} (namely, for every $v\in M$, there exists $k\in\amsbb{Z}_{\geq 0}$ such that $\varphi^{ k}(v)=0$).
	
	\begin{pf}
		Let $K_k:=\ker\left(\varphi^{k}\right)$ and $I_k:=\im\left(\varphi^{k}\right)$ for each $k=0,1,2,\ldots$ (here, $\varphi^0$ is the identity map $\text{id}_M$).  We observe that $K_0\subseteq K_1\subseteq K_2 \subseteq \ldots$ and $I_0\supseteq I_1\supseteq I_2\supseteq \ldots$.  Set $K:=\bigcup_{k=0}^\infty\,K_k$ and $I:=\bigcap_{k=0}^\infty\,I_k$.  
		
		Fix $\lambda\in\hfrak^*$.  The restriction $\psi_\lambda$ of $\varphi$ onto $M^\lambda$ is a linear map on a finite-dimensional vector space.  Hence, $M^\lambda$ decomposes as $\im\left(\psi_\lambda^{ k}\right)\oplus \ker\left(\psi_\lambda^{k}\right)$ for every $k=0,1,2,\ldots$.  Since $M^\lambda$ is a finite-dimensional vector space and 
		\begin{align}
			\im\left(\psi_\lambda\right)\supseteq \im\left(\psi_\lambda^{2}\right)\supseteq \im\left(\psi_\lambda^{3}\right)\supseteq \ldots\,,
	\end{align} the submodules $\im\left(\psi_\lambda^{k}\right)$, where $k=0,1,2,\ldots$, must stabilize.  Assume that, for some $j\in\amsbb{Z}_{\geq 0}$, we have 
	\begin{align}\im\left(\psi_\lambda^{ j}\right)=\im\left(\psi_\lambda^{(j+1)}\right)=\im\left(\psi_\lambda^{(j+2)}\right)=\ldots\,.
	\end{align}  That is, the kernels must also stabilize at the same index: 
	\begin{align}
	\ker\left(\psi_\lambda^{ j}\right)=\ker\left(\psi_\lambda^{(j+1)}\right)=\ker\left(\psi_\lambda^{(j+2)}\right)=\ldots\,.
	\end{align}  This shows that $K^\lambda$ and $I^\lambda$ are equal to $K_j^\lambda$ and $I_j^\lambda$ for some $j\in\amsbb{Z}_{\geq 0}$, depending on $\lambda$.  Therefore, the sum 
	\begin{align}
		(K+L)^\lambda=K^\lambda+L^\lambda =K_j^\lambda+I_j^\lambda
	\end{align} is direct and equals $M^\lambda$, as $M^\lambda =\im\left(\psi_\lambda^{ j}\right)\oplus\ker\left(\psi_\lambda^{j}\right)=I_j^\lambda\oplus K_j^\lambda$.  Since $\lambda$ is arbitrary, we obtain $M=K\oplus I$.
		
		As $M$ is an indecomposable object, we have either $K=0$ or $I=0$.  In the former case, we conclude that $\varphi$ is an isomorphism, and in the latter case, we see that $\varphi$ is locally nilpotent.
	\end{pf}
\end{prop}

\begin{prop}
	For every indecomposable object $M\in\bar\bggO$, the endomorphism ring $\End_{\bar\bggO}{(M)}$ is a local ring.
	
	\begin{pf}
		Let $J\subseteq R:=\End_{\bar\bggO}(M)$ be the set of all locally nilpotent endomorphisms of $M$.  By the previous proposition, $J$ is the set of all non-invertible elements of $R$.  We must prove that $J$ is an ideal of $R$.  
		
		First, if $\varphi\in J$ and $\psi\in R$, then $\varphi\circ\psi$ cannot be an epimorphism because $\varphi$ is not surjective on any weight space of $M$, and $\psi\circ\varphi$ is not a monomorphism because $\varphi$ is not injective on any weight space of $M$.  That is, $\varphi\circ\psi$ and $\psi\circ\varphi$ are both in $J$.  
		
		Now, we assume that $\varphi,\psi\in J$.  We must show that $\varphi+\psi$ belongs to $J$ too.  Suppose on the contrary that $\varphi+\psi\notin J$.  Then, $\varphi+\psi$ is invertible.  Hence, $\varphi+\psi=\phi$ for some automorphism $\phi:M\to M$.  Let $\alpha:=\varphi\circ\phi^{-1}$ and $\beta:=\psi\circ\phi^{-1}$.  Then,  $\alpha+\beta=\id_M$ and $\alpha,\beta\in J$.  Note that $\alpha\circ\beta=\beta\circ\alpha$.  Fix a weight $\lambda$ of $M$.  Suppose that $\alpha^{k}$ and $\beta^{l}$ vanish on $M^\lambda$, for some $k,l\in\amsbb{Z}_{> 0}$.  Then, 
		\begin{align}
			(\alpha+\beta)^{ (k+l)}=\sum_{r=0}^{k+l}\,\binom{k+l}{r}\,\alpha^{ r}\circ \beta^{(k+l-r)}
			\end{align} must vanish on $M^\lambda$ as well.  Ergo, the endomorphism $\alpha+\beta$ cannot equal $\id_M$, which is a contradiction.
	\end{pf}
\end{prop}

Next, we recall the following important theorem.  The proof is omitted here, but can be found in~\cite{KSRA}.

\begin{thm}[\bf Krull-Schmidt-Remak-Azuyama]
	Let $R$ be a unital ring and $M$ a unitary left $R$-module.  Suppose that $M$ is a direct sum of modules whose endomorphism rings are local rings.  Then, any two (not necessarily finite) direct sum decompositions of $M$ into indecomposable direct summands are isomorphic.
\end{thm}

The Krull-Schmidt-Remak-Azuyama Theorem immediately implies the following corollary.

\begin{cor}
	Every object in $\bar\bggO$ admits a unique, up to isomorphism, decomposition into a direct sum of indecomposable objects.
	\label{cor:KSRA}
\end{cor}

\begin{prop}
	Every indecomposable object $M\in\bar\bggO$ satisfies $\text{supp}(M)\subseteq \lambda+\Lambda$, where $\Lambda$ is the root lattice.  In particular, $\text{supp}(M)$ is countable and $M$ is countable dimensional.
	
	\label{prop:indecobj_weight}
	
	\begin{pf}
		For each equivalence class $[\lambda]\in\hfrak^*/\Lambda$, we define $M^{[\lambda]}$ to be the submodule of $M$ consisting of weight vectors whose weights lie in $[\lambda]$.  It is trivial that $M=\bigoplus_{[\lambda]\in\hfrak^*/\Lambda}\,M^{[\lambda]}$.  Being indecomposable, $M=M^{[\lambda]}$ for some $\lambda\in\hfrak^*$.
	\end{pf}
\end{prop}
	
\subsection{Generalized Composition Series}

In this subsection, $\bfrak$ is assumed to be a Dynkin Borel subalgebra, unless otherwise specified.  Following the approach of V. Kac in~\cite{kac}, we have the following theorem.

\begin{thm}
	Let $M\in\bar\bggO$ and $\lambda\in\hfrak^*$.  Suppose that all weights $\xi$ of $M$ satisfy $\xi\preceq \upsilon$ for some fixed upper bound $\upsilon\in\hfrak^*$.  Then, there exist a $\gfrak$-module filtration 
	\begin{align}
		0=M_0\subseteq M_1\subseteq M_2\subseteq \ldots \subseteq M_{k-1}\subseteq M_k=M
	\end{align} 
	and a subset $J\subseteq \{1,2,\ldots,k\}$ such that
	\begin{enumerate}[(i)]
		\item if $j\in J$, then $M_{j}/M_{j-1}\cong \Llie\big(\xi(j)\big)$ for some $\xi(j)\in \hfrak^*$ with $\xi(j)\succeq \lambda$,
		\item if $j\notin J$, then $\left(M_{j}/M_{j-1}\right)^\mu=0$ for every $\mu\succeq \lambda$.
	\end{enumerate}
	\label{thm:kac_series}
	
	\begin{pf}
		For any $\lambda\in\hfrak^*$, set $d(M,\lambda):=\sum_{\mu\succeq \lambda}\,\dim_\amsbb{K}\left(M^\mu\right)$ and note that $d(M,\lambda)<\infty$.  This follows from to the facts that $M$ has an upper bound $\upsilon$ and that $\bfrak$ is Dynkin.  We shall prove the theorem by induction on $d(M,\lambda)$.  
		
		The base case $d(M,\lambda)=0$ is done by considering the filtration $0=M_0\subseteq M_1=M$, with $J=\emptyset$.  Now, suppose that $d(M,\lambda)>0$. Choose a singular weight $\tilde{\lambda}$ of $M$ with $\tilde{\lambda}\succeq \lambda$.   Let $v\in M^{\tilde{\lambda}}$ be nonzero.  Set $N:=\Ulie(\gfrak)\cdot v$.  Take $\tilde{N}$ to be the maximal proper submodule of $N$.  Then we have
	\begin{align}
		0\subseteq \tilde{N}\subseteq N\subseteq M
	\end{align}
	with $N/\tilde{N}\cong \Llie\left(\tilde{\lambda}\right)$ and $\tilde{\lambda}\succeq \lambda$.  Since $d\left(\tilde{N},\lambda\right)<d(M,\lambda)$ and $d\left(M/N,\lambda\right)<d(M,\lambda)$, we can apply the induction hypothesis to get a filtration of $M$ as required.
	\end{pf}
\end{thm}

The previous theorem raises a question whether every object $M$ in $\bbggO$ satisfies \begin{equation}
	d(M,\lambda)=\sum_{\mu\preceq \lambda}\,\dim_\amsbb{K}\left(M^\mu\right) <\infty\,.
\end{equation}  In addition, does every indecomposable object $M$ in $\bbggO$ have an upper bound $\upsilon\in\hfrak^*$?  However, these are still open questions.

\begin{cor}
	Let $M\in\bar\bggO$ and $\lambda,\nu\in\hfrak^*$ with $\lambda\preceq \nu$.    Then there exist a $\gfrak$-module filtration $0=M_0\subseteq M_1\subseteq M_2\subseteq \ldots \subseteq M_{k-1}\subseteq M_k =M$ and a subset $J\subseteq \{1,2,\ldots,k\}$ such that
	\begin{enumerate}[(i)]
		\item if $j\in J$, then $M_j/M_{j-1}\cong \Llie\big(\xi(j)\big)$ for some $\xi(j)\in\hfrak^*$ with $\lambda\preceq \xi(j)\preceq \nu$,
		\item if $j\notin J$, then  either $\left(M_j/M_{j-1}\right)^\mu=0$ for every $\mu\in\hfrak^*$ satisfying $\lambda\preceq\mu\preceq\nu$, or $M_j/M_{j-1}\cong \Llie\big(\xi(j)\big)$ for some $\xi(j)\succ\lambda$ such that $\xi(j)\not\preceq \nu$.
	\end{enumerate}	Such a filtration is called a \emph{composition series of $M$ with bounds $\lambda$ and $\nu$}.  The set $J$ is called the \emph{relevant index set} of such a filtration.
	\label{cor:kac_series}
	
	\begin{pf}
		Since the interval $[\lambda,\nu]:=\big\{\zeta\in\hfrak^*\,\suchthat{\lambda\preceq\zeta\preceq \nu}\big\}$ is finite (as $\bfrak$ is Dynkin) and $M$ is $\nfrak$-locally finite, the submodule
		\begin{align}
			\tilde{M}:=\sum_{\substack{{\mu\in\hfrak^*}\\{\lambda\preceq \mu\preceq\nu}}}\,\Ulie(\gfrak)\cdot M^\mu =\sum_{\substack{{\mu\in\hfrak^*}\\{\lambda\preceq \mu\preceq\nu}}}\,\Ulie(\nfrak)\cdot M^\mu
			\label{eq:mtilde}
		\end{align}
		has finitely many weights $\zeta$ with $\zeta\succeq \lambda$.  Therefore, $\tilde{M}$ has an upper bound $\upsilon\in\hfrak^*$.  We apply Theorem~\ref{thm:kac_series} on $\tilde{M}$ and obtain a filtration
		\begin{align}
			0=M_0\subseteq M_1\subseteq M_2\subseteq \ldots\subseteq M_{k-2}\subseteq M_{k-1}=\tilde{M}
		\end{align}
		along with a subset $\tilde{J}\subseteq\{1,2,\ldots,k-1\}$ satisfying the condition that, if $j\in\tilde{J}$, then $M_j/M_{j-1}\cong \Llie\big(\xi(j)\big)$ for some $\xi(j)\in\hfrak^*$ with $\xi(j)\succeq \lambda$, and if $j\notin\tilde{J}$, then $\left(M_j/M_{j-1}\right)^\mu=0$ for every $\mu\succeq\lambda$.  Then, by setting $M_{k}:=M$, we have the filtration 
		\begin{align}
			0=M_0\subseteq M_1\subseteq M_2\subseteq \ldots\subseteq M_{k-2}\subseteq M_{k-1}\subseteq M_{k}=M\,.
			\label{eq:filtrationM}
		\end{align}
		Let $J:=\big\{j\in\tilde{J}\,|\,\xi(j)\preceq \nu\big\}$.  The filtration (\ref{eq:filtrationM}) clearly satisfies (i) and (ii), with the relevant index set $J$, noting that $\left(M_k/M_{k-1}\right)^\mu=\left(M/\tilde{M}\right)^\mu=0$ for all $\mu\in\hfrak^*$ such that $\lambda\preceq\mu\preceq \nu$ holds.
	\end{pf}
\end{cor}

\begin{define}
	Let $\lambda,\nu\in\hfrak^*$ satisfy $\lambda\preceq \nu$.  Suppose that 
	\begin{align} 0=M_0\subseteq M_1\subseteq M_2\subseteq \ldots \subseteq M_{k-1}\subseteq M_k=M\label{eq:com1}\end{align} and \begin{align} 0=M'_0\subseteq M'_1\subseteq M'_2\subseteq \ldots \subseteq M'_{k'-1}\subseteq M'_{k'}=M\label{eq:com2}\end{align} are two composition series of $M\in\bar\bggO$ with bounds $\lambda$ and $\nu$, and with relevant index sets $J$ and $J'$, respectively.  We say that these filtrations are \emph{equivalent} if there exists a bijection $f:J\to J'$ such that $M_{j}/M_{j-1}\cong M'_{f(j)}/M'_{f(j)-1}$ for all $j\in J$.
\end{define}

\begin{lem}
	Let $\lambda,\nu\in\hfrak^*$ be such that $\lambda\preceq \nu$.  Denote by $\tilde{M}$ the submodule of $M$ given by (\ref{eq:mtilde}).  Suppose that $0=M_0\subseteq M_1 \subseteq M_2 \subseteq \ldots \subseteq M_{k-1}\subseteq M_k=M$.  Define $\tilde{M}_j:=M_j\cap\tilde{M}$ for every $j=0,1,2,\ldots,k$.  Then, $0=\tilde{M}_0\subseteq \tilde{M}_1\subseteq \tilde{M}_2\subseteq \ldots\subseteq \tilde{M}_{k-1}\subseteq \tilde{M}_k=\tilde{M}$ is a composition series of $\tilde{M}$ with bounds $\lambda$ and $\nu$.
	
	\begin{pf}
		For each $j=1,2,\ldots,k$, define $\varphi_j:\tilde{M}_j/\tilde{M}_{j-1}\to M_j/M_{j-1}$ via $v+\tilde{M}_{j-1}\mapsto v+M_{j-1}$ for every $v\in\tilde{M}_j$.  Clearly, $\varphi_j$ is well defined, and it is injective because 
		\begin{align} 	\tilde{M}_{j}\cap M_{j-1}=\tilde{M}\cap M_j\cap M_{j-1}=\tilde{M}\cap M_{j-1}=\tilde{M}_{j-1}\,.\end{align}
		If $M_j/M_{j-1}\cong \Llie(\mu)$ for some $\mu\in\hfrak^*$ with $\mu\preceq\lambda$, then 
		\begin{align} \dim_\amsbb{K}\left(\tilde{M}_j^\mu\right)=\dim_\amsbb{K}\left(M_j^\mu\right)=\dim_\amsbb{K}\left(M_{j-1}^\mu\right)+1=\dim_\amsbb{K}\left(\tilde{M}_{j-1}^\mu\right)+1\,.\end{align}
		Hence, $\dim_\amsbb{K}\left(\left(\tilde{M}_{j}/\tilde{M}_{j-1}\right)^\mu\right)=1$, so $\varphi_j$ is nonzero.  As $M_j/M_{j-1}$ is simple, $\varphi_j$ must be surjective, whence it gives an isomorphism $\tilde{M}_j/\tilde{M}_{j-1} \cong M_j/M_{j-1}\cong \Llie(\mu)$.
		
		Let $J$ be the relevant index set of the composition series $M_0\subseteq M_1 \subseteq \ldots \subseteq M_{k-1}\subseteq M_k$ of $M$ with bounds $\lambda$ and $\nu$.  By the observation above, if $j\in J$, then $\tilde{M}_j/\tilde{M}_{j-1}\cong M_j/M_{j-1}$ is simple with highest weight $\mu$ satisfying $\lambda\preceq\mu\preceq\nu$. If an index $j\in\{1,2,\ldots,k\}\smallsetminus J$ is such that $M_j/M_{j-1}$ is a simple module with highest weight $\xi\succ \lambda$ with $\xi\not\preceq\lambda$, then using the embedding  $\varphi_j: \tilde{M}_j/\tilde{M}_{j-1}\to M_j/M_{j-1}$, we  conclude that either $\tilde{M}_j/\tilde{M}_{j-1}\cong\Llie(\xi)$ or $\tilde{M}_j/\tilde{M}_{j-1}=0$.  Finally, if $j\in\{1,2,\ldots,k\}\smallsetminus J$ is such that $\left(M_j/M_{j-1}\right)^\mu =0$ for every $\mu$ with $\lambda\preceq\mu\preceq\nu$, using the embedding $\varphi_j: \tilde{M}_j/\tilde{M}_{j-1}\to M_j/M_{j-1}$, we see that $\left(\tilde{M}_j/\tilde{M}_{j-1}\right)^\mu =0$ for every $\mu\in\hfrak^*$ with $\lambda\preceq\mu\preceq\nu$.
	\end{pf}
\end{lem}

\begin{thm}
	Let $\lambda,\nu\in\hfrak^*$ be such that $\lambda\preceq \nu$.  Then, any two composition series of $M\in\bar\bggO$ with bounds $\lambda$ and $\nu$ are equivalent.
	\label{thm:kacequiv}
	
	\begin{pf}
		Suppose that (\ref{eq:com1}) and (\ref{eq:com2}) are two composition series of $M$ with bounds $\lambda$ and $\nu$.  Let $\tilde{M}$ be the submodule of $M$ defined by (\ref{eq:mtilde}).  From the lemma above, it suffices to assume that $M=\tilde{M}$.  
		
		From the assumption $M=\tilde{M}$, there are finitely many weights $\mu$ of $M$ satisfying $\mu\preceq\lambda$.  Thus, we can refine (\ref{eq:com1}) and (\ref{eq:com2}) in the same manner as in Theorem~\ref{thm:kac_series} to get index sets $\tilde{J}\subseteq\{1,2,\ldots,k\}$ and $\tilde{J}'\subseteq\{1,2,\ldots,k'\}$ such that the following three conditions are met:
		\begin{enumerate}[(i)]
			\item $J\subseteq \tilde{J}$ and $J'\subseteq \tilde{J}'$,
			\item for $j\in \tilde{J}$ and $j'\in \tilde{J}'$, $M_j/M_{j-1}$ and $M'_{j'}/M'_{j'-1}$ are simple modules with highest weights greater than or equal to $\lambda$,
			\item for $j\notin\tilde{J}$ and $j'\notin\tilde{J}'$, all the weight spaces $\left(M_j/M_{j-1}\right)^\mu$ and $\left(M'_{j'}/M'_{j'-1}\right)^\mu$ with $\mu\succeq\lambda$ are the zero vector space.
		\end{enumerate}
		As a result, we can instead show that there exists a bijection $\tilde{f}:\tilde{J}\to\tilde{J}'$ such that $M_j/M_{j-1}\cong M'_{f(j)}/M'_{f(j)-1}$ for all $j\in\tilde{J}$.  The restriction $f:=\tilde{f}|_J$ yields a bijection $f:J\to J'$ as required.
		
		For each $\xi\in\hfrak^*$ with $\lambda \preceq \xi \preceq \nu$, let $t(\xi)$ denote the maximum possible value of the length of the positive root $\mu-\xi$, where $\mu\succeq\xi$ is a singular weight of $M$.  For each $l=0,1,2,\ldots$, write $T_l$ for the set $\big\{\xi\in\hfrak^*\suchthat{}\xi\succeq\lambda\text{ and }t(\xi)=l\big\}$.  We shall instead prove that, for a fixed $l=0,1,2,\ldots$, the number of $j\in J$ with $M_{j}/M_{j-1}\cong \Llie(\xi)$ is the same as the number of $j'\in J'$ with $M'_{j'}/M'_{j'-1}\cong \Llie(\xi)$ for every $\xi \in T_l$.
		
		The proof goes by induction on $l$.  For the base case $l=0$, every $\xi\in T_l$ is a singular weight of $M$, whence the weight space $M^\xi$ comes from  $\dim_\amsbb{K}\left(M^\xi\right)$ copies of $\Llie(\xi)$ in any composition series with  bounds $\lambda$ and $\nu$.   
		
		Let now assume that $l>0$ and $\xi\in T_l$.  By the induction hypothesis, the multiplicities of each factor $\Llie(\tilde{\xi})$ with $\tilde{\xi}\in T_0\cup T_1 \cup\ldots \cup T_{l-1}$ in the filtrations (\ref{eq:com1}) and (\ref{eq:com2}) are equal, and let $m\left(\tilde{\xi}\right)$ denote the common value.  For each $j\in J$ with $M_j/M_{j-1} \not\cong \Llie\left(\tilde{\xi}\right)$ with $\tilde{\xi}\in T_0\cup T_1\cup\ldots\cup T_{l-1}$, we observe that either $\left(M_j/M_{j-1}\right)^{\xi}=0$, or $\xi$ is a singular weight of $M_j/M_{j-1}$ (making $M_j/M_{j-1}\cong \Llie(\xi)$).  Hence, there are exactly 
		\begin{align}
		m(\xi):= \dim_\amsbb{K}\left(M^\xi\right)-\sum_{r=0}^{l-1}\,\sum_{\tilde{\xi}\in T_r}\,m(\tilde{\xi})\,\dim_\amsbb{K}\Big(\big(\Llie(\tilde{\xi})\big)^\xi\Big)
		\label{eq:multcount}
		\end{align} values of such $j\in J$ with $M_j/M_{j-1}\cong \Llie(\xi)$.  Therefore, $m(\xi)$ is the multiplicity of $\Llie(\xi)$ in (\ref{eq:com1}).  Since the value $m(\xi)$ as shown in (\ref{eq:multcount}) depends only on the previously known values $m\left(\tilde{\xi}\right)$ with $\tilde{\xi}\in T_0\cup T_1 \cup\ldots \cup T_{l-1}$, $m(\xi)$ is also the multiplicity of $\Llie(\xi)$ in (\ref{eq:com2}).  The induction is now complete and the claim follows.
	\end{pf}
\end{thm}

Corollary~\ref{cor:kac_series} and Theorem~\ref{thm:kacequiv} form a partial extension of the usual Jordan-H\"{o}lder Theorem for modules of finite length.  Based on this, we now extend the usual definition of composition factors and composition factor multiplicities as follows.

\begin{cor}
	Let $M\in\bar\bggO$ and $\mu\in\hfrak^*$ be such that $\lambda\preceq\mu\preceq \nu$.  The number of times $\Llie(\mu)$ occurs as a factor in any composition series of $M$ with bounds $\lambda$ and $\nu$ is independent of the choice of the composition series with bounds and the choice of the bounds $\lambda,\nu\in\hfrak^*$ (as long as $\lambda\preceq \mu \preceq \nu$).    This number is known as the \emph{composition factor multiplicity} of $\Llie(\mu)$ in $M$, and is denoted by $\big[M:\Llie(\mu)\big]$.  If $\big[M:\Llie(\mu)\big]>0$, then we say that $\Llie(\mu)$ is a \emph{composition factor} of $M$.
	\label{cor:compositionfactormult}
	
	\begin{pf}
		For given weights $\lambda,\nu\in\hfrak^*$, Theorem~\ref{thm:kacequiv} guarantees that the number $m^M_\mu(\lambda,\nu)$  of times $\Llie(\mu)$ occurs as a factor does not depend on the choice of the composition series of $M$ with bounds $\lambda$ and $\nu$.  We have to show that $m^M_\mu(\lambda,\nu)$ is also independent of $\lambda$ and $\nu$, provided that $\lambda\preceq\mu\preceq\nu$.
		
		Let $\lambda,\nu\in\hfrak^*$ be such that $\lambda\preceq\mu\preceq\nu$.  We choose an arbitrary composition series $0=M_0\subsetneq M_1\subsetneq M_2\subsetneq \ldots \subsetneq M_k=M$ of $M$ with bounds $\lambda$ and $\nu$.  Then, this filtration is also a composition series with bounds $\mu$ and $\mu$.  Again, by Theorem~\ref{thm:kacequiv}, this filtration is equivalent to any composition series with bounds $\mu$ and $\mu$, which immediately implies that $m^M_\mu(\lambda,\nu)=m^M_\mu(\mu,\mu)$.
	\end{pf}
\end{cor}

\begin{rmk}
	From the proof above, there are at most two possible values of $m^M_\mu(\lambda,\nu)$.  If the condition $\lambda\preceq\mu\preceq \nu$ is not met, then $m_\mu^M(\lambda,\nu)=0$.  Otherwise, $m_\mu^M(\lambda,\nu)=m_\mu^M(\mu,\mu)$.  The composition factor multiplicity can be simply defined as $\big[M:\Llie(\mu)\big]:=m_\mu^M(\mu,\mu)$.  
	
	From this observation, it is possible to make yet another generalization of the usual composition series.  We can fix a finite subset $S\subseteq\hfrak^*$, and then create a \emph{composition series of $M\in\bbggO$ with weight reference set $S$}, that is, a filtration  
	\begin{align}
		0=M_0\subseteq M_1\subseteq M_2\subseteq \ldots \subseteq M_{k-1}\subseteq M_k=M
	\end{align} and a subset $J\subseteq \{1,2,\ldots,k\}$ such that
	\begin{enumerate}[(i)]
		\item if $j\in J$, then $M_{j}/M_{j-1}\cong \Llie\big(\xi(j)\big)$ for some $\xi(j)\in S$,
		\item if $j\notin J$, then either $\left(M_{j}/M_{j-1}\right)^\mu=0$ for every  $\mu\in S$, or $M_j/M_{j-1}\cong\Llie\big(\xi(j)\big)$ for some $\xi(j) \notin S$ such that $\xi(j)\succ \lambda$ for some $\lambda\in S$.
	\end{enumerate}
Composition series with a weight reference set are useful especially when the splitting Borel subalgebras are not Dynkin.  However, as our focus lies with Dynkin Borel subalgebras, properties of composition series with a weight reference set will not be discussed in detail, except that there is a notion of equivalence of such filtrations, and therefore, a composition factor multiplicity $\big[M:\Llie(\mu)\big]$ can also be defined.  The number $m_\mu^M(\lambda,\nu)$ in the proof of Corollary \ref{cor:compositionfactormult} is replaced by the number $m_\mu^M(S)$ of times $\Llie(\mu)$ shows up as a factor in a given composition series with weight reference set $S$.  Then, we can set $\big[M:\Llie(\mu)\big]:=m_\mu^M\big(\{\mu\}\big)$.
\end{rmk}

We return to the case where $\bfrak$ is a Dynkin Borel subalgebra.  Now we shall use the composition series with bounds to study generalized composition series, as introduced below.  With the restriction that the modules in $\bbggO$ have finite-dimensional weight spaces, we shall see that these generalized composition series behave similarly to the composition series of modules of finite length. 

\begin{define}
	A \emph{generalized composition series} of a module $M\in\bbggO$ is a family of submodules $\left(M_j\right)_{j\in \mathcal{J}}$ satisfying the following conditions:
	\begin{enumerate}[(i)]
		\item the index set $\mathcal{J}$ is equipped with a total order $<$,
		\item $\bigcap_{j\in\mathcal{J}}\,M_{j}=0$ and $\bigcup_{j\in\mathcal{J}}\,M_j=M$,
		\item $M_j/M_{<j}$ is a simple module for all $j\in \mathcal{J}$, where $M_{<j}:=\bigcup_{k<j}\,M_k$.
	\end{enumerate}
	\label{def:gencompseries}
\end{define}

Note that $\left(M_j\right)_{j\in\mathcal{J}}$ is a generalized composition series of an object $M\in\bbggO$, then the family $\mathcal{G}:=\left(M_j\right)_{j\in\mathcal{J}}\cup\left(M_{<j}\right)_{j\in\mathcal{J}}$ satisfies the property that each module $N\in\mathcal{G}$, $N$ either has an immediate successor or an immediate predecessor (considering $\mathcal{G}$ a totally ordered set with respect to inclusion).  If $\left\{N',N''\right\}$ is a predecessor-successor pair in $\mathcal{G}$, then $N''/N'$ is a simple module.  Therefore, a generalized composition series of $M$ a \emph{generalized flag} on $M$ in the sense of~\cite{genflags}.

\begin{thm}
	Every $M\in\bar\bggO$ has a generalized composition series. 
	\label{thm:gencomp} 
	
	\begin{pf}
		First, we shall prove this theorem when $M$ is indecomposable.  We start with arbitrary weights $\lambda(1)$ and $\nu(1)$ of $M$ with $\lambda(1)\preceq \nu(1)$.  Let $\mathcal{J}(1)$ be the relevant index set of a composition series of $M$ with bounds $\lambda(1)$ and $\nu(1)$.  We create two sequences of weights $\big\{\lambda(k)\big\}_{k\in\mathbb{Z}_{>0}}$ and $\big\{\nu(k)\big\}_{k\in\mathbb{Z}_{>0}}$ in such a way that
		\begin{align}
			\lambda(1)\succ\lambda(2)\succ\lambda(3)\succ\ldots\,,
			\label{eq:lambdakk}
		\end{align}
		\begin{align}
			\nu(1)\prec\nu(2)\prec\nu(3)\prec\ldots\,,
			\label{eq:nukk}
		\end{align}
		and, for every weight $\zeta\in \text{supp}(M)$, there exists $k\in\amsbb{Z}_{>0}$ (depending on $\zeta$) such that $\lambda(k)\preceq \zeta\preceq \nu(k)$.  Note that $\lambda(k)$ and $\nu(k)$ do not have to be weights of $M$.  Therefore, it is always possible to find an infinite strictly decreasing sequence (\ref{eq:lambdakk}) and an infinite strictly increasing sequence (\ref{eq:nukk}).

		Suppose that $\mathcal{J}(k)$ is known.  We extend the filtration in the $k$-th step to obtain a composition series of $M$ with bounds $\lambda(k+1)$ and $\mu(k+1)$.  To be precise, suppose that
		\begin{align}
			0=M^k_0\subsetneq M^k_1\subsetneq M^k_2 \subsetneq \ldots \subsetneq M^k_{l(k)-1}\subsetneq M^k_{l(k)}=M
			\label{eq:lambdak}
		\end{align}
		is a composition series with bounds $\lambda(k)$ and $\mu(k)$.  If $i$ is in the relevant index set $\mathcal{J}(k)$, then $M^k_{i}/M^k_{i-1}$ is simple with highest weight $\mu$ with 
		\begin{align}
			\lambda(k+1)\prec\lambda(k)\preceq \mu\preceq \nu(k)\prec \nu(k+1)\,.
		\end{align}  
		If $i>0$ is not in the relevant index set, then we find a composition series of $M^k_i/M^k_{i-1}$ with bounds $\lambda(k+1)$ and $\nu(k+1)$.  Then, take the preimages of the submodules that occur this composition series of $M^k_i/M^k_{i-1}$ for each $i>0$ not in $\mathcal{J}(k)$.  Using these preimages, we then refine the composition series (\ref{eq:lambdak}) and obtain a composition series 
		\begin{align}
			0=M^{k+1}_0\subsetneq M^{k+1}_1\subsetneq M^{k+1}_2 \subsetneq \ldots \subsetneq M^{k+1}_{l(k+1)-1}\subsetneq M^{k+1}_{l(k+1)}=M
		\end{align}
		with bounds $\lambda(k+1)$ and $\nu(k+1)$, along with an inclusion $\iota_k:\mathcal{J}(k)\to\mathcal{J}(k+1)$ of totally ordered sets satisfying $M^k_i/M^k_{i-1}\cong M^{k+1}_{\iota_k(i)}/M^{k+1}_{\iota_k(i)-1}$ for every $j\in\mathcal{J}(k)$.

		We now take $\mathcal{J}$ to be the direct limit $\lim_{\underset{k}{\boldsymbol{\longrightarrow}}}\,\mathcal{J}(k)$.  By construction, there is a total order $<$ on $\mathcal{J}$ extending the total orders on the sets $\mathcal{J}(k)$.  Each $j\in\mathcal{J}$ corresponds to an element in $j(k)\in \mathcal{J}(k)$ for some large enough $k$, and to a submodule $M_j:=M^k_{j(k)}$ of $M$ in the composition series from the $k$-th step.  Note that $M_{<j}=\bigcup_{j'<j}\,M_j=M^k_{j(k)-1}$, whence $M_j/M_{<j}=M^k_{j(k)}/M^k_{j(k)-1}$ is simple.   Clearly, the index set $\mathcal{J}$ and  the family of submodules $\left(M_j\right)_{j\in\mathcal{J}}$ satisfy the requirements.
		
		Finally, suppose that $M$ has a direct sum decomposition $M=\bigoplus_{t\in \pazocal{I}}\,D_t$, where each $D_t$ is indecomposable (by Theorem~\ref{thm:directsumdecomposition}).  We first equip $\pazocal{I}$ with a well order $\triangleleft$ (which exists by the Well-Ordering Principle).  Then, we create a generalized filtration series $\left\{D_t(j)\right\}_{j\in \mathcal{J}_t}$ for each $D_t$.  Write $<_t$ for the total order on $\mathcal{J}_t$.  Let $\pazocal{J}$ be the totally ordered set
		\begin{align}
			\pazocal{J}:=\big\{(t,j)\,\suchthat{t\in\pazocal{I}\text{ and }j\in \mathcal{J}_t}\big\}
		\end{align}
		with the total order $<$ defined via the lexicographic ordering as follows:
		\begin{align}
			(t,j)<(\tilde{t},\tilde{j})\text{ if and only if }t\triangleleft \tilde{t}\,,\text{ or }t=\tilde{t}\text{ and }j\,{<_t}\,\tilde{j}\,.
		\end{align}
		Then, we take $M_{(t,j)}:=D_t(j)\oplus \left(\bigoplus_{\tilde{t}\triangleleft t}\,D_{\tilde{t}}\right)$ for every $(t,j)\in\pazocal{I}$.  Then, it is obvious that $\left\{M_{(t,j)}\right\}_{(t,j)\in \pazocal{I}}$ is a generalized composition series of $M$.
	\end{pf}
	
\end{thm}

\begin{define}
	Two generalized composition series $\left(M_j\right)_{j\in\mathcal{J}}$ and $\left(M'_{j'}\right)_{j'\in\mathcal{J}'}$ are \emph{equivalent} if there exist a bijective function $f:\mathcal{J}\to\mathcal{J}'$ such that $M_j/M_{<\,j}\cong M'_{f(j)}/M'_{<'\,f(j)}$ for each $j\in\mathcal{J}$.  Here, $<$ is the total order on $\mathcal{J}$, whereas $<'$ is the total order on $\mathcal{J}'$.  In addition, $M_{< \,j}:=\bigcup_{k\,<\, j}\,M_k$ as well as $M'_{<'\,j'}:=\bigcup_{k'\,<'\,j'}\,M'_{k'}$ for all $j\in\mathcal{J}$ and $j'\in\mathcal{J}'$.
	\label{def:equivcomp}
\end{define}


\begin{thm}
	Any two generalized composition series of $M\in\bar\bggO$ are equivalent.  
	\label{thm:uniquegencomp}
	
	\begin{pf}
		We may first assume that $M$ is indecomposable.  Let $\left(M_j\right)_{j\in \mathcal{J}}$ and $\left(M'_{j'}\right)_{j'\in\mathcal{J}'}$ be two generalized composition series of an object $M\in\bbggO$. We create a decreasing sequence of weights $\big\{\lambda(k)\big\}_{k\in\mathbb{Z}_{>0}}$ and an increasing sequence of weights $\big\{\nu(k)\big\}_{k\in\mathbb{Z}_{>0}}$ such that every weight $\mu \in\text{supp}(M)$ satisfies $\lambda(k)\preceq\mu\preceq\nu(k)$ for some $k$. 
		
		For each $k$, define
		\begin{align}
			\mathcal{J}(k):= \Big\{j\in \mathcal{J}\,\suchthat{M_j/M_{<\,j}\cong \Llie(\xi)\text{ for some }\xi\text{ with }\lambda(k)\preceq\xi\preceq\nu(k)}\Big\}
		\end{align}
		and
		\begin{align}			
			\mathcal{J}'(k):= \Big\{j'\in \mathcal{J}'\,\suchthat{M'_j/M'_{<'\,j'}\cong \Llie(\xi)\text{ for some }\xi\text{ with }\lambda(k)\preceq\xi\preceq\nu(k)}\Big\}
		\end{align}
		From Theorem~\ref{thm:kacequiv}, we have a bijection $f_k:\mathcal{J}(k)\to\mathcal{J}'(k)$ with the property that $M_j/M_{<\,j}\cong M'_{f_k(j)}/M'_{<'\,f_k(j)}$ for every $j\in\mathcal{J}(k)$.  We claim that the bijections $f_k$ can be aligned so that $\left.f_{k+1}\right|_{\mathcal{J}(k)}=f_k$.   Taking the direct limit $f:=\lim_{\underset{k}{\boldsymbol{\longrightarrow}}}\,f_k$ then yields a bijection $f:\mathcal{J}\to\mathcal{J}'$ satisfying the requirement of Definition~\ref{def:equivcomp}.
		
		To prove the claim above, assume that $\left.f_{k+1}\right|_{\mathcal{J}(k)}\neq f(k)$.  Then we define the function $\tilde{f}_{k+1}:\mathcal{J}(k+1)\to\mathcal{J}'(k+1)$ as follows: 
		\begin{align}
			\tilde{f}_{k+1}(j)=\left\{\begin{array}{ll}
				f_k(j) &\text{if }j\in \mathcal{J}(k)\,,\\
				\tilde{f}_{k+1}(j)&	\text{if }j\in\mathcal{J}(k+1)\smallsetminus\mathcal{J}(k)\,.
			\end{array}\right.
		\end{align}  
		Replacing $f_{k+1}$ by $\tilde{f}_{k+1}$ and continuing this process for each positive integer $k$, we obtain a set of aligned bijections as desired.

		When $M$ is decomposable, we note that it has a unique direct sum decomposition into indecomposable direct summands (Theorem~\ref{thm:directsumdecomposition} and Corollary~\ref{cor:KSRA}).  From this, we can easily conclude that two generalized composition series of $M$ are equivalent.
	\end{pf}
\end{thm}

\begin{rmk}
	Theorem~\ref{thm:uniquegencomp} does not hold in general for a module of infinite length over an arbitrary ring $R$.  For example, $\amsbb{Z}$ as a module over itself has the following filtration
	\begin{align}
		\amsbb{Z}\supsetneq 2\,\amsbb{Z}\supsetneq 2^2\,\amsbb{Z}\supsetneq 2^3\,\amsbb{Z}\supsetneq\ldots\,,
		\label{eq:gencompZ}
	\end{align}
	which is a generalized composition series of $\amsbb{Z}$ in the same sense as in Theorem~\ref{thm:gencomp}.  However, $\amsbb{Z}$ admits a different generalized composition series inequivalent to (\ref{eq:gencompZ}):
	\begin{align}
		\amsbb{Z}\supsetneq 3\,\amsbb{Z}\supsetneq 3^2\,\amsbb{Z}\supsetneq 3^3\,\amsbb{Z}\supsetneq\ldots\,.
	\end{align}
\end{rmk}

\begin{exm}
	If a module $M\in\bbggO$ has finite length, then any generalized composition series of $M$ is a composition series.  As we shall see from the filtration (\ref{eq:vermainfinite}),  the Verma module $\Mlie\big((2,2,3,4,5,6,7,\ldots)\big)$ over $\gllie_\infty$ has a generalized composition series with an index set isomorphic to $\amsbb{Z}_{<0}$ as an ordered set.  Its dual will then have a generalized composition series with an index set isomorphic to $\amsbb{Z}_{>0}$.
	
	In Section~\ref{sec:injectiveobjects}, injective objects in $\bbggO$ and in some of its subcategories are introduced.  A generalized composition series of such an injective object is given by an index set isomorphic to a totally ordered set isomorphic to a disjoint union of finitely many copies of $\amsbb{Z}_{>0}$ on top of one another.
\end{exm}

\begin{openq}
	What are possible index sets (up to isomorphism of ordered sets) of generalized composition series of indecomposable objects in $\bbggO$?  In particular, is it true that every indecomposable object $M$ of $\bbggO$ admits a generalized composition series ordered by a subset of $\amsbb{Z}$ with the usual order?
\end{openq}

\begin{cor}
	Let $M\in\bar\bggO$ and $\lambda\in\hfrak^*$.  If $M=\bigoplus_{j\in J}\,M_j$ for objects $M_j\in\bar\bggO$ and some index set $J$, then $\big[M:\Llie(\lambda)\big]=\sum_{j\in J}\,\big[M_j:\Llie(\lambda)\big]$.  Note that this sum is always finite.
\end{cor}

\begin{rmk}
	We can also define a generalized composition series of a module $M\in\bbggO$ when the splitting Borel subalgebra is non-Dynkin.  The idea is to also first deal with the case where $M$ is indecomposable, so there exist finite subsets $S_1,S_2,S_3,\ldots$ of $\text{supp}(M)$ such that $S_1\subseteq S_2\subseteq S_3\subseteq \ldots$ and that $\bigcup_{r=1}^\infty\,S_r=\text{supp}(M)$.  Then, we first construct a composition series of $M$ with weight reference set $S_1$.  For each $r>1$, we refine the previously obtained composition series of $M$ with  weight reference set $S_{r-1}$ to get a composition series of $M$ with  weight reference set $S_{r}$.  By taking the direct limit as in the proof of Theorem~\ref{thm:gencomp}, we obtain a \emph{generalized composition series} of $M$.  Two generalized composition series of $M$ are also equivalent, and the corollary above holds as well. The details are omitted here as our focus is on Dynkin Borel subalgebras. 
\end{rmk}

We now assume again that $\bfrak$ is a Dynkin Borel subalgebra.

\begin{define}
Let $\lambda\in\hfrak^*$.  Let $W[\lambda]$ be the subgroup of $W$ containing all $w\in W$ such that $w\cdot\lambda-\lambda \in \Lambda$.  Write $W_n$ for the Weyl group of $\gfrak_n$.    We similarly define $W_n\left[\lambda_n\right]$ for each $n\in\amsbb{Z}_{>0}$ and $\lambda_n\in\hfrak_n^*$.  These subgroups are known as the \emph{integral Weyl groups}.
\label{def:Wlambda}
\end{define}

\begin{thm}
	Let $\lambda,\mu\in\hfrak^*$.  If the simple module $\Llie(\mu)$ is a composition factor of the Verma module $\Mlie(\lambda)$, then $\mu\preceq \lambda$ and $\mu \in W[\lambda]\cdot\lambda$.
	\label{thm:compositionfactors}
	
	\begin{pf}
		Suppose that $\Mlie(\lambda)$ has $\Llie(\mu)$ as a composition factor.   Let $v$ be a highest-weight vector of $\Mlie(\lambda)$.  Then, for all sufficiently large $n\in\mathbb{N}$, $M_n:=\Ulie\left(\gfrak_n\right)\cdot v \in\bggO^{\gfrak_n}_{\hfrak_n}$ is a Verma module and must have $\gfrak_n$-submodules $N_n$ and $N'_n$ with $N_n\subseteq N'_n$ and $N'_n/N_n$ has $\mu_n:=\mu|_{\hfrak_n}$ as a highest weight.  Hence, $\Llie\left(\mu_n\right)$ is a composition factor of $M_n\cong\Mlie\left(\lambda_n\right)$, where $\lambda_n:=\lambda|_{\hfrak_n}$.  Due to the finite-dimensional theory, $\mu_n\preceq\lambda_n$ and $\mu_n\in W_n\left[\lambda_n\right]\cdot\lambda_n$.  The result follows immediately.
	\end{pf}
\end{thm}

\pagebreak

\section{Verma Modules}
\label{ch:verma}

In this section, $\gfrak$ is an arbitrary root-reductive Lie algebra.

\subsection{Fundamental Properties of Verma Modules}
	\label{ch:fundamentals}
	
	In this subsection, $\bfrak$ is an arbitrary splitting Borel subalgebra of $\gfrak$ containing a splitting maximal toral subalgebra $\hfrak$.  
	
	\begin{thm} 
		A Verma module has at most one simple Verma submodule.  If $\mathfrak{b}$ is Dynkin, then a Verma module has at most one simple submodule.
		\label{thm:simplesubmodules}
	
	\begin{pf}
	Let $\lambda \in \mathfrak{h}^*$ and $M:=\Mlie(\lambda)$.  Suppose that $N_1$ and $N_2$ are highest-weight $\Ulie(\gfrak)$-submodules of $M$ with highest-weight vectors $v_1 \neq 0$ and $v_2 \neq 0$, respectively. 
	Let $R$ denote the universal enveloping algebra $\Ulie\left(\mathfrak{n}^-\right)$.  As we have seen before (in Section~\ref{ch:wtmodules}), $M$ is isomorphic to $R$ as an $R$-module.  We can identify $M$ with $R$, making $v_1$ and $v_2$ elements of $R$.  Ergo, $N_1$ and $N_2$ are left $S$-ideals $R\cdot v_1$ and $R\cdot v_2$, respectively.

	Because $\mathfrak{g}$ is locally finite, the Lie subalgebra $\mathfrak{n}^-$ is also locally finite.  Thus, there exists a finite-dimensional Lie subalgebra $\underline{\mathfrak{n}}^-$ of $\mathfrak{n}^-$ that contains the elements of $\mathfrak{n}^-$ involved in the PBW polynomial expressions for $v_1$ and $v_2$.  We take $\underline{R}$ to be the universal enveloping algebra of $\underline{\nfrak}^-$.  Consequently, $\underline{R}$ is a noetherian ring.\index{noetherian ring}  We want to show that the left $\underline{R}$-ideals $\underline{R}\cdot v_1$ and $\underline{R}\cdot v_2$ intersect nontrivially.
	
	A more general statement is true.  Let $\mathscr{R}$ be a left noetherian ring.  If $x \in \mathscr{R}$  is not a right zero-divisor, then the left ideal $\mathscr{R}\boldsymbol{\cdot} x$ intersects every nonzero left ideal of $\mathscr{R}$ nontrivially.  In particular, if $\mathscr{R}$ has no right zero-divisors, then any two nonzero left ideals of $\mathscr{R}$ intersect nontrivially.  (For a proof, see Lemma 4.1 in \cite{bggo}.)
	
	From the paragraph above, we conclude that $\underline{R}\cdot v_1$ must intersect $\underline{R}\cdot v_2$ nontrivially.  Thence, $N_1$ and $N_2$  intersect nontrivially as well.  As a result, if $N_1$ and $N_2$ are both simple highest-weight $\Ulie(\gfrak)$-submodules of $M$, then $N_1=N_2$.  In other words, every Verma module over $\mathfrak{g}$ has at most one simple Verma submodule. 
	\end{pf}
	\end{thm}
	
	\begin{cor}
		If a Verma module $M$ has a simple Verma submodule $L$, then every Verma submodule of $M$ contains $L$.  In particular, if $\bfrak$ is Dynkin and $M$ has a simple submodule $L$, then all nontrivial submodules of $M$ must contain $L$, i.e., $L$ equals the socle of $M$.
		\label{cor:simplevermasubmodule}
	\end{cor}
	
	In parallel to the finite-dimensional case, the theorem below gives a generalized version to items (b) and (c) of Theorem~\ref{thm:findim}.  The next subsection partially offers a condition under which there exists an embedding of a Verma module into another Verma module.
	
	\begin{thm}
		Let $\lambda,\mu\in\hfrak^*$.  Then,
		\begin{align}
		\dim_\amsbb{K}\Big(\Hom_{\Ulie(\gfrak)}\big(\Mlie(\lambda),\Mlie(\mu)\big)\Big)\in\{0,1\}\,.
		\end{align}  Furthermore, all nonzero elements of  $\Hom_{\Ulie(\gfrak)}\big(\Mlie(\lambda),\Mlie(\mu)\big)$ are embeddings.  If a nonzero homomorphism exists, then $\lambda \preceq \mu$.
	\label{thm:embeddings}
		
	\begin{pf}
	Let us suppose that $\lambda,\mu \in \mathfrak{h}^*$ are such that there exist two nonzero homomorphisms $\phi_1,\phi_2:\Mlie(\lambda)\to\Mlie(\mu)$ of $\Ulie(\gfrak)$-modules.  We shall prove that $\phi_2=\kappa \, \phi_1$ for some $\kappa \in \amsbb{K}$.  Let $v_\lambda \neq 0$ and $v_\mu \neq 0$ be highest-weight vectors\index{vector!highest-weight vector} of $\Mlie(\lambda)$ and $\Mlie(\mu)$, respectively.  Write $w_i:=\phi_i\left(v_\lambda\right)$ for $i \in \{1,2\}$.   We identify $\Mlie(\mu)$ as a $\Ulie\left(\mathfrak{n}^-\right)$-module which is isomorphic to $\Ulie\left(\mathfrak{n}^-\right)$ itself.  Ergo, $w_1$ and $w_2$ are now elements of $\Ulie\left(\mathfrak{n}^-\right)$.   By the local finiteness of $\mathfrak{g}$, there exists a finite-dimensional  Lie subalgebra $\underline{\mathfrak{g}}$ with a Borel subalgebra $\underline{\mathfrak{b}}:=\mathfrak{b}\cap\underline{\mathfrak{g}}$ that contains a maximal toral subalgebra $\underline{\mathfrak{h}}:=\mathfrak{h}\cap\underline{\mathfrak{g}}$.  Then, $\underline{M}:=\Ulie\left(\underline{\mathfrak{g}}\right)\cdot v_\mu$ is a Verma module over $\underline{\mathfrak{g}}$.  Now, $\Ulie\left(\underline{\mathfrak{g}}\right)\cdot w_1$ and $\Ulie\left(\underline{\mathfrak{g}}\right)\cdot w_2$ are isomorphic Verma modules over $\underline{\mathfrak{g}}$, both of which are embedded into $\underline{M}$.  Since, in the finite-dimensional case, the homomorphism space between two Verma modules is either trivial or one-dimensional.  Therefore, we must have $\Ulie\left(\underline{\mathfrak{g}}\right)\cdot w_1=\Ulie\left(\underline{\mathfrak{g}}\right)\cdot w_2$.  Consequently, $w_2=\kappa\, w_1$ for some nonzero $\kappa \in \amsbb{K}$.  This means $\phi_2=\kappa\, \phi_1$.  Hence, $\text{Hom}_{\Ulie(\gfrak)}\big(\Mlie(\lambda),\Mlie(\mu)\big)$ is of dimension $0$ or $1$ over $\amsbb{K}$. 
	
	To show that any nonzero homomorphism in  $\text{Hom}_{\Ulie(\gfrak)}\big(\Mlie(\lambda),\Mlie(\mu)\big)$ must be an embedding, let $\phi$ be such a map.  Via the identification of $\Mlie(\lambda)$ and $\Mlie(\mu)$ with $\Ulie\left(\nfrak^-\right)$ as left $\Ulie\left(\nfrak^-\right)$ modules, we can easily see that $\phi$ is the multiplication map $x\mapsto u\boldsymbol{\cdot}x$ for some $ u \in \Ulie\left(\nfrak^-\right)$ and for all $x \in \Ulie\left(\nfrak^-\right)$.  Because $\Ulie\left(\nfrak^-\right)$ lacks zero divisors, the map $\phi$ must be injective.
	\end{pf}
	\end{thm}
	
	Unlike in the finite-dimensional case, we do not know whether there exists a Verma module over $\gfrak$ which has more than one simple submodules.  However, the uniqueness of simple submodules is guaranteed if $\bfrak$ is a Dynkin Borel subalgebra.  We shall see later that, in contrast with the finite-dimensional theory, there are Verma modules with no simple submodules (see Section~\ref{ch:stdborel}).\index{simple module|see {Lie module}}\index{Lie module!simple module}   
			
\subsection{Verma Modules for Dynkin Borel Subalgebras}
	\label{ch:simpleroots}
	\index{sufficiency of simple roots}
	
	It can be easily seen that all Verma modules are objects of the category $\bar\bggO^\gfrak_\bfrak$ if and only if $\bfrak$ is a Dynkin Borel subalgebra of $\gfrak$.  This is one of the reasons why this dissertation focuses on the category $\bbggO^\gfrak_\bfrak$ for a Dynkin Borel subalgebra $\bfrak$ of a root-reductive Lie algebra $\gfrak$.  
	
	Suppose now that $\bfrak$ is Dynkin.  We have the following trivial proposition (which has been used in Theorem~\ref{thm:simplesubmodules}, Corollary~\ref{cor:simplevermasubmodule}, and Theorem~\ref{thm:embeddings}).
	
	\begin{prop}
		Let $\lambda\in\hfrak^*$.  
		\begin{itemize}
			\item[(a)] Every weight space of $\Mlie(\lambda)$ is finite dimensional.  
			\item[(b)] Every submodule of $\Mlie(\lambda)$ has a singular vector. 
			\item[(c)]  Any simple submodule of $\Mlie(\lambda)$ is also a Verma module.  
			\item[(d)] For $\mu\in\hfrak^*$ such that $\lambda \preceq \mu$, there exists only finitely many $\xi\in\hfrak^*$ such that $\lambda \preceq \xi \preceq \mu$.
		\end{itemize}
	\end{prop}
	
	
	
	
	\begin{define}
		For a root $\alpha \in \Delta$, let $h_\alpha$ be the unique element of $\left[\gfrak^{+\alpha},\gfrak^{-\alpha}\right]$ such that $\alpha\left(h_\alpha\right)=2$.\nomenclature[B]{$h_\alpha$}{the unique element of $\left[\gfrak^{+\alpha},\gfrak^{-\alpha}\right]$ such that $\alpha\left(h_\alpha\right)=2$, where $\alpha$ is a root}
			\begin{enumerate}
				\item[(a)] We say that $\lambda \in \hfrak^*$ is \emph{integral}\index{weight!integral weight} if $\lambda\left(h_\alpha\right) \in \amsbb{Z}$ for every $\alpha \in \Delta$.
				\item[(b)] We say that $\lambda \in \hfrak^*$ is \emph{antidominant}\index{weight!antidominant weight} if $(\lambda+\rho)\left(h_\alpha\right) \notin \amsbb{Z}_{> 0}$ for any $\alpha \in \Delta^+$.\nomenclature[B]{$\amsbb{Z}_{>0}$}{the set of positive integers}
				\item[(c)] We say that $\lambda \in \hfrak^*$ is \emph{almost antidominant}\index{weight!almost antidominant weight} if $(\lambda+\rho)\left(h_\alpha\right) \in \amsbb{Z}_{> 0}$ for only finitely many $\alpha \in \Delta^+$.
			\end{enumerate}
			\label{def:wttypes}
		\end{define}

	It turns out that Verma's Theorem (Theorem \ref{thm:verma}) and  the BGG Theorem (Theorem \ref{thm:stronglinkage}) hold also for root-reductive Lie algebras.   The generalizations below shall be called Verma's Theorem and BGG Theorem as well.\index{Verma's Theorem}\index{BGG Theorem}
	
	\begin{thm}[\bf Verma's Theorem]
			For $\lambda \in \hfrak^*$ and for a given a positive root $\alpha$ such that $s_\alpha\cdot \lambda \preceq \lambda$, there exists an embedding $\Mlie\left(s_\alpha\cdot\lambda\right)\overset{\subseteq}{\longrightarrow}\Mlie(\lambda)$.
			\label{thm:vermainf}
			\begin{pf}
				For $n \in \amsbb{Z}_{> 0}$, write $\bfrak_n$ and $\hfrak_n$ for $\bfrak\cap\gfrak_n$ and $\hfrak\cap\gfrak_n$, respectively.   Let $\lambda_n$ be the restriction of $\lambda$ onto $\hfrak_n$.  Denote by $M$ the Verma module $\Mlie(\lambda;\gfrak,\bfrak,\hfrak)$, while $M_n$ is the Verma module $\Mlie\left(\lambda_n;\gfrak_n,\bfrak_n,\hfrak_n\right)$.  If $u$ is a highest-weight vector of $M$\index{vector!highest-weight vector}, then by identifying a highest-weight vector of $M_n$ with $u$, we have $M_1 \subseteq M_2 \subseteq M_3 \subseteq \ldots$.  In other words, $M$ is the direct limit of $\left(M_n\right)_{n\in\amsbb{Z}_{> 0}}$ under inclusion maps.
				
				 Define
				 \begin{align}
				 	\Psi:=\left\{\gamma\in\Delta^+\suchthat{} \gamma \preceq \lambda -s_\alpha\cdot\lambda\right\}\,.
				\end{align}  The set $\Psi$ is clearly finite.  Therefore, for sufficiently large values of $n$, say $n \geq m$ for some $m \in \amsbb{Z}_{> 0}$, we have $\gfrak^{-\beta}\subseteq \gfrak_n$ for all $\beta \in \Psi$, which further implies that $\alpha|_{\hfrak_n}$ is a positive root of $\gfrak_n$ and that $s_{\alpha} \in W_n$.  Thus, for such $n\in\amsbb{Z}_{> 0}$, the Verma module $\Mlie\left(s_\alpha\cdot\lambda_n;\gfrak_n,\bfrak_n,\hfrak_n\right)$ is isomorphic to a unique $\gfrak_n$-submodule $N_n$ of $M_n$, where we have applied Theorem~\ref{thm:findim} and Verma's Theorem (Theorem \ref{thm:verma}).  

				Now, observe that, for $n \geq m$, $s_\alpha \cdot \lambda_n \in \hfrak_n^*$ is identical to the restriction of $s_\alpha \cdot\lambda_{n+1} \in \hfrak_{n+1}^*$ onto $\hfrak_n$.  Furthermore, the weight space associated to the weight $s_\alpha\cdot\lambda_n$ of $M_n$ (where the dot action is done in $\hfrak_n^*$) is precisely $M^{s_\alpha\cdot\lambda}$.    This means that the highest-weight spaces of $N_n$ and of $N_{n+1}$, which correspond to the weights $s_\alpha\cdot\lambda_n$ and $s_\alpha\cdot\lambda_{n+1}$, respectively, are identical for $n \geq m$.  That is, $N_n \subseteq N_{n+1}$ for every integer $n \geq m$.  The direct limit $N:=\lim_{\underset{n}{\boldsymbol{\longrightarrow}}}\, N_n$ of $\left(N_n\right)_{n  \geq m}$ under inclusion maps is thus a $\gfrak$-submodule of $M$ isomorphic to $\Mlie(s_\alpha\cdot\lambda;\gfrak,\bfrak,\hfrak)$.
			\end{pf}
			
	\end{thm}
	
	
	\begin{thm}[\bf BGG Theorem]
		For $\lambda,\mu\in\hfrak^*$, there exists a nontrivial $\gfrak$-module homomorphism from $\Mlie(\lambda)$ to $\Mlie(\mu)$ if and only if $\lambda$ is \emph{strongly linked}\index{strong linkage} to $\mu$, namely, there exist positive roots $\alpha_1,\alpha_2,\ldots,\alpha_l$ such that
		\begin{align}
			\lambda = \left(s_{\alpha_l}\cdots s_{\alpha_2}s_{\alpha_1}\right)\cdot\mu \preceq \left(s_{\alpha_{l-1}}\cdots s_{\alpha_2}s_{\alpha_1}\right)\cdot\mu \preceq \ldots \preceq s_{\alpha_1}\cdot\mu \preceq \mu\,.
			\label{eq:stronglinkage}
		\end{align}
		That is, for $\mu \in \hfrak^*$, all Verma submodules of $\Mlie(\mu)$ is of the form $\Mlie(w\cdot \mu)$, where $w$ is an element of the Weyl group.
		\label{thm:stronglinkagegen}
		
		\begin{pf}
			The converse (i.e., that a strong linkage implies the existence of an embedding) is clear, so we prove the other direction (i.e., that an embedding implies the existence of a strong linkage).  Suppose that $\Mlie(\lambda)$ is a submodule $N$ of $M:=\Mlie(\mu)$.  Let $u$ and $v$ be highest-weight vectors of $M$ and $N$, respectively.  For each $n \in \amsbb{Z}_{> 0}$, write $M_n$ for $\Ulie\left(\gfrak_n\right)\cdot u$.  As $M=\lim_{\underset{n}{\boldsymbol{\longrightarrow}}}\,M_n$, there exists $m \in \amsbb{Z}_{> 0}$ such that $n \geq m$ implies $v \in M_n$.  
			
			For $n \geq m$, write $N_n$ for $\Ulie\left(\gfrak_n\right)\cdot v$.  Then, $N_n$ is a Verma submodule of $M_n$ (over $\gfrak_n$).  The finite-dimensional BGG Theorem (Theorem \ref{thm:stronglinkage}) guarantees that $\lambda_n:=\lambda|_{\hfrak_n}$ is strongly linked to $\mu_n:=\mu|_{\hfrak_n}$.  The positive roots $\alpha_n^j$, $j=1,2,\ldots,l_n$, involved in the $n$-th linkage
			\begin{align}
				\lambda_n  = \left(s_{\alpha_n^{l_n}}\cdots s_{\alpha_n^2}s_{\alpha_n^1}\right)\cdot \mu_n \preceq \left(s_{\alpha_n^{l_n-1}}\cdots s_{\alpha_n^2}s_{\alpha_n^1}\right)\cdot \mu_n \preceq \ldots \preceq s_{\alpha_n^1}\cdot \mu_n \preceq \mu_n
			\end{align} 
			between $\lambda_n$ and $\mu_n$ must belong to the set $\left\{\alpha \in \Delta^+\suchthat{} \alpha \preceq \mu-\lambda\right\}$, which is a finite set.  
			
			If $\mu-\lambda =\sum_{\alpha \in \Sigma}t_\alpha \alpha$, where $t_\alpha \in \amsbb{Z}_{\geq 0}$ for each $\alpha \in \Sigma$, then the lenth $l_n$ of the $n$-th linkage is at most $\sum_{\alpha \in \Sigma}t_\alpha <\infty$.  Using the Pigeonhole Principle\index{Pigeonhole Principle}, it follows that there are infinitely many $n \geq m$ with the same linkage pattern, say
			\begin{align}
				\lambda_n  = \left(s_{\alpha_{l}}\cdots s_{\alpha_2}s_{\alpha_1}\right)\cdot \mu_n \preceq \left(s_{\alpha_{l-1}}\cdots s_{\alpha_2}s_{\alpha_1}\right)\cdot \mu_n \preceq \ldots \preceq s_{\alpha_1}\cdot \mu_n \preceq \mu_n\,,
			\end{align} 
			where $\alpha_1,\alpha_2,\ldots,\alpha_l$ are positive roots.  Hence, (\ref{eq:stronglinkage}) holds.
		\end{pf}
	\end{thm}

	\begin{thm}
			For $\lambda \in\hfrak^*$, $\Mlie(\lambda)$ is simple if and only if $\lambda$ is antidominant.\index{Lie module!simple module}
			\label{thm:simpleverma}
			
			\begin{pf}
				For each root $\alpha$, $s_\alpha$ is the reflection with respect to $\alpha$ and $h_{\alpha}$ is as defined in Definition~\ref{def:wttypes}.
				
				\begin{itemize}
					\item[($\Rightarrow$)]  Suppose that $\lambda$ is not antidominant. Then, there exists a positive root $\alpha$ such that $(\lambda+\rho)\left(h_\alpha\right) \in \amsbb{Z}_{>0}$.  This means $s_\alpha \cdot \lambda \precneq \lambda$ and
					$
						0 \subsetneq \Mlie\left(s_\alpha\cdot\lambda\right)\subsetneq \Mlie(\lambda)
					$, where we have applied Verma's Theorem (Theorem \ref{thm:vermainf}); as a result, $\Mlie(\lambda)$ is not simple.
					\item[($\Leftarrow$)]  
					Suppose that $\Mlie(\lambda)$ is not simple.  Then, it has a proper nonzero submodule, which must have a highest-weight vector whose weight is $\mu \in \hfrak^*$.\index{vector!highest-weight vector}  Then, $\Mlie(\mu)$ is a proper Verma submodule of $\Mlie(\lambda)$, so $\mu \precneq \lambda$.\nomenclature[B]{$\precneq$}{the strict version of $\preceq$}  Using the BGG Theorem (Theorem \ref{thm:stronglinkagegen}), there are positive roots $\alpha_1$, $\alpha_2$, $\ldots$, $\alpha_l$ with $l \in \amsbb{Z}_{> 0}$ such that 
					\begin{align}
						\mu= \left(s_{\alpha_l}\cdots s_{\alpha_2}s_{\alpha_1}\right)\cdot\lambda \preceq \left(s_{\alpha_{l-1}}\cdots s_{\alpha_2}s_{\alpha_1}\right)\cdot\lambda \preceq \ldots \preceq s_{\alpha_1}\cdot\lambda \preceq \lambda\,.
					\end{align}
					Because $s_{\alpha_1}\cdot\lambda \preceq \lambda$, we have $(\lambda+\rho)\left(h_{\alpha_1}\right) \in \amsbb{Z}_{>0}$.  Thence, $\lambda$ is not antidominant.
				\end{itemize}
			\end{pf}
		\end{thm}
		
		\begin{thm}
			Let $\lambda  \in \hfrak^*$.  The following conditions are equivalent:
			\begin{itemize}
				\item[(a)]  The module $\Mlie(\lambda)$ is of finite length.
				\item[(b)] The module $\Mlie(\lambda)$ has a simple submodule.
				\item[(c)] There exists an antidominant weight $\xi$ such that $\xi$ is strongly linked to $\lambda$.
				\item[(d)] The weight $\lambda$ is almost antidominant.
			\end{itemize}
			\label{thm:finitelength}
			
			\begin{pf}
				For simplicity, we shall denote $M$ for $\Mlie(\lambda)$.
				
				\begin{itemize}
				\itemindent1cm
					\item[$\!\!\!\!\!\big($(a)$\Leftrightarrow$(b)$\big)$]  For the direct implication, let 
						$
							0=M_0\subsetneq M_1 \subsetneq M_2 \subsetneq \ldots \subsetneq M_{l-1}\subsetneq M_l=M
						$ be a composition series of $M$ for some $l \in \amsbb{Z}_{>0}$.  Then, $M_1$ must be simple.  
						
						Conversely, let $L$ be a simple submodule of $M$.  Then, $L$ is a Verma module of highest weight $\mu \in \hfrak^*$ with $\mu \preceq \lambda$.  Note that every nonzero submodule of $M$ must include $L$.  Any singular weight $\xi$ of $M$ must then satisfy $\mu\preceq\xi\preceq\lambda$.  There are only finitely many weights $\xi$ for which $\mu \preceq \xi \preceq \lambda$, and  the weight space with weight $\xi$ is finite dimensional for each $\xi \in \hfrak^*$.  If $m$ is the sum of the dimensions of all the weight spaces with weight $\xi \in \hfrak^*$ such that $\mu \preceq \xi \preceq \lambda$, then we can immediately see that $M$ has a composition series of length at most $m$.  Hence, $M$ is of finite length.
						
					\item[$\!\!\!\!\!\big($(b)$\Leftrightarrow$(c)$\big)$] We can easily apply Verma's Theorem (Theorem \ref{thm:vermainf}), the BGG Theorem (Theorem \ref{thm:stronglinkagegen}), and Theorem~\ref{thm:simpleverma} to verify that (b) and (c) are equivalent. 
					
						\item[$\!\!\!\!\!\big($(b)$\Rightarrow$(d)$\big)$]
						Suppose $L$ is a simple submodule of $M$.  Then, $L$ is a Verma module with the highest weight $\mu \preceq \lambda$, for some $\mu \in \hfrak^*$.  By the BGG Theorem (Theorem \ref{thm:stronglinkagegen}), $w\cdot\lambda = \mu$ for some $w \in W$.  Since $\mu$ is antidominant (by Theorem~\ref{thm:simpleverma}), this means $\lambda$ is almost antidominant.
						
						\item[$\!\!\!\!\!\big($(d)$\Rightarrow$(c)$\big)$]	Suppose that $\lambda$ is almost antidominant.  For each $\mu\in\hfrak^*$, let $\Xi(\mu)$ denotes the set of positive roots $\alpha$ such that $h_\alpha(\lambda+\rho)$ is a positive integer.  We say that $\alpha \in \Xi(\mu)$ is \emph{minimal} if $\alpha$ cannot be written as a sum of at least two elements of $\Xi(\mu)$.  Let $\Upsilon(\mu)$ be the cardinality of $\Xi(\mu)$.  
		 	
		 	By the assumption, $\Upsilon(\lambda)<\infty$.  We shall prove by induction on $\Upsilon(\lambda)$.  Pick a minimal $\alpha\in \Xi(\lambda)$.  Then, we have $s_\alpha\cdot \lambda\preceq \lambda$ and $\Upsilon(s_\alpha\cdot \lambda)<\Upsilon(\lambda)$.  By the induction hypothesis, there exists an antidominant weight $\xi$ such that $\xi$ is strongly linked to $s_\alpha\cdot\lambda$.  That is, there are positive roots $\alpha_2,\alpha_3,\ldots,\alpha_l$ such that 
				\begin{align}
					\xi = \left(s_{\alpha_l}\cdots s_{\alpha_2}\right)\cdot\left(s_{\alpha_1}\cdot \lambda\right)& \preceq \left(s_{\alpha_{l-1}}\cdots s_{\alpha_2}\right)\cdot\left(s_{\alpha_1}\cdot \lambda\right) \nonumber\\&\preceq \ldots \preceq s_{\alpha_2}\cdot\left(s_{\alpha_1}\cdot \lambda\right) \preceq s_{\alpha_1}\cdot\lambda\,.
				\end{align}
				It follows immediately that				
				\begin{align}
			\xi = \left(s_{\alpha_l}\cdots s_{\alpha_2}s_{\alpha_1}\right)\cdot\lambda \preceq \left(s_{\alpha_{l-1}}\cdots s_{\alpha_2}s_{\alpha_1}\right)\cdot\lambda \preceq \ldots \preceq s_{\alpha_1}\cdot\lambda \preceq \lambda\,,
		\end{align}
			and our proof is now complete.
				\end{itemize}
			\end{pf}
		\end{thm}

\subsection{Examples for $\gllie_\infty$}
	\label{ch:stdborel}

	In this subsection, $\gfrak:=\gllie_\infty$.  We take $\hfrak$ to be the subalgebra $\hfrak_\text{diag}$ and $\bfrak$ to be the one-sided Dynkin Borel subalgebra $\bfrak_{\text{1st}}$ of $\gfrak$. For $n\in\amsbb{Z}_{> 0}$, $\gfrak_n$ is the subalgebra of $\gfrak$ spanned over $\amsbb{K}$ by $\textbf{E}_{i,j}$ with $i,j\in\{1,2,\ldots,n\}$.  
	
	Recall that $\bfrak$ is a Dynkin Borel subalgebra of $\gfrak$.  Its simple roots are $\epsilon_i-\epsilon_{i+1}$ for all $i=1,2,3,\ldots$, where $\epsilon_i\in\hfrak^*$ is the map sending $\textbf{E}_{i,i}$ to $1$ and $\textbf{E}_{j,j}$ to $0$ for every $j \in \amsbb{Z}_{> 0}\smallsetminus\{i\}$.   The global half sum of positive roots with respect to $\bfrak$ is given by the same formula as the half sum of positive root for $\sllie_\infty$ with the Dynkin Borel subalgebra $\bfrak_{\text{1A}}$ shown in (\ref{eq:rho}). 
	

		The Weyl group  of $W$ is the \emph{infinite symmetric group} $\mathfrak{S}_\infty$ (i.e., the group of permutations on $\amsbb{Z}_{>0}$ which fix all but finitely many numbers).  Now, if $\lambda=\left(\lambda^k\right)_{k\in\amsbb{Z}_{> 0}}\in\hfrak^*$, then we can see that the actions of $w \in W$ are given by
			\begin{align}
				w(\lambda) = \left(\lambda^{w(k)}\right)_{k\in\amsbb{Z}_{> 0}}\text{ and }w\cdot\lambda = \left(\lambda^{w(k)}-w(k)+k\right)_{k\in\amsbb{Z}_{> 0}}\,.
			\end{align}
		
		\begin{rmk}
			Let $\lambda=\left(\lambda^k\right)_{k\in\amsbb{Z}_{> 0}} \in \hfrak^*$.  In light of Definition~\ref{def:wttypes}, we observe the followings:
			\begin{enumerate}
				\item[(a)] The weight $\lambda$ is integral if and only if $\lambda^i-\lambda^j\in\amsbb{Z}$ for all $i,j\in\amsbb{Z}_{> 0}$, or equivalently, there exists $c \in \amsbb{K}$ such that $\lambda^k-c$ is an integer for every $k \in \amsbb{Z}_{> 0}$.
				\item[(b)] If $\lambda$ is integral, then it is antidominant if and only if the sequence $\left(\lambda^k\right)_{k\in\amsbb{Z}_{> 0}}$ is strictly increasing.
				\item[(c)] If $\lambda$ is integral, then it almost antidominant if and only if there exists $N \in \amsbb{Z}_{> 0}$ such that $\lambda^k-k \leq \lambda^N-N$ for every $k=1,2,\ldots,N-1$ and that $\left(\lambda^k\right)_{k \geq N}$ is strictly increasing.
			\end{enumerate}
		\end{rmk}

	Since the zero weight $\mathrel{\ooalign{\hss{\color{Gray}{$\cdot$}}\hss\cr0}} \,=(0,0,0,\ldots)$ is not antidominant (in fact, it is \emph{integral dominant}), $\Mlie(\mathrel{\ooalign{\hss{\color{Gray}{$\cdot$}}\hss\cr0}})$ is of infinite length and has no simple submodules.  On the other hand, the weight $\lambda:=(1,2,3,\ldots)$ is antidominant, so $\Mlie(\lambda)$ is simple (whence $\Mlie(\lambda)=\Llie(\lambda)$).  A small modification of $\lambda$, say $\mu:=(1,0,-1,4,5,6,\ldots)$, is almost antidominant.  The antidominant weight $\xi$ which is strongly linked to $\mu$ is $\xi:=(-3,0,3,4,5,6,\ldots)$.
			
			However, there are some weights $\nu \in \hfrak^*$ whose sequential forms $\nu=\left(\nu^1,\nu^2,\nu^3,\ldots\right)$ are eventually strictly increasing and whose corresponding Verma modules $\Mlie(\nu)$ are of infinite length.  An example is the weight $\nu:=(2,2,3,4,5,6,\ldots)$.  We can see that, with $\nu_0:=\nu$ and $\nu_k:=(k\;\;k+1)\cdot\nu_{k-1}$ for each $k=1,2,\ldots$ (where $(k\;\;k+1)\in W$ is the transposition that swaps the $k$-th coordinate and the $(k+1)$-st coordinate), we have the infinite filtration $\Mlie(\nu)=\Mlie\left(\nu_0\right) \supsetneq \Mlie\left(\nu_1\right) \supsetneq \Mlie\left(\nu_2 \right) \supsetneq \ldots $, or 
			\begin{align}
			 	\Mlie\big((2,2,3,4,5,6,7,\ldots)\big) &\supsetneq \Mlie\big((1,3,3,4,5,6,7,\ldots) \big)
			 	\nonumber
			 	\\
			 	&\supsetneq \Mlie\big((1,2,4,4,5,6,7,\ldots)\big) \supsetneq \ldots\,.
			 	\label{eq:vermainfinite}
			\end{align}
	This filtration is a generalized composition series of $\Mlie(\nu)$, with composition factors $\Llie\left(\nu_0\right)$, $\Llie\left(\nu_1\right)$, $\Llie\left(\nu_2\right)$, $\ldots$, each with multiplicity $1$.  Additionally, observe that the Verma module $\Mlie(\nu)$ has no socle.
	
\subsection{Verma Modules for General Splitting Borel Subalgebras}
 	
 	In this subsection, $\bfrak$ is an arbitrary splitting Borel subalgebra of a root-reductive Lie algebra $\gfrak$.
 	
 	\begin{define}
 		The \emph{$\bfrak$-foundational subalgebra} of $\gfrak$ is the subalgebra $\tilde{\gfrak}$ of $\gfrak$ generated by $\hfrak$ and all simple root spaces (both positive and negative).
 	\end{define}
	
	\begin{prop}
		The $\bfrak$-foundational subalgebra $\tilde{\gfrak}$ of $\gfrak$ is given by
		\begin{align}
			\tilde{\gfrak}= \hfrak\oplus\bigoplus_{\alpha \in \Delta_{\text{fl}}}\gfrak^\alpha\,,
		\end{align}
		where $\Delta_{\text{fl}}$ is the subset of $\Delta$ consisting of all roots of finite length.  The set of positive roots of finite length and the set of negative roots of finite length are denoted by $\Delta_\text{fl}^+$ and $\Delta_\text{fl}^-$, respectively.
	\end{prop}
	
	Observe that the $\bfrak$-foundational subalgebra $\tilde{\gfrak}$ of a root-reductive Lie algebra $\gfrak$ is a root-reductive Lie algebra with a splitting maximal toral subalgebra $\hfrak$.  Let $\tilde{\bfrak}$ denote $\bfrak\cap\tilde{\gfrak}$.  Then, $\tilde{\bfrak}$ is a Dynkin Borel subalgebra of $\tilde{\gfrak}$.  Consequently, $\tilde{\gfrak}$ possesses a half sum of positive roots $\tilde{\rho}$ with respect to $\tilde{\bfrak}$.

	\begin{exm}
		Let $\gfrak$ be the Lie algebra $\sllie_\infty$.  If $\bfrak$ is given by the ordering 
		\begin{align}
			1\prec 3 \prec 5 \prec \ldots \prec 6 \prec 4 \prec 2\,,
		\end{align}
		then the $\bfrak$-foundational subalgebra $\tilde{\gfrak}$ of $\gfrak$ is the direct sum $\gfrak[1]\oplus\gfrak[2]$, where $\gfrak[1]$ is given by odd indices and $\gfrak[2]$ is given by even indices.  This Borel subalgebra $\bfrak$ shall be denoted by $\bfrak_\text{si}$ (and it is called the \emph{simplest ideal Borel subalgebra}).  Note that $\gfrak[1]$ is isomorphic to $\sllie_\infty$ with the one-sided Dynkin Borel subalgebra $\bfrak_{\text{1A}}$, whilst $\gfrak[2]$ is isomorphic to $\sllie_\infty$ with the opposite one-sided Dynkin Borel subalgebra $\bfrak^-_\text{1A}$.  The global half sum of positive roots of $\tilde{\gfrak}$ is
		\begin{align}
			\tilde{\rho}=(-1,+1,-2,+2,-3,+3,\ldots)\,.
		\end{align}
	\end{exm}
	
	Let $W_\text{fl}$ be the subgroup of the Weyl group $W$ generated by the reflections with respect to roots of finite length.  Then, the  \emph{dot action}  of $W_\text{fl}$ on $\gfrak$ is given by
	\begin{align} 
		s_\alpha\cdot\lambda \defeq \lambda - \left(\lambda+\tilde{\rho}\right)\left(h_\alpha\right)\,\alpha\,,
	\end{align}
	where $s_\alpha$ is the reflection with respect to $\alpha \in \Delta_{\text{fl}}$.  If $l \in \amsbb{Z}_{\geq 0}$ and $\alpha_1,\alpha_2,\ldots,\alpha_l$ are positive roots of finite length, then
	\begin{align}
		\left(s_{\alpha_l}s_{\alpha_{l-1}}\cdots s_{\alpha_2}s_{\alpha_1}\right)\cdot\lambda\defeq s_{\alpha_l}\cdot\Bigg(s_{\alpha_{l-1}}\cdot\Big(\ldots\cdot \big(s_{\alpha_2}\cdot\left(s_{\alpha_1}\cdot\lambda\right)\big)\ldots\Big)\Bigg)\,.
	\end{align}
	 The dot action defined above is clearly a well defined group action of $W_\text{fl}$ on $\hfrak^*$.  It is easy to see that $W_{\text{fl}}$ is a (generalized) Coxeter group in the sense of Definition~\ref{def:gencox} with respect to the generating set $\left\{s_\alpha\suchthat{}\alpha\in \Sigma^+\right\}$.

	\begin{define}
		Let $\lambda \in \hfrak^*$.  We say that $\lambda$ is \emph{antidominant} if  $\lambda$ is antidominant as a weight of $\tilde{\gfrak}$ with respect to $\tilde{\bfrak}$, or equivalently, if $\left(\lambda+\tilde{\rho}\right)\left(h_\alpha\right) \notin \amsbb{Z}_{>0}$ for every positive root $\alpha$ of finite length.
	\end{define}
	
	\begin{define}
		Let $\lambda \in \hfrak^*$.   We say that $\lambda$ is \emph{almost antidominant} if  $\left(\lambda+\tilde{\rho}\right)\left(h_\alpha\right) \in \amsbb{Z}_{>0}$ for only finitely many positive roots $\alpha$ of finite length.
	\end{define}

	\begin{thm}[\bf Verma's Theorem Revisited]
		Let $\lambda\in\hfrak^*$.  For a positive root $\alpha$ of finite length such that $s_\alpha\cdot\lambda \preceq \lambda$, there exists an embedding $\Mlie\left(s_\alpha\cdot\lambda\right)\overset{\subseteq}{\longrightarrow}\Mlie(\lambda)$.
		
		\begin{pf}
			Let $u$ be a highest-weight vector of $M:=\Mlie(\lambda)$.  Consider $\tilde{M}:=\Ulie\left(\tilde{\gfrak}\right)\cdot u$, which is a Verma module with highest weight $\lambda$ of $\tilde{\gfrak}$.  Note that $\tilde\gfrak$ is locally semisimple and $\tilde{\bfrak}\subseteq\tilde{\gfrak}$ is a Dynkin Borel subalgebra of $\gfrak$.  Using Verma's Theorem (Theorem~\ref{thm:vermainf}) for the Dynkin Borel subalgebra $\tilde{\bfrak}$ of $\tilde{\gfrak}$, we conclude that $\tilde{M}$ has a Verma $\tilde{\gfrak}$-submodule $\tilde{N}$ with highest weight $\mu:=s_\alpha\cdot\lambda$.  Suppose that $v$ is a highest-weight vector of $\tilde{N}$.  Define $N:=\Ulie\left(\gfrak\right)\cdot v$.  We claim that $N$ is a Verma submodule of $M$ with highest weight $\mu$.  
			
			For a root $\gamma$, recall that $h_\gamma$ is the coroot of $\gamma$.  There exist $x_{+\gamma}\in \gfrak^{+\gamma}$ and  $x_{-\gamma}\in \gfrak^{-\gamma}$ such that $\left[x_{+\gamma},x_{-\gamma}\right]=h_{\gamma}$.
			
			First, the weight of $v$ is trivially $\mu$.  To check that $\nfrak\cdot v=0$, we only need to show that $x_\beta\cdot v=0$ for every positive root $\beta$.  Recall from the PBW Theorem that $v$ is a linear combination of elements of the form $\left(x_{-\beta_1}\cdot x_{-\beta_2}\cdot\ldots\cdot x_{-\beta_l}\right)\cdot u$ for some positive roots of finite length $\beta_1,\beta_2,\ldots,\beta_l$.  If $\beta$ is a positive root of finite length, then $x_{\beta}\cdot v=0$ since $v$ is a highest-weight vector of the $\Ulie\left(\tilde{\gfrak}\right)$-module $\tilde{N}$.   If $\beta$ is a positive root of infinite length, then we can easily prove by induction on $l$ that there exist $s\in \amsbb{Z}_{\geq 0}$ and positive roots $\beta'_1,\beta'_2,\ldots,\beta'_s$ of infinite length such that 
			\begin{align}
				\big(x_{\beta}\left(x_{-\beta_1}x_{-\beta_2}\ldots x_{-\beta_l}\right)\big)\cdot u = \left(\sum_{i=1}^s\,A_ix_{\beta'_i}\right)\cdot u\,,
			\end{align}
			where $A_1,A_2,\ldots,A_s$ are scalar multiples of subexpressions of $x_{-\beta_1} x_{-\beta_2} \ldots x_{-\beta_l}$ (which are elements of $\Ulie(\gfrak)$ of the form $x_{-\beta_{i_1}}x_{-\beta_{i_2}}\ldots x_{-\beta_{i_t}}$ where $t\in\amsbb{Z}_{\geq 0}$, $i_1,i_2,\ldots,i_t\in\{1,2,\ldots,l\}$, and $i_1<i_2<\ldots<i_t$). Therefore, 
			\begin{align}
			x_{\beta}\cdot\Big(\left(x_{-\beta_1}x_{-\beta_2}\cdots x_{-\beta_l}\right)\cdot u\Big) =\sum_{i=1}^l\,A_i\cdot\left(x_{\beta'_i}\cdot u\right)=0\,.
			\end{align} 
			The proof is now complete.
		\end{pf}
	\end{thm}
	
	\begin{define}
		A weight $\lambda\in\hfrak^*$ is \emph{of finite length} if it can be written as $\sum_{i=1}^k\,t_i\,\alpha_i$, where $t_1,t_2,\ldots,t_k\in\amsbb{K}$ and each $\alpha_i$ is a root of finite length.
	\end{define}
	
	\begin{thm}[\bf BGG Theorem Revisited]
		Let $\lambda,\mu\in\hfrak^*$ is such that $\lambda\preceq\mu$ and $\mu-\lambda$ is a weight of finite length.  Then, there exists a nontrivial embedding $\Mlie(\lambda)\overset{\subseteq}{\longrightarrow}\Mlie(\mu)$ if and only if $\lambda$ is \emph{strongly linked} to $\mu$, that is, there exist $l \in \amsbb{Z}_{\geq 0}$ and positive roots $\alpha_1,\alpha_2,\ldots,\alpha_l$ of finite length such that
		\begin{align}
			\lambda=\left(s_{\alpha_l}s_{\alpha_{l-1}}\cdots s_{\alpha_2}s_{\alpha_1}\right)\cdot\mu \preceq \left(s_{\alpha_{l-1}}\cdots s_{\alpha_2}s_{\alpha_1}\right)\cdot\mu \preceq \ldots \preceq s_{\alpha_1}\cdot \mu \preceq \mu\,.
		\end{align}
		
		\begin{pf}
			The converse is clear by Verma's Theorem above.  We only need to justify the direct implication.  
			
			To prove the direct implication, suppose that $\Mlie(\lambda)$ is a submodule of $\Mlie(\mu)$.  Write $u$ and $v$ for highest-weight vectors of $\Mlie(\lambda)$ and $\Mlie(\mu)$, respectively.   As $\mu-\lambda$ is of finite length, we can see that $\tilde{\Mlie}(\lambda):=\Ulie\left(\tilde\gfrak\right)\cdot u$ is a $\tilde\gfrak$-submodule of $\tilde{\Mlie}(\mu):=\Ulie\left(\tilde{\gfrak}\right)\cdot v$.  The BGG Theorem for Dynkin Borel subalgebras (Theorem~\ref{thm:stronglinkagegen}) applies and the claim follows immediately.
		\end{pf}
	\end{thm}
	
	\begin{cor}
		If $\Mlie(\lambda)$ is simple, then $\lambda$ is antidominant.  If $\Mlie(\lambda)$ is of finite length, then $\lambda$ is almost antidominant.
	\end{cor}
	
	\begin{rmk}
		Unfortunately, the converse of the first part of the corollary above does not hold.  Unlike in the case of Dynkin Borel subalgebras, an antidominant weight does not necessarily give rise to a simple Verma module.  
		
		Consider a dense total ordering $\prec$ on $\amsbb{Z}_{> 0}$ (i.e., an ordering $\prec$ that makes $(\amsbb{Z}_{> 0},\prec)$ isomorphic to $(\amsbb{Q},<)$, $\big(\amsbb{Q}\cup\{-\infty\},<\big)$, $\big(\amsbb{Q}\cup\{+\infty\},<\big)$, or $\big(\amsbb{Q}\cup\{-\infty,+\infty\},<\big)$ as an ordered set).  In this case, there are no roots of finite length, and we get $\tilde{\gfrak}=\hfrak$ with $\tilde{\rho}=(0,0,0,\ldots)$.  Hence, every $\lambda\in\hfrak^*$ is antidominant.  In particular, the zero weight $\textbf{0}$ is antidominant, but $\Mlie(\textbf{0})$ is not simple.  To see this, we note that there exists a nontrivial $\Ulie(\gfrak)$-module homomorphism $\Mlie(\textbf{0})\to\amsbb{K}$.  The kernel of this homomorphism is a proper submodule of $\Mlie(\textbf{0})$.
	\end{rmk}
	
	\begin{define}
		We say that a weight $\lambda \in \hfrak^*$ is \emph{primitive} if there are only finitely many positive roots $\alpha$ of infinite length such that $\lambda\left(h_\alpha\right) \in \amsbb{Z}$.
	\end{define}
	
	\begin{thm}
		Let $\lambda \in \hfrak^*$ be an antidominant weight.  If $\lambda$ is primitive, then $\Mlie(\lambda)$ is a simple $\Ulie(\gfrak)$-module.
		
		\begin{pf}
			Suppose to the contrary that $M:=\Mlie(\lambda)$ has a nonzero proper submodule $N$.  Define $M_n$ to be $\Ulie\left(\gfrak_n\right)\cdot v$ for each $n \in \amsbb{Z}_{> 0}$, where $v$ is a highest-weight vector of $M$.  Take $N_n$ to be $N\cap M_n$ for $n \in \amsbb{Z}_{> 0}$.  There exists $m\in \amsbb{Z}_{> 0}$ such that $N_n$ is a nonzero proper submodule of $M_n$ for every $n \in \amsbb{Z}_{> 0}$ with $n\geq m$.  Hence, $M_n$ has a proper Verma $\Ulie\left(\gfrak_n\right)$-submodule $\tilde{M}_n$ for every $n \in \amsbb{Z}_{> 0}$ with $n \geq m$.  We can assume that 
			\begin{align}\tilde{M}_m \subseteq \tilde{M}_{m+1}\subseteq \tilde{M}_{m+2} \subseteq \ldots\,.\end{align}
			
			The highest weight $\lambda_n:=\lambda|_{\hfrak_n}$ of $M_n$ and the highest weight $\tilde{\lambda}_n$ of $\tilde{M}_n$ must satisfy the equality: $\tilde{\lambda}_n = \lambda_n-\left(\lambda_n+\rho_n\right)\left(h_{\alpha_n}\right)\,\alpha_n$, where $\alpha_n$ is a $\bfrak_n$-positive root of $\gfrak_n$.  Note also that $\tilde{\lambda}_n \preceq \tilde{\lambda}_{n+1}|_{\hfrak_n}$ for every $n \in \amsbb{Z}_{> 0}$ such that $n \geq m$.  If $\alpha_n$ is a positive root of finite length for some integer $n\geq m$, then $s_{\alpha_n}\cdot \lambda$ restricted to $\hfrak_n$ is precisely $\tilde{\lambda}_n$.  It follows immediately that $\tilde{\lambda}_n= \left(s_{\alpha_n}\cdot\lambda\right)|_{\hfrak_n} \preceq \lambda_n$, which then means $s_{\alpha_n}\cdot \lambda \preceq \lambda$.  We now have a contradiction, as $\lambda$ is antidominant.  Hence, $\alpha_n$ is a positive root of infinite length for every $n \in \amsbb{Z}_{> 0}$ such that $n \geq m$.
			
			Because there are only finitely many positive roots $\alpha$ of infinite length such that $\lambda\left(h_\alpha\right)$ is an integer, there exists a positive root $\beta$ of infinite length such that $\alpha_n =\beta|_{\hfrak_n}$ for infinitely many integers $n \geq m$, say $n\in\left\{n_1,n_2,n_3,\ldots\right\}$ where $n_1,n_2,n_3,\ldots$ are integers such that $m \leq n_1<n_2<n_3<\ldots$.  Note that $\lambda_{n_k}\left(h_\beta\right)=\lambda\left(h_\beta\right)$ for every $k \in \amsbb{Z}_{> 0}$.  For $k \in \amsbb{Z}_{> 0}$, the condition $\tilde{\lambda}_{n_k} \preceq \left.\tilde{\lambda}_{n_{k+1}}\right|_{\hfrak_{n_k}}$ translates into $\rho_{n_k}\left(h_\beta\right)\,\alpha_{n_k} \succeq \left(\left.\rho_{n_{k+1}}\right|_{\hfrak_{n_k}}\right)\left(h_\beta\right)\,\alpha_{n_k}$, or equivalently, $\rho_{n_k}\left(h_\beta\right) \geq \rho_{n_{k+1}}\left(h_\beta\right)$ for all sufficiently large $k \in \amsbb{Z}_{> 0}$.  We now observe that $\rho_{n_k}\left(h_\beta\right)\tendsto \infty$ as $k\tendsto \infty$ (using the hypothesis that $\beta$ is a positive root of infinite length).  This is a contradiction.  Therefore, $M$ must be a simple $\Ulie(\gfrak)$-module.
		\end{pf}
	\end{thm}

	\begin{thm}
		Let $\lambda \in \hfrak^*$ be an almost antidominant weight.  If $\lambda$ is  primitive, then $\Mlie(\lambda)$ is a $\Ulie(\gfrak)$-module of finite length.
		
		\begin{pf}

			Denote by $M$ the Verma module $\Mlie(\lambda)$.  Suppose that $v$ is a highest-weight vector of $M$.  Define $M_n:=\Ulie\left(\gfrak_n\right)\cdot v$.   We claim that there exists $m\in\amsbb{Z}_{>0}$ such that each $M_n$ is a module of length at most $m$.  As $M=\lim\limits_{\underset{n}{\boldsymbol{\longrightarrow}}}$, we conclude that $M$ is also a module of length at most $m$.
			
			Since $\lambda$ is primitive, the set $A:=\big\{\alpha\in\Delta^+\suchthat{}\lambda\left(h_\alpha\right)\in\amsbb{Z}\big\}$ is finite.  Furthermore, as $\lambda$ is almost antidominant,  the set $B:=\big\{\alpha\in \Delta^+_\text{fl}\suchthat{}\left(\lambda+\tilde{\rho}\right)\left(h_\alpha\right)\in\amsbb{Z}_{>0}\big\}$ is finite.  Hence, there exists a positive integer $n_0$ such that the elements of $A\cup B$ are roots of $\gfrak_n$ with respect to $\hfrak_n$ for every integer $n\geq n_0$.
			
			Now, consider the Weyl group $W_{n_0}$ of $\gfrak_{n_0}$.  Using the finite-dimensional theory, we see that the length of the module $M_{n}$, where $n\geq n_0$, is at most 
			\begin{align}
				m:=\sum_{\substack{{x,y\in W_{n_0}}\\{x\preceq y}}}\,P_{x,y}(1)\,,
			\end{align} 
			where $P_{g_1,g_2}(q)\in\amsbb{Z}[q]$ is the Kazhdan-Lusztig polynomial of the Coxeter group $W_{n_0}$ corresponding to elements $g_1$ and $g_2$ of the Weyl group $W_{n_0}$, and $w_{n_0}^0$ is the longest element of $W_{n_0}$.
			
			Thus, the claim is justified.  That is, the Verma module $M$ is a module of finite length with length at most $m$.
		\end{pf}
	\end{thm}
	
	If a splitting Borel subalgebra $\bfrak$ is not Dynkin, then no Verma modules are in $\bbggO$.  However, it is not clear which simple quotients $\Llie(\lambda)$ lie in $\bbggO$.  While it is true that every simple module in $\bbggO$ is of the form $\Llie(\lambda)$ for some $\lambda\in\hfrak^*$, it remains an open question whether all such simple modules $\Llie(\lambda)$ are in $\bbggO$.
	
\pagebreak

\section{Kazhdan-Lusztig Theory}
\label{sec:KLTheory}

	The results in this section follow the ideas established in~\cite{Coxeter}.
	
	\subsection{Generalized Coxeter Groups}

		\begin{define}
			Let $G$ be a group with identity $1_G$.  For a (not necessarily finite) subset $X$ of $G$, we say that $G$ is a \emph{(generalized) Coxeter group} with respect to $S$ or that $(G,S)$ is \emph{(generalized) Coxeter system} iff $G$ is generated by $S$ with a presentation of the form 
			\begin{align}
				G=\left\langle S\suchthat{} (st)^{m_{s,t}}=1_G\text{ for }s,t\in S\text{ with }m_{s,t}\in\amsbb{Z}_{> 0}\right\rangle\,,
			\end{align} 
			where, for each $s,t\in S$,
			\begin{align}
				m_{s,t}=m_{t,s}
			\end{align} 
		is a positive integer or $\infty$, and, for all $s \in S$,
			\begin{align}
				m_{s,s}=1\,.
			\end{align}
		The \emph{Coxeter matrix} of $G$ is given by $\left[m_{s,t}\right]_{s,t\in S}$.  We write $\bar{S}$ for the set 
		\begin{align}
			\left\{gsg^{-1}\suchthat{}s\in S\text{ and }g\in G\right\}\,.
		\end{align}  
		Elements of $\bar{S}$ are known as \emph{reflections} in $G$.
		We say that $(G,S)$ is \emph{crystallographic} iff $m_{s,t} \in \{2,3,4,6,\infty\}$ for every $s,t\in S$ with $s \neq t$.
		\label{def:gencox}
		\end{define}
		
		\begin{thm}[\bf Universality Property]
			Let $(G,S)$ be a Coxeter system with Coxeter matrix $\left[m_{s,t}\right]_{s,t\in S}$.  For a group $\tilde{G}$ and a function $f:S\to \tilde{G}$ such that $\big(f(s)\,f(t)\big)^{m_{s,t}}=1_{\tilde{G}}$ for every $s,t\in S$ such that $m_{s,t}<\infty$, then there exists a unique extension of $f$ to a group homomorphism $f:G\to \tilde{G}$.
			\label{thm:coxuniv}
		\end{thm}
		
		\begin{cor}
			Let  $(G,S)$ be a Coxeter system with Coxeter matrix $\left[m_{s,t}\right]_{s,t\in S}$.  Then, $G$ is isomorphic as a group to $F_S/N$ where $F_S$ is the free group on $S$ and $N$ is the normal subgroup generated by $(st)^{m_{s,t}}$ for $s,t\in S$ with $m_{s,t}<\infty$.
		\end{cor}
		
		\begin{cor}
			The function $\epsilon^G:S\to \{-1,+1\}$ sending $s\mapsto -1$ for all $s\in S$ extends to a group homomorphism $\epsilon^G:G\to\{-1,+1\}$ (where $\{-1,+1\}$ is equipped with the usual multiplication).  We call $\epsilon^G$ the \emph{signature map} of $G$.
		\end{cor}
		
		\begin{define}
			The \emph{Bruhat length} $\ell^G$ of a Coxeter system $(G,S)$ is given by the function $\ell^G:G\to\amsbb{Z}_{\geq 0}$ such that, for all $g\in G$, $\ell^G(g)$ is the smallest integer $k\geq 0$ such that $g=s_1s_2\cdots s_k$ for some $s_1,s_2,\ldots,s_k\in S$.   We say that $g=s_1s_2\cdots s_k$ is a \emph{reduced expression}  for $g\in G$ if $s_1,s_2,\ldots,s_k\in S$ and $k=\ell^G(g)$.
		\end{define}
		
		\begin{prop}
			Let $(G,S)$ be a Coxeter system.    Then,
			\begin{enumerate}[(a)]
				\item $\epsilon^G(g)=(-1)^{\ell^G(g)}$ for every $g\in G$;
				\item $\ell(gh)\equiv\ell(g)+\ell(h)\pmod{2}$ for any $g,h\in G$;
				\item $\ell^G(sg)=\ell(g)\pm 1$ for all $s\in S$ and $g\in G$;
				\item $\ell^G\left(g^{-1}\right)=\ell^G(g)$ for every $g\in G$;
				\item if $\ell^G(g,h)\defeq \ell^G\left(gh^{-1}\right)$ for all $g,h\in G$, then $\ell^G(\bullet,\bullet)$ is a \emph{metric} on $G$.
			\end{enumerate}
		\end{prop}
		
		\begin{thm}[\bf Strong Exchange Property]
			Let $(G,S)$ be a Coxeter system.  For an element $g=s_1s_2\cdots s_k\in G$ with $s_1,s_2,\ldots,s_k\in S$ and $t\in \bar{S}$ such that $\ell^G(tg)<\ell^G(g)$, then $t=s_1\cdots s_{i-1}s_is_{i-1}\cdots s_1$ and $tg=s_1\cdots\slashed{s}_i\cdots s_k$ for some $i\in [k]$, where every slashed term is disregarded from the product. (Here, $[n]$ denotes $\{1,2,\ldots,n\}$ for each positive integer $n$.  Furthermore, $[0]$ denotes $\emptyset$.)  The index $i$ is unique of $g=s_1s_2\cdots s_k$ is a reduced expression for $g$.
			
			\begin{pf}
				Let $\hat{S}$ denote $\bar{S}\times\{-1,+1\}$.  For each $s\in S$, define $\varpi_s:\hat{S}\to\hat{S}$ to be the map $(t,e)\mapsto \left(sts,e\,\eta(s,t)\right)$ for each $t\in \bar{S}$ and $e\in\{-1,+1\}$, where 
				\begin{align}	
					\eta(s,t):=\left\{\begin{array}{ll}
						-1\,,&\text{if }t=s\,,
						\\
						+1\,,&\text{if }t\neq s\,.
					\end{array}\right.
				\end{align}
			We extend the definition of $\eta$ to $G\times\bar{S}$ as follows: if $g=s_1s_2\cdots s_k\in G$ with $s_1,s_2,\ldots,s_k\in S$, then
			\begin{align}
				\eta(g,t):=\prod_{j=1}^k\,\eta\left(s_j,s_{j-1}\cdots s_1ts_1\cdots s_{j-1}\right)\,.
			\end{align}
			Note that $\eta$ is well defined.	
				
			For a set $X$, $\text{Perm}(X)$ denotes the group of permutations on $X$.  Let $\left[m_{s,t}\right]_{s,t\in S}$ be the Coxeter matrix of $(G,S)$.  It can be easily shown that the function $\varpi:S\to \text{Perm}\left(\hat{S}\right)$ sending $s\mapsto \varpi_s$ for every $s\in S$ satisfies $\left(\varpi_s\circ \varpi_t\right)^{m_{s,t}}=\text{id}_{\hat{S}}$ for every $s,t\in S$ such that $m_{s,t}<\infty$.  Therefore, by the Universality Property (Theorem~\ref{thm:coxuniv}), $\varpi$ extends to a group homomorphism $\varpi:G\to \text{Perm}\left(\hat{S}\right)$.
			
			It follows immediately that, for any $g\in G$, $t\in\bar{S}$, and $e\in\{-1,+1\}$, we have
			\begin{align}
				\varpi_{g}(t,e) =\Big(gtg^{-1},e\,\eta\left(g^{-1},t\right)\Big)\,.
			\end{align}
			Now, we assume that $g=s_1s_2\cdots s_k\in G$ with $s_1,s_2,\ldots,s_k\in S$ and $t\in \bar{S}$.  We claim that $\ell^G(tg)<\ell^G(g)$ if and only if $\eta(g,t)=-1$.  
			
			If $\eta(g,t)=-1$, then $\eta\left(s_i,s_{i-1}\cdots s_1ts_i\cdots s_{i-1}\right)=-1$ for some $i\in[k]$.  Hence, we obtain $t=s_1\cdots s_{i-1}s_is_{i-1}\cdots s_1$, so that $tg=s_1\cdots\slashed{s}_i\cdots s_k$, leading to $\ell^G(tg)<\ell^G(g)$.  The uniqueness of $i$ is trivial if $k=\ell^G(g)$.  
			
			If $\eta(g,t)=+1$, then
			\begin{align}
				\varpi_{(tg)^{-1}}(t,e)=\varpi_{g^{-1}}\big(\varpi_t(t,e)\big)=\varpi_{g^{-1}}(t,-e)=\left(g^{-1}tg,-e\,\eta(g,t)\right)=\left(g^{-1}tg,-e\right)\,.
			\end{align}
			That is, $\eta(tg,t)=-1$.  Consequently, $\ell^G(g)=\ell^G\big(t(tg)\big)<\ell^G(tg)$.
			\end{pf}
		\end{thm}
		
\subsection{Bruhat Ordering}
		
		\begin{define}
			Let $(G,S)$ be a Coxeter system and $g,h\in G$.  Then, we write $g\underset{G}{\preccurlyeq } h$ if there is a reduced expression $h=s_1s_2\cdots s_k$ where $k\in\amsbb{Z}_{\geq 0}$ and $s_1,s_2,\ldots,s_k\in S$ such that $g$ is a \emph{subword} of $s_1s_2\cdots s_k$, namely, there exist $j\in\amsbb{Z}_{\geq 0}$ and integers $i_1,i_2,\ldots,i_j$ with $1\leq i_1<i_2<\ldots <i_j\leq k$ such that $g=s_{i_1}s_{i_2}\cdots s_{i_j}$.  The relation $\underset{G}{\preccurlyeq}$ is called the \emph{Bruhat order} on $G$.
		\end{define}
		
		\begin{prop}
			Let $(G,S)$ be a Coxeter system and $g,h\in G$.
			\begin{enumerate}[(a)]
				\item The condition $g\underset{G}{\preccurlyeq} h$ implies that $\ell^G(g)\leq \ell^G(h)$.
				\item If $g\underset{G}{\preccurlyeq} h$ and $\ell^G(g)< \ell^G(h)$, then there exists an element $\tilde{g}\in G$ such that $g\underset{G}{\preccurlyeq} \tilde{g}\underset{G}{\preccurlyeq} h$ and $\ell^G(g)+1=\ell^G\left(\tilde{g}\right)$.
			\end{enumerate}
			
			\begin{pf}
				Part (a) is trivial.  To prove Part (b), let $h=s_1s_2\cdots s_k$ with $s_1,s_2,\ldots,s_k\in S$ be a reduced expression for $h$ such that $g$ is a subword of $s_1s_2\cdots s_k$.  Suppose that $g$ is obtained by removing $s_{i_1}$, $s_{i_2}$, $\ldots$, $s_{i_q}$ from the expression $s_1s_2\cdots s_k$, where $1\leq i_1<i_2<\ldots<i_q\leq k$ and $i_1$ is chosen to be the largest possible.  Take $t:=s_1\cdots s_{i_1-1}s_{i_1}s_{i_1-1}\cdots s_1 \in \bar{S}$.  Then, $tg=s_1\cdots s_{i_1}\cdots \slashed{s}_{i_2}\cdots \slashed{s}_{i_q}\cdots s_k$.  Hence, $\ell^G(tg)\leq \ell^G(g)+1$.
				
				If $\ell^G(tg)<\ell^G(g)+1$, then by the Strong Exchange Property, then 
				\begin{align}
					t=s_1\cdots \slashed{s}_{i_1}\cdots\slashed{s}_{i_j} \cdots  s_{r-1}s_rs_{r-1}\cdots \slashed{s}_{i_j} \cdots \slashed{s}_{i_1}\cdots s_1
				\end{align} 
				for some $i_{j-1}<r<i_j$ with $i_0=0$ and $j\in[q]$.
				
				If $j=1$, then we have $h=t^2h=s_1\cdots \slashed{s}_{r}\cdots\slashed{s}_{i_1}\cdots s_k$ contradicting the assumption that $h=s_1s_2\cdots s_k$ is a reduced expression for $h$.  If $j>1$, then 
				\begin{align}
					g=t^2g=s_1\cdots s_{i_1} \cdots \slashed{s}_{i_2}\cdots \slashed{s}_r \cdots \slashed{s}_{i_q}\cdots s_k\,,
				\end{align} 
				which contradicts the maximality of $i_1$.
			\end{pf}
		\end{prop}
		
		\begin{thm}
			Let $(G,S)$ be a Coxeter system.  For $g,h\in G$, we write $g\overset{t}{\rightharpoondown} h$, where $t\in\bar{S}$, if $tg=h$ and $\ell^G(g)<\ell^G(h)$.  We write $g\rightharpoondown h$ if there exists $t\in\bar{S}$ such that $g\overset{t}{\rightharpoondown} h$.    Then, $g\underset{G}{\preccurlyeq} h$ if and only if there exists $k\in\amsbb{Z}_{\geq 0}$ and $u_0,u_1,\ldots,u_k\in G$ such that $g=u_0\rightharpoondown u_1\rightharpoondown u_2 \rightharpoondown \ldots \rightharpoondown u_{k-1} \rightharpoondown u_k=h$.
			
			\begin{pf}
				Suppose that $g\underset{G}{\preccurlyeq} h$.  Then, we can easily prove the direct implication by induction on $\ell^G(h)-\ell^G(g)$, with the application of the previous proposition. 
				
				Conversely, assume that there exists $k\in\amsbb{Z}_{\geq 0}$ and $u_0,u_1,\ldots,u_k\in G$ such that the chain $g=u_0\rightharpoondown u_1\rightharpoondown u_2 \rightharpoondown \ldots \rightharpoondown u_{k-1} \rightharpoondown u_k=h$ holds.  Then, it follows immediately from the Strong Exchange Property that $g\underset{G}{\preccurlyeq} h$.
			\end{pf}
		\end{thm}
		
		\begin{cor}
			Let $(G,S)$ be a Coxeter system.  For $g,h\in G$, the following statements are equivalent:
			\begin{enumerate}[(i)]
				\item $g\underset{G}{\preccurlyeq} h$;
				\item some reduced expression for $h$ has a subword that is a reduced expression for $g$;
				\item every reduced expression for $h$ has a subword that is a reduced expression for $g$.
			\end{enumerate}
		\end{cor}
		
		\begin{cor}
			Let $(G,S)$ be a Coxeter system.
			\begin{enumerate}[(a)]
				\item The Bruhat order on $G$ is indeed a partial order on $G$.
				\item For $g,h\in G$, the \emph{Bruhat interval} $[g,h]^G\defeq \left\{w\in G\suchthat{}g\underset{G}{\preccurlyeq} w \underset{G}{\preccurlyeq} h\right\}$ is finite with at most $2^{\ell^G(h)}$ elements.
				\item The Bruhat order on $G$ is locally finite.
			\end{enumerate}
		\end{cor}

		Note that the Weyl group $W$ is a crystallographic Coxeter group with respect to the simple reflections.  Hence, it is equipped with a Bruhat ordering $\preccurlyeq$ with the Bruhat length function $\ell$.  We shall write 
		\begin{align}
			\ell(x,y) \defeq \ell(y)-\ell(x)
		\end{align}
	for $x,y \in W$ with $x \preccurlyeq y$.  Below is the list of some properties of this Bruhat ordering on $W$.
	\begin{enumerate}
		\item Each $w \in W$ is determined by the set of $\alpha \in \Delta^+$ for which $w\alpha \succ 0$.  Furthermore, $\ell(w)$ is precisely the cardinality of this set.
		\item For $w \in W$, $\ell(w)=\ell\left(w^{-1}\right)$.  Hence, $\ell(w) = \Big|\Delta^+\cap w\left(\Delta^-\right)\Big|$.
		\item If $\alpha \in \Delta^+$ and $w \in W$ satisfy $\ell\left(w s_\alpha\right) > \ell(w)$, then $w \alpha \succ 0$.  On the other hand, $\ell\left(w s_\alpha\right) < \ell(w)$ implies $w\alpha\prec 0$.  Consequently, $\ell\left(w s_\alpha \right) > \ell(w)$ if and only if $w\alpha \succ 0$.
		\item Similarly, for $\alpha \in \Delta^+$ and $w \in W$, $\ell\left(s_\alpha w\right) > \ell(w)$ if and only if $w^{-1}\alpha \succ 0$.
	\end{enumerate}
	
\subsection{Parabolic Subgroups of a Coxeter Group}
	
	\begin{define}
		Let $(G,S)$ be a Coxeter system.  For $J\subseteq S$, the \emph{parabolic subgroup} $G_J$ of $G$ is the subgroup of $G$ generated by $J$.
	\end{define}
	
	\begin{thm}
		Let $(G,S)$ be a Coxeter system and $J\subseteq S$.  Then, $\left(G_J,J\right)$ is also a Coxeter system.
	\end{thm}
	
	\begin{thm}
		Let $(G,S)$ be a Coxeter system and $\tilde{G}$ a parabolic subgroup of $G$.  Then, the following statements hold:  
		\begin{enumerate}[(a)]
			\item The Bruhat length of $\tilde{G}$ coincides with the restriction of the Bruhat length of $G$ on $\tilde{G}$ (i.e., $\ell^{\tilde{G}}=\ell^G|_{\tilde{G}}$);
			\item For $g\in G$ and $h\in \tilde{G}$, $g\underset{\tilde{G}}{\preccurlyeq} h$ implies $g\in\tilde{G}$;
			\item For $g,h\in\tilde{G}$, $g\underset{\tilde{G}}{\preccurlyeq} h$ if and only if $g\underset{{G}}{\preccurlyeq} h$;
			\item For $g,h\in\tilde{G}$, $[g,h]^{\tilde{G}}=[g,h]^G$.
		\end{enumerate}
	\end{thm}
	
	The theorem above allows us to omit the subscripts and the superscripts when discussing terms such as the Bruhat length, the Bruhat order, and the Bruhat interval of a Coxeter group.  That is, if we use notations such as $\ell$, $\preccurlyeq$, or $[\bullet,\bullet]$, there is no ambiguity whether the notations are restricted to only a parabolic subgroup of a Coxeter group, or to the whole group.

\subsection{Hecke Algebras}
	
This subsection is based on the work~\cite{KL} by Kazhdan and Lusztig.

	\begin{define}
		Let $(G,S)$ be a generalized Coxeter system and $q$ an indeterminate.   The ring $\amsbb{Z}\left[q^{-\frac{1}{2}},q^{+\frac{1}{2}}\right]$ of Laurent polynomials in $q^{\frac{1}{2}}$ is denoted by $\mathcal{A}$.  The \emph{Hecke algebra} $\mathcal{H}$ is an associative algebra which is a free module over $\mathcal{A}$ with the generating set $\left\{T_g\suchthat{}g\in G\right\}$ such that the multiplicative identity of $\mathcal{H}$ is $1_{\mathcal{H}}=T_{1_G}$ and that the following multiplicative relations are satisfied:
		\begin{align}
			T_s^2 = (q-1)\,T_s + q\,T_{1_G}\,,
		\end{align}
		\begin{align}
			T_g\,T_s = T_{gs} \,\text{   if }g \prec gs\,,
		\end{align}
		and
		\begin{align}
			T_s\,T_g = T_{sg} \,\text{   if }g \prec sg\,,
		\end{align}
		for each $s \in S$ and $g \in G$.
	\end{define}
	
	  \begin{prop} 
	  	Let $\mathcal{H}$ be the Hecke algebra of a given Coxeter system $(G,S)$ with indeterminate $q$.  
	  	\begin{itemize}
	  		\item[(a)] For each $g \in G$, $T_g$ is invertible, as $T_s^{-1}= q^{-1}\,T_s+\left(q^{-1}-1\right)\,T_{1_G}$.  In general, for $g \in G$, there are polynomials $R^G_{x,g}(q) \in \amsbb{Z}[q]$ such that
		\begin{align}
			T_g^{-1}= q^{-\ell(g)}\sum_{x \preccurlyeq g}\,(-1)^{\ell(x,g)}\,R_{x,g}(q)\,T_x\,,
		\end{align}  
		for all $s \in S$.
			\item[(b)] Let $g \in G$ and $s \in S$.  If $g \succ gs$, then we have $T_gT_s=q\,T_{gs}+(q-1)\,T_g$.  If $g \succ sg$, then $T_sT_g=q\,T_{sg}+(q-1)T_g$.
			\item[(c)] There exists an \emph{involution} $\iota:\mathcal{H}\to\mathcal{H}$ sending $q^{+\frac{1}{2}}\mapsto q^{-\frac{1}{2}}$ and $T_g \mapsto T_{g^{-1}}^{-1}$ for all $g \in G$.  It is customary to write $\bar{X}$ for $\iota(X)$, where $X\in \mathcal{H}$.
		\end{itemize}
		
		\begin{pf}
			For each $g\in G$, there exists a finite subset $J$ of $S$ such that $g\in G_J$.  Applying the usual theory on the Coxeter group $G_J$ with finite generating set $J$, we can prove this proposition easily.
		\end{pf}
	\end{prop}
	
	\begin{thm}
		The \emph{$R$-polynomials} for a Coxeter system $(G,S)$ in the previous proposition satisfy the conditions below.
		\begin{itemize}
			\item[(a)]  If $x,y \in G$ with $x \not\preccurlyeq y$, we have $R^G_{x,y}(q)=0$.
			\item[(b)] For every $x \in G$, $R^G_{x,x}(q)=1$.
			\item[(c)] For $x,y \in G$ and $s \in S$ such that $ys \prec y$, we have
				\begin{align}
					R^G_{x,y}(q) = \left\{
					\begin{array}{ll}
						R^G_{xs,ys}(q)	&	, \text{if }xs \prec x\,,
						\\
						q\,R^G_{xs,ys}(q)+(q-1)\,R^G_{x,ys}(q)	&	, \text{if }xs \not\prec x\,.
					\end{array}
					\right.
				\end{align}
		\end{itemize}
	\end{thm}
	
	\begin{prop}
		Here are some properties of the $R$-polynomial $R^G_{x,y}$ for a Coxeter system $(G,S)$, where $x,y \in G$ are such that $x \preccurlyeq y$.
		\begin{itemize}
			\item[(a)] $R^G_{x,y}(q)$ is monic in $q$ of degree $\ell(x,y)$.
			\item[(b)] $R^G_{x,y}(0)=(-1)^{\ell(x,y)}$.
			\item[(c)] $R^G_{x,y}(1) = 0$, provided that $x\neq y$.
			\item[(d)] $q^{\ell(x,y)}\,R^G_{x,y}\left(\frac{1}{q}\right)=(-1)^{\ell(x,y)}\,R^G_{x,y}(q)$.
			\item[(e)] $\sum_{w\in[x,y]}\,(-1)^{\ell(x,y)}\,R^G_{x,w}(q)\,R_{w,y}^G(q)=\deltaup_{x,y}$.
		\end{itemize}
	\end{prop}
	
	\begin{thm}
		For a Coxeter system $(G,S)$, the Bruhat order on $G$ makes $G$ a graded partially ordered set with the Bruhat length as the rank function.  Furthermore, the \emph{M\"{o}bius function} with respect to the Bruhat order on $G$ is given by $\muup(x,y)\defeq(-1)^{\ell(x,y)}$ for all $x,y\in G$ with $x \preccurlyeq y$.
	\end{thm}

\subsection{Kazhdan-Lusztig Polynomials}

\begin{thm}[\bf Kazhdan-Lusztig]
		Let $(G,S)$ be a Coxeter system with the associated Hecke algebra $\mathcal{H}$.  For each $g \in G$, there is a unique $C_g \in \mathcal{H}$ fixed by the involution on $\mathcal{H}$ and satisfying the condition
		\begin{align}
			C_g = q^{-\frac{\ell(g)}{2}}\, \sum_{x \preccurlyeq g}\,(-1)^{\ell(x,g)}q^{\ell(x,g)}\,\overline{P^G_{x,g}(q)}\,T_x\,,
		\end{align}
		where, for each $x,y \in G$,
		\begin{itemize}
			\item[(i)] $P^G_{x,y}(q) \in \amsbb{Z}[q]$,
			\item[(ii)] $P^G_{x,y}(q)=0$ if $x \not\preccurlyeq y$,
			\item[(iii)] $P^G_{x,x}(q)=1$, and
			\item[(iv)] $\deg\left(P^G_{x,y}(q)\right) \leq \frac{\ell(x,y)-1}{2}$ if $x \prec y$.
		\end{itemize}
		The polynomials $P^G_{x,y}(q)$, where $x,y \in G$, are known as the \emph{Kazhdan-Lusztig (KL) polynomials} of $G$.
	\end{thm}
	
	\begin{thm}
		Let $(G,S)$ be a Coxeter system.  Then, for all $x,y\in G$ such that $x \preccurlyeq y$, we have
		\begin{align}
		q^{\ell(x,y)}\,P^G_{x,y}\left(\frac{1}{q}\right) = \sum_{a \in [x,y]}\,R^G_{x,a}(q)\,P^G_{a,y}(q)\,.
		\end{align}
	\end{thm}
	
		\begin{thm}
		Let $(G,S)$ be a Coxeter system.  For $u,v \in G$ such that $u \preccurlyeq v$, let $\varkappa^G(u,v)$ be the coefficient of $q^{\frac{\ell(u,v)-1}{2}}$ in $P^G_{u,v}(q)$.  Now, for a given pair $x,y\in G$ such that $x \preccurlyeq y$ and for $s \in S$ such that $ys \prec y$, we have
		\begin{align}
			P^G_{x,y}(q) = \left\{
				\begin{array}{ll}
					P^G_{xs,ys}(q)+q\,P^G_{x,ys}(q)-\sum_{\substack{{z \in [x,ys]}\\{zs \prec z}}}\,q^{\frac{\ell(z,y)}{2}}\,\varkappa^G(z,ys)\,P^G_{x,z}(q)	&	, \text{if }xs\prec x\,,
					\\
					q\,P^G_{xs,ys}(q)+P^G_{x,ys}(q)-\sum_{\substack{{z \in [x,ys]}\\{zs \prec z}}}\,q^{\frac{\ell(z,y)}{2}}\,\varkappa^G(z,ys)\,P^G_{x,z}(q)	&	, \text{if }xs\not\prec x\,.
				\end{array}
			\right.
		\end{align}
	\end{thm}
	
	\begin{prop}
		For a Coxeter system $(G,S)$, let $x,y \in G$ be such that $x \preccurlyeq y$.
		\begin{itemize}
			\item[(a)]	If $\ell(x,y) \leq 2$, then $P^G_{x,y}(q)=1$.
			\item[(b)] 	$P^G_{x,y}(0)=1$.
			\item[(c)]	If $s \in S$ is such that $ys \prec y$, then $P^G_{xs,y}(q) = P^G_{x,y}(q)$.
			\item[(d)]	If $G$ is finite and $g_{\max}$ is the longest element of $G$, then $P_{x,g_{\max}}(q)=1$.
			\item[(e)] 	If $G$ is a finite dihedral group (such as the Weyl groups in types $\textsf{A}_2$, $\textsf{B}_2$, and $\textsf{G}_2$), then $P^G_{x,y}(q)=1$.
		\end{itemize}
	\end{prop}
	
	In 2004, Soergel proved that the coefficients of the Kazhdan-Lusztig polynomials are nonnegative, for the case of finite generating sets (see~\cite{KLnonnegative}).  It follows immediately that the coefficients of Kazhdan-Lusztig polynomials of generalized Coxeter groups are also nonnegative.

	\begin{exm}
		For a Coxeter system $(G,S)$, we have $C_{1_G}=T_{1_G}$ and 
		\begin{align}
			C_{s}=q^{-\frac{1}{2}}\left(T_s - q\,T_{1_G}\right)
		\end{align} 
		for each $s \in S$.
	\end{exm}
	
	\begin{thm}
		For any Coxeter group $G$ and an arbitrary parabolic subgroup $\tilde{G}$, it holds that $R^G_{x,y}=R^{\tilde{G}}_{x,y}$ and $P^G_{x,y}=P^{\tilde{G}}_{x,y}$ for every $x,y\in \tilde{G}$ with $x\preccurlyeq y$.  That is, there is no ambiguity if the superscripts $G$ and $\tilde{G}$ are dropped.
		\label{thm:KLpolystab}
	\end{thm}

\pagebreak

\section{Block Decomposition and Kazhdan Lusztig Multiplicities}
\label{ch:blocks}

In this section, $\gfrak$ is a root-reductive Lie algebra and $\bfrak$ is a Dynkin Borel subalgebra.

\subsection{Block Decomposition}

Let $M$ be an indecomposable object of $\bbggO$.  We shall construct an ordered countable set $\Gamma(M)=\left(v_0,v_1,v_2,\ldots\right)$ such that $\Gamma(M)\subseteq M$ consists of weight vectors of $M$ which generate $M$ as a $\Ulie(\gfrak)$-module, and the set $\Gamma(M)$ has certain desirable properties.  If $M=0$, then we set $\Gamma(M):=(0,0,0,0,0,\ldots)$.

For $M\neq 0$, we let $u\neq 0$ be a singular vector of $M$, and $\xi\in\hfrak^*$ the weight of $u$.  Let $[\xi]$ denote the set of all weights $\zeta\in\hfrak^*$ such that $\zeta-\xi$ is in the $\amsbb{Z}$-span of the simple ($\bfrak$-positive) roots of $\gfrak$.  For a weight $\zeta\in[\xi]$, the \emph{distance} between $\zeta$ and $\xi$, denoted by $\text{dist}(\zeta,\xi)$, is defined to be the sum $\sum_{i=1}^k\,\left|t_i\right|$, if $\zeta-\xi=\sum_{i=1}^k\,t_i\,\alpha_i$, where $\alpha_1,\alpha_2,\ldots,\alpha_k$ are ($\bfrak$-positive) simple roots of $\gfrak$ and $t_1,t_2,\ldots,t_k\in\amsbb{Z}$.  Furthermore, the height of $\zeta-\xi$, denoted by $\text{ht}(\zeta-\xi)$, is the smallest integer $n\geq 0$ such that $h_{\alpha_1},h_{\alpha_2},\ldots,h_{\alpha_k}$ are all in $\gfrak_n$ (noting that $\text{ht}(\zeta-\xi)=0$ if and only if $\zeta =\xi$).  

We start with $m:=0$; then we set $d(0):=0$ and $v_0:=u$. Now, for $m>0$, suppose that the value $d(m-1)$ is known and that the vectors $v_0,v_1,\ldots,v_{d(m-1)}$ have been defined.  The set $\pazocal{S}_m$ of weights $\zeta \in [\xi]$ such that $\text{dist}(\zeta,\xi)\leq m$ and $\text{ht}(\zeta-\xi)\leq m$ is a finite set.  Let $\pazocal{V}^1_m$  denote the $\amsbb{K}$-span of all weight vectors $v\in M$ with weights in $\pazocal{S}_m$ such that $\nfrak\cdot v=0$.  Let $u^1_1,u^1_2,\ldots,u^1_{l_1}$ be weight vectors of $M$ which form a $\amsbb{K}$-basis of $\pazocal{V}^1_m$.  

Assume that the collections $\left(u^1_j\right)_{j=1}^{l_1}$, $\left(u^2_j\right)_{j=1}^{l_2}$, $\ldots$, $\left(u^{r}_j\right)_{j=1}^{l_{r}}$ of weight vectors of $M$ have been obtained. Consider the module 
\begin{align}
	\underline{M}(m,r):=M/\left(\sum_{i=0}^{d(m-1)}\,\Ulie(\gfrak)\cdot v_i+ \sum_{p=1}^r\,\sum_{j=1}^{l_p}\,\Ulie(\gfrak)\cdot u^p_j\right)\,.
\end{align}
Let $\pazocal{V}^{r+1}_m$ denote the $\amsbb{K}$-span of all weight vectors $v\in \underline{M}(m,r)$ with weights in $\pazocal{S}_m$ such that $\nfrak\cdot v=0$.  Take $\tilde{u}^{r+1}_1,\tilde{u}^{r+1}_2,\ldots,\tilde{u}^{r+1}_{l_{r+1}}$ to be weight vectors of $M$ which form a $\amsbb{K}$-basis of $\pazocal{V}^{r+1}_m$.  Now, there exist weight vectors $u^{r+1}_1,u^{r+1}_2,\ldots,u^{r+1}_{l_{r+1}}$, whose respective images under the projection $M\to \underline{M}(m,r)$ are $\tilde{u}^{r+1}_1,\tilde{u}^{r+1}_2,\ldots,\tilde{u}^{r+1}_{l_{r+1}}$.

The process in the previous paragraph must end with $\pazocal{V}_m^{\bar{r}(m)+1}=0$ for some nonnegative integer $\bar{r}(m)$ because the vector subspace of $M$ spanned by the weight vectors with weights in $\pazocal{S}_m$ is finite dimensional.  Then, we take 
\begin{align}
	d(m):=d(m-1)+\sum_{p=1}^{\bar{r}(m)}\,l_p\,,
\end{align}
and
\begin{align}
	v_{d(m-1)+j}:=u^{l_p}_{j-\sum_{\tilde{p}=1}^{p-1}\,l_{\tilde{p}}}
\end{align}
if $\sum_{\tilde{p}=1}^{p-1}\,l_{\tilde{p}}<j\leq \sum_{\tilde{p}=1}^p\,l_{\tilde{p}}$.

Note that $d(m)>d(m-1)$ for every $m=1,2,\ldots$ because $\pazocal{V}_m^1$ always contains $u$.  When $M$ is a $\gfrak$-module of finite length, it is possible that $\Gamma(M)$ is eventually periodic (that is, there exist positive integers $n_0$ and $n_1$ such that $v_{n}=v_{n+n_1}$ for every integer $n\geq n_0$).  In particular, the ordered set $\Gamma(M)$ may take the form $(u,u,u,u,\ldots)$ when $M$ is a highest-weight module with $u$ as a highest-weight vector. 

 We claim that the ordered set $\Gamma(M):=\left(v_0,v_1,v_2,v_3,\ldots\right)$ generate $M$ as a $\Ulie(\gfrak)$-module.  For a fixed weight $\zeta$ of $M$, consider the vector subspace $M^\zeta$.  

Let $\pazocal{T}_\zeta$ denote the set of all weights $\tilde{\zeta}$ of the $\Ulie(\gfrak)$-module $\Ulie(\gfrak)\cdot M^\zeta$ which satisfy $\tilde{\zeta}\succeq \zeta$.  Note that $\pazocal{T}_\zeta$ is finite as $\nfrak$ acts locally finitely on $M^\zeta$.  Let $m_\zeta$ denote the maximum value of the two numbers
	\begin{align}
		\max\big\{\text{dist}(\tilde{\zeta},\xi)\suchthat{}\tilde{\zeta} \in \pazocal{T}_\zeta\big\}\text{ and }\max\big\{\text{ht}(\tilde{\zeta}-\xi)\suchthat{}\tilde{\zeta} \in \pazocal{T}_\zeta\big\}\,.
	\end{align}
Then, in the $m_\zeta$-th step of the procedure (from which $d\left(m_\zeta\right)$ is obtained), the $\Ulie(\gfrak)$-module
\begin{align}
	\underline{M}\Big(m_\zeta,\bar{r}\left(m_\zeta\right)\Big)=M/\left(\sum_{i=0}^{d(m-1)}\,\Ulie(\gfrak)\cdot v_i + \sum_{p=1}^{\bar{r}(m)}\,\sum_{j=1}^{l_p}\,\Ulie(\gfrak)\cdot u^p_j\right)
\end{align}
cannot have $\Llie(\zeta)$ as a composition factor.  To elaborate, if such a composition factor exists, it must arise from the weight space of $\underline{M}\Big(m_\zeta,\bar{r}\left(m_\zeta\right)\Big)$ with weight $\zeta$.  However, all composition factors of $M$ isomorphic to $\Llie(\zeta)$ come from subquotients of  $\Ulie(\gfrak)\cdot M^\zeta$, and by the definition of $m_\zeta$, the image of $\Ulie(\gfrak)\cdot M^\zeta$ under the canonical projection $M\to \underline{M}\Big(m_\zeta,\bar{r}\left(m_\zeta\right)\Big)$ has no composition factors isomorphic to $\Llie(\zeta)$.

Therefore, for every composition factor $\Llie(\zeta)$ of $M$, it is exhausted in the quotient module $\underline{M}\Big(m_\zeta,\bar{r}\left(m_\zeta\right)\Big)$.  Thus, the sum $\sum_{i=0}^\infty\,\Ulie\left(\gfrak\right)\cdot v_i$ is precisely $M$ itself.  

It is important to note the following properties of the set $\Gamma(M)$.  For every $n\in\amsbb{Z}_{>0}$, let $\Xi_n\subseteq \hfrak^*$ be the support (as a semisimple $\hfrak$-module) of the $\gfrak'_n$-module $M_n:=\sum_{i=0}^n\,\Ulie\left(\gfrak'_n\right)\cdot v_i$, where $\gfrak_n':=\hfrak+\gfrak_n$.   Then, for every $\xi\in \Xi_n$ and for any integer $\tilde{n}\geq n$, we have
	\begin{align}
		\dim_\amsbb{K}\left(M_n^{\xi}\right)=\dim_\amsbb{K}\Big(\big(\Ulie\left(\gfrak'_{\tilde{n}}\right)\cdot M_n\big)^{\xi}\Big)\,.
		\label{eq:equaldim}
	\end{align}
	This is because our construction of $\Gamma(M)$ ensures that $\Ulie\left(\gfrak'_{\tilde{n}}\right)\cdot M_n = \Ulie\left(\nfrak_{\tilde{n}}^-\right)\cdot M_n$.  Note that, if $x$ is in a positive root space of $\gfrak'_{\tilde{n}}$ that is not in $\gfrak'_n$, then $x\cdot M_n=0$.

	\begin{define}
		Let $\lambda\in\hfrak^*$.  Define $\bbggO\llbracket\lambda\rrbracket$ to be the full subcategory of $\bbggO$ consisting of modules $M$ whose composition factors are of the form $\Llie(\mu)$ with $\mu\in\llbracket\lambda\rrbracket$, where $\llbracket\lambda\rrbracket$ is the integral Weyl dot-orbit $\llbracket\lambda\rrbracket:=W[\lambda]\cdot\lambda$.
	\end{define}
	
	\begin{prop}
		Let $M\in\bbggO$ be indecomposable and $\lambda\in\hfrak^*$ be such that $\Llie(\lambda)$ is a composition factor of $M$.  Then, all composition factors of $M$ are of the form $\Llie(\mu)$ for some $\mu\in \llbracket\lambda\rrbracket$.
		\label{prop:blocks}
		
		\begin{pf}
			Let $\lambda$ and $\mu$ be on different integral Weyl dot-orbits.  Suppose there exists an indecomposable $M\in\bar\bggO$ with $\Llie(\lambda)$ and $\Llie(\mu)$ as composition factors.  Since $M\neq0$, we can apply the algorithm discussed earlier in this subsection and obtain an ordered set $\Gamma(M)=\left(v_0,v_1,v_2,\ldots\right)$ which generates $M$ as a $\Ulie(\gfrak)$-module.  For every $n\in\amsbb{Z}_{>0}$, let $\gfrak'_n$ denote the subalgebra $\hfrak+\gfrak_n$ and set $\bfrak_n'\subseteq\gfrak_n'$ to be the Borel subalgebra $\hfrak+\bfrak_n$ of $\gfrak_n'$.  Then, the $\Ulie\left(\gfrak'_n\right)$-module  $M_n$ is given by
			\begin{align}
				M_n:=\sum_{i=0}^n\,\Ulie\left(\gfrak'_n\right)\cdot v_i\,.
			\end{align}
	
	Note that the finite-dimensional theory carries trivially over to $\gfrak_n'$ (see Remark~\ref{rmk:reductive}), and we use the notation $\bggO^{\gfrak'_n}_{\bfrak'_n}$ for the category $\bggO$ of $\gfrak'_n$ with respect to the Borel subalgebra $\bfrak_n'$.  Denote by $\bggO^{\gfrak'_n}_{\bfrak'_n}\left\llbracket\lambda\right\rrbracket$ the block of $\bggO^{\gfrak'_n}_{\bfrak'_n}$ containing $\Llie_n\left(\lambda\right):=\Llie\left(\lambda;\gfrak'_n,\bfrak'_n,\hfrak\right)$.  Note that we have the direct sum decomposition $M_n=X_n\oplus Y_n$, where $X_n\in \bggO^{\gfrak_n}_{\bfrak_n}\left\llbracket\lambda\right\rrbracket$ and all composition factors of $Y_n$ are not in $ \bggO^{\gfrak'_n}_{\bfrak'_n}\left\llbracket\lambda\right\rrbracket$.  The submodules $X_n$ and $Y_n$ are unique.  Furthermore, if $N_n$ is an indecomposable submodule of $M_n$, then $N_n$ must lie entirely in $X_n$ or in $Y_n$.
	
	 Define 
	\begin{align} X:=\big\{x\in M\,\boldsymbol{|}\,x\in X_n\text{ for all sufficiently large }n\big\}\end{align}
	and
	\begin{align}
		Y:=\text{span}_\amsbb{K}\big\{y\in M\,\boldsymbol{|}\,y\in Y_n\text{ for infinitely many }n\big\}\,.
	\end{align}
	Then, it is evident that $X$ and $Y$ are $\gfrak$-submodules of $M$.  We shall  prove  that $X+Y=M$ and that $X\cap Y=0$.  
	
	First, let the $\Ulie\left(\gfrak'_n\right)$-module $X_n'$ be an indecomposable direct summand of $X_n$.  Fix $\tilde{n}\geq n$.   Note that we have either $X'_n\subseteq X_{\tilde{n}}$ or $X'_n\subseteq Y_{\tilde{n}}$.   Likewise, if $Y_n'$ is an indecomposable direct summand of $Y_n$, then either $Y'_n\subseteq Y_{\tilde{n}}$ or $Y'_n\subseteq X_{\tilde{n}}$.	
	
	To justify the statement in the paragraph above, consider the following $\Ulie\left(\gfrak'_{\tilde{n}}\right)$-module 
	\begin{align}
		\tilde{X}_{\tilde{n}}:=\Ulie\left(\gfrak'_{\tilde{n}}\right)\cdot X'_n\,.
	\end{align}
	Let $\Xi\subseteq\hfrak^*$ denote the set of weights of $X'_n$.  By (\ref{eq:equaldim}), $\Xi$ is also an indecomposable weight set of $\tilde{X}_{\tilde{n}}$.  Therefore, $\tilde{X}_{\tilde{n}}$ can be decomposed as a direct sum $\tilde{X}_{\tilde{n}}^1\oplus\tilde{X}_{\tilde{n}}^2\oplus\ldots\oplus\tilde{X}^l_{\tilde{n}}$, where each $\tilde{X}^i_{\tilde{n}}$ is an indecomposable $\Ulie\left(\gfrak'_{\tilde{n}}\right)$-module, but as $\Xi$ is an indecomposable weight set of $\tilde{X}_{\tilde{n}}$, we must have $\Xi\subseteq \text{supp}\left(\tilde{X}^i_{\tilde{n}}\right)$ for some $i$.  However, this means $X'_n\subseteq \tilde{X}^i_{\tilde{n}}$.  Now, being indecomposable, $\tilde{X}^i_{\tilde{n}}$ must lie entirely either in $X_{\tilde{n}}$ or in $Y_{\tilde{n}}$.  Ergo, $X'_n$ is a subspace of $X_{\tilde{n}}$ or $Y_{\tilde{n}}$ for every $\tilde{n}\geq n$.

	The paragraph above proves that $X\cap X_n$ is given by a direct sum of some indecomposable direct summands of $X_n$.  Indeed, for a fixed direct summand $X'_n$ of $X_n$, we have only two possible scenarios: either $X'_n$ lies in $X_{\tilde{n}}$ for all sufficiently large $\tilde{n}\geq n$, or $X'_n$ lies in $Y_{\tilde{n}}$ for infinitely many $\tilde{n}$.  In the former case, $X'_n\subseteq X$, whereas, in the latter case, $X'_n\cap X=0$ and $X'_n\subseteq Y$.  In other words, $A_n:=X\cap X_n$ is a direct summand of $M_n$.  Write $B_n:=Y_n\oplus Z_n$, where $Z_n$ is the direct sum of indecomposable direct summands $X'_n$ of $X_n$ which intersect $X$ trivially. 
	
	Next, we fix $\xi\in\text{supp}(M)$.  We shall verify that $M^\xi=X^\xi+Y^\xi$.  For a given $v\in M^\xi$, $v=a_n+b_n$ for some $a_n\in \left(A_n\right)^{\xi}$ and $b_n\in \left(B_n\right)^{\xi}$.  Suppose that $n_0$ is a positive integer such that $\left(M_n\right)^{\xi}=M^\xi$ for all $n\geq n_0$.    We claim that there exists a positive integer $n_1\geq n_0$ such that $a_{n_1}=a_{n_1+1}=a_{n_1+2}=\ldots$.   This claim follows from the observation that $A_n\subseteq A_{n+1}\subseteq A_{n+2}\subseteq \ldots$ for all $n\geq n_0$.  The finite-dimensionality assumption implies that $\left(A_{n_1}\right)^{\xi}=\left(A_{n_1+1}\right)^{\xi}=\left(A_{n_1+2}\right)^{\xi}=\ldots$ for some $n_1\geq n_0$.  Furthermore, we note that $B_n\supseteq B_{n+1}\supseteq B_{n+2}\supseteq\ldots$; consequently, the finite-dimensionality assumption yields $\left(B_{n_1}\right)^{\xi}=\left(B_{n_1+1}\right)^{\xi}=\left(B_{n_1+2}\right)^{\xi}=\ldots$,    The claim follows immediately. 
	
	We write $a$ for the common value $a_{n_1}=a_{n_1+1}=a_{n_1+2}=\ldots$.  Set $b:=v-a$.  We shall now justify that $b$ is an element of $Y$.  Recall that $B_n=Y_n\oplus Z_n$, where $Z_n$ is the direct sum of indecomposable direct summands of $X_n$ that intersect $X$ trivially.  For $n\geq n_1$, we can write 
	\begin{align}
		b=\sum_{i=1}^{k_n}\,y_n^i+\sum_{j=1}^{l_n}\,z_n^j\,,
	\end{align}
	where $y_n^i$ and $z_n^j$ are nonzero elements of indecomposable direct summands of $Y_n$ and $Z_n$.  We shall now prove that each $y_n^i$ and each $z_n^j$ belong in $Y$.  For $\tilde{n}\geq n$, note that each $y_n^i$ lies either in $X_{\tilde{n}}$ or in $Y_{\tilde{n}}$.  If the former scenario occurs for all sufficiently large $\tilde{n}\geq n$, then $y_n^i\in X$, but this immediately implies $y_n^i=0$, which is a contradiction.  Ergo, the latter scenario occurs for infinitely many values $\tilde{n}\geq n$, whence $y_n^i \in Y$.  The same argument applies to each $z_n^j$.  Thus, we conclude that $y_n^i,z_n^j\in Y$ for every $i=1,2,\ldots,k_n$ and $j=1,2,\ldots,l_n$ with $n\geq n_1$.  Thence, $b\in Y$.  This proves that $M^\xi=X^\xi+Y^\xi$, leading to $M=X+Y$.
	
	Now, we shall check that $X\cap Y=0$.  Let $y_1,y_2,\ldots,y_k $ be linearly independent elements of $Y$ such that $x:=y_1+y_2+\ldots+y_k$ is in $X$ and that, for each $i=1,2,\ldots,k$, there are infinitely many positive integers $n$ for which $y_j\in Y_n$.  We may assume that there exists $\xi\in\text{supp}(M)$ with $y_i\in M^\xi$ for every $i=1,2,\ldots,k$.   Additionally, there exists a positive integer $m$ such that $x\in M_m$ and that $x\in X_n$ for every integer $n\geq m$.  
	
	Let $\bar{n}(j) \geq m$ be a positive integer such that $y_j\in Y_{\bar{n}(j)}$.  We decompose $y_j$ as
	\begin{align}
		y_j = y_j^1+y_j^2+\ldots+y_j^{r_j}\,,
	\end{align}
	where each $y_j^i$ is nonzero and in an indecomposable direct summand of $Y_{\bar{n}(j)}$.
	Pick an arbitrary $n\geq\bar{n}(j)$.  We note that each $y_j^i$ must lie in $X_n$ or in $Y_n$.  However, as $x\in X_n$, we conclude that $y_j^i$ is in  $X_n$, whence $y^i_j\in X_n$ for every $n\geq \bar{n}(j) $.  As a result, $y_j=\sum_{i=1}^{r_j}\,y_j^i$ is in $X$, which means $y_j=0$, and a contradiction is reached.

	Finally, we have the following equalities of composition series factor multiplicities: $\big[X:\Llie(\lambda)\big]=\big[M:\Llie(\lambda)\big]$, $\big[X:\Llie(\mu)\big]=0$, $\big[Y:\Llie(\lambda)\big]=0$, and $\big[Y:\Llie(\mu)\big]=\big[M:\Llie(\mu)\big]$.  That is, $M=X\oplus Y$ with $X\neq 0$ and $Y\neq 0$ .  This contradicts the assumption that $M$ is indecomposable.
		\end{pf} 
	\end{prop}
	
	\begin{rmk}
		The finite-dimensionality of the weight spaces of objects in $\bbggO$ plays an essential role in the proof of the proposition above.  Without this assumption, the proposition may be false.  As an example, $\gllie_\infty$ is an $\sllie_\infty$-module with the adjoint representation.  With respect to the splitting maximal toral subalgebra $\hfrak$ of diagonal matrices of $\sllie_\infty$, the $\sllie_\infty$-module $\gllie_\infty$ has a weight space decomposition
		\begin{align}
			\gllie_\infty=\left(\gllie_\infty\right)^0\oplus\left(\bigoplus_{\substack{{i,j\in\amsbb{Z}_{>0}} \\ {i\neq j}}}\,\amsbb{K}\,\textbf{E}_{i,j}\right)\,,
		\end{align}
		where $\left(\gllie_\infty\right)^0=\hfrak\oplus \amsbb{K}\,\textbf{E}_{1,1}$ is not a finite-dimensional weight space.  As the short exact sequence $0\to \sllie_\infty \to \gllie_\infty \to \amsbb{K}\to 0$ of $\sllie_\infty$-modules does not split, the trivial module $\amsbb{K}$ is in the same block as $\sllie_\infty$, despite the fact that, at each finite level $n$, $\amsbb{K}$ and $\sllie_n$ do not belong in the same block of $\bggO^{\sllie_n}_{\bfrak_n}$ (or $\bbggO^{\sllie_n}_{\bfrak_n}$), for any Borel subalgebra $\bfrak_n$ of $\sllie_n$.   In the proof above, with $\lambda:=0$, it can be seen that $X_n\cong \amsbb{K}$ and $Y_n=\sllie_n$ for every $n$, whilst $X=0$ and $Y=\sllie_\infty$.  Hence, the equality $\gllie_\infty=X\oplus Y$ does not hold.
	\end{rmk}

	\begin{prop}
		A block of $\bbggO$ containing $\Llie(\lambda)$ contains $\bbggO\llbracket\lambda\rrbracket$ as a subcategory.
	
		\begin{pf}
			Using the indecomposability of the Verma modules, we conclude that the block containing $\Llie(\lambda)$ must have $\bar\bggO\llbracket\lambda\rrbracket$ as a subcategory.  In other words, let $\mu,\nu\in \llbracket\lambda\rrbracket$.  Let $n\in\amsbb{N}$ be sufficiently large that $\mu=w \cdot \lambda$ for some $w\in W_n\left[\lambda_n\right]$.   (Here, $\xi_n$ denotes $\xi|_{\hfrak_n}$ for all $\xi\in\hfrak^*$.)
			
			From the finite-dimensional theory (see~\cite{bggo}), $W_n\left[\lambda_n\right]\cdot\lambda$ has a unique maximal element $\upsilon$ (with respect to the order $\preceq$ given by $\bfrak$).  Then, the Verma module $\Mlie(\upsilon)$ has $\Mlie(\mu)$ and $\Mlie(\lambda)$ as submodules due to the BGG Theorem (Theorem~\ref{thm:stronglinkagegen}).  Therefore, we have nonzero homomorphisms $\Mlie(\mu)\to \Mlie(\upsilon)$ and $\Mlie(\lambda)\to \Mlie(\upsilon)$.   Thus, the indecomposable modules $\Mlie(\mu)$ and $\Mlie(\lambda)$ are in the same block.  Furthermore, with  nontrivial homomorphisms $\Mlie(\mu)\to \Llie(\mu)$ and $\Mlie(\lambda)\to\Llie(\lambda)$, we conclude that $\Llie(\mu)$ and $\Llie(\lambda)$ are in the same block.
			
			Now, suppose that $M\in\bbggO$ is indecomposable with $\Llie(\mu)$ as a composition factor, where $\mu\in \llbracket\lambda\rrbracket$.  By Theorem~\ref{thm:gencomp} and Proposition~\ref{prop:blocks}, we see that $M$ has a submodule $N$ such that $M/N\cong \Llie(\nu)$ for some $\nu \in \llbracket\lambda\rrbracket$.  Thus, the nonzero homomorphism $M\to M/N$ establishes that $M$ is in the same block as $\Llie(\nu)$, which is also in the same block as $\Llie(\lambda)$.  Thus, every indecomposable object $M$ whose composition factors are of the form $\Llie(\mu)$ with $\mu\in \llbracket\lambda\rrbracket$ is in the same block as $\Llie(\lambda)$.  The proposition follows immediately.
		\end{pf}
	\end{prop}

	\begin{thm}
		Let $\Omega$ denote the set of integral Weyl dot-orbits.  Then, the full abelian subcategories $\bar\bggO\llbracket\lambda\rrbracket$, where $\llbracket\lambda\rrbracket\in\Omega$, are the blocks of $\bar\bggO$, and
		\begin{align}
		\bar	\bggO =\bigoplus_{\llbracket\lambda\rrbracket\in \Omega}\,\bar\bggO\llbracket\lambda\rrbracket\,.
		\end{align}
		\label{thm:blocks}
		
		\begin{pf}
			From the proposition above, we know that each block of $\bbggO$ contains $\bbggO\llbracket\lambda\rrbracket$ for some $\lambda\in\hfrak^*$.  We shall prove that the block containing $\bbggO\llbracket\lambda\rrbracket$ must then coincide with $\bbggO\llbracket\lambda\rrbracket$.  If the block contains an indecomposable module $M$ not in $\bbggO\llbracket\lambda\rrbracket$, then there exists a finite sequence $M=M_0,M_1,M_2,\ldots,M_k$ of indecomposable modules in this block such that $M_k\in\bbggO\llbracket\lambda\rrbracket$, $M_{k-1}\notin\bbggO\llbracket\lambda\rrbracket$, and for each $i=0,1,2,\ldots,k-1$, either $\Hom_{\bbggO}\left(M_i,M_{i+1}\right)\neq 0$ or $\Hom_{\bbggO}\left(M_{i+1},M_i\right)\neq 0$.
			
			If $\Hom_{\bbggO}\left(M_{k-1},M_{k}\right)\neq 0$ or $\Hom_{\bbggO}\left(M_{k},M_{k-1}\right)\neq 0$, then $M_{k-1}$ has a composition factor $\Llie(\mu)$ (which is also a composition factor of $M_k$) for some $\mu\in \llbracket\lambda\rrbracket$, which then means that $M_{k-1}\in\bbggO\llbracket\lambda\rrbracket$ by Proposition~\ref{prop:blocks}.  This contradicts the assumption that $M_{k-1}\notin \bbggO\llbracket\lambda\rrbracket$.  Therefore, the blocks of $\bbggO$ are precisely $\bbggO\llbracket\lambda\rrbracket$.
			
			To complete the proof, let now $M$ be an arbitrary object in $\bbggO$.  By Theorem~\ref{thm:directsumdecomposition}, $M$ has a direct sum decomposition with indecomposable summands.  Write $M\llbracket\lambda\rrbracket$ for the (direct) sum of the direct summands of $M$ that belong to $\bbggO\llbracket\lambda\rrbracket$.  Then, we can clearly see that
		\begin{align}
			M=\bigoplus_{\llbracket\lambda\rrbracket\in \Omega}\,M\llbracket\lambda\rrbracket\,.
		\end{align}
		Note that this direct sum may be an uncountable direct sum.
		\end{pf}	
	\end{thm}

\subsection{Formal Characters}

\begin{define}
	For each $M\in\bar\bggO$, the \emph{formal character} of $M$ is the element $\ch(M)$ of the set $\mathcal{S}$ of functions $f:\hfrak^*\to\amsbb{Z}$ defined by 
	\begin{align}
		\ch(M)\defeq \sum_{\lambda\in\hfrak^*}\,\dim_\amsbb{K}\left(M^\lambda\right)\,e^\lambda\,,
	\end{align}
	where $e^\lambda:\hfrak^*\to\amsbb{Z}$ is defined by $e^\lambda(\mu):=\deltaup_{\lambda,\mu}$ and $\deltaup$ is the Kronecker delta.  
\end{define}

\begin{prop}
	Let $M,N\in\bar\bggO$.  Then, $\ch(M\oplus N)=\ch(M)+\ch(N)$.  If $M\otimes N\in\bggO$, then $\ch(M\otimes N)=\ch(M)\cdot\ch(N)$, where the coefficients of $\ch(M\otimes N)$ are given by \emph{convolutions}, namely, for  $f,g\in\mathcal{S}$,
	\begin{align}
		(f\cdot g)(\lambda)=\sum_{\mu\in\hfrak^*}\,f(\mu)\,g(\lambda-\mu)\,,
	\end{align}
	when all such summations are finite.  (Note that $e^\lambda \cdot e^\mu=e^{\lambda+\mu}$ for all $\lambda,\mu\in\hfrak^*$.)
\end{prop}

\begin{prop}
	Let the \emph{Kostant function} $p\in\mathcal{S}$ be defined via $p=\sum_{\mu\succeq 0}\,\varp(\mu)\,e^{-\mu}$, where the \emph{partition function} $\varp$ is defined as follows: $\varp(\mu)$ denotes the number of ways to express $\mu\succeq0$ as a sum of positive roots.  Then, $\ch\big(\Mlie(\lambda)\big)=p\cdot e^\lambda$ for every $\lambda\in\hfrak^*$.
\end{prop}

\begin{thm}
	For every $M\in\bar\bggO$, $\ch(M)=\sum_{\lambda\in\hfrak^*}\,\big[M:\Llie(\lambda)\big]\,\ch\big(\Llie(\lambda)\big)$.
	\label{thm:character}
\end{thm}

\begin{cor}
	If $0\to M_0\to M_1\to M_2\to\ldots \to M_k \to 0$ is an exact sequence of objects in $\bar\bggO$, then $\sum_{r=0}^k\,(-1)^r\,\ch\left(M_r\right)=0$.
\end{cor}

Let $M$ be an indecomposable object in $\bar\bggO$.  Suppose that $u_1,u_2,\ldots\in M$ are weight vectors that generate $M$.  Write $M_n:=\sum_{i=1}^n\,\Ulie\left(\gfrak_n'\right)\cdot u_i$, where $\gfrak'_n:=\gfrak_n+\hfrak$.  For any simple object $L:=\Llie(\lambda) \in \bar\bggO$, write $L_n:=\Llie_n\left(\lambda\right)=\Llie\left(\lambda;\gfrak'_n,\bfrak'_n,\hfrak\right)$, where $\bfrak_n':=\bfrak_n+\hfrak$.

\begin{prop}
	For all sufficiently large $n$, 
	\begin{align}
		[M:L]=\big[M_n:L_n\big]\,.
	\end{align}
	\label{prop:MnLn}
	
	\begin{pf}
	
	To prove this, if the index $n$ is so large that $\dim_\amsbb{K}\left(M_n^{\lambda}\right)=\dim_\amsbb{K}\left(M^\lambda\right)$, then it is immediate that $[M:L]\leq \left[M_n:L_n\right]$.  Furthermore, the sequence $\big(\left[M_n:L_n\right]\big)_{n\in\amsbb{Z}_{> 0}}$ is nonincreasing at some point.  Therefore, for some $m\geq [M:L]$, we have $\left[M_n:L_n\right]=m$ for every large $n$, say $n\geq n_0$.

	Now, we look at the character $\ch(M)$ and $\ch\left(M_n\right)$.  For sufficiently large $n\geq n_0$, we must have $\dim_\amsbb{K}\left(M^\mu\right)=\dim_\amsbb{K}\left( M_n^{\mu}\right)$ for every $\mu \succeq \lambda$.  Theorem~\ref{thm:character} guarantees that $m=[M:L]$.
	\end{pf}
\end{prop}

\subsection{Kazhdan-Lusztig Multiplicities}

Let $\gfrak_n':=\gfrak_n+\hfrak$ and $\bfrak'_n:=\bfrak_n+\hfrak$.  Note that the Weyl group of $\gfrak_n'$ is still the Weyl group $W_n$ of $\gfrak_n$.  For each $\xi\in\hfrak^*$, write $\Mlie_n(\xi)$ and $\Llie_n(\xi)$ for the Verma module $\Mlie\left(\xi;\gfrak'_n,\bfrak_n',\hfrak\right)$ and the simple module $\Llie\left(\xi;\gfrak_n',\bfrak_n',\hfrak\right)$, respectively.

Fix a regular integral weight $\lambda$.  Take $\mu\in W\cdot \lambda$.  For each $n\in\amsbb{Z}_{> 0}$, write $\nu_n$ for the antidominant weight in $\hfrak^*$ that is strongly linked to $\lambda$ with respect to $\bfrak_n'$.  In addition, there exist elements $x_n$ and $y_n$ of $W_n$ such that $x_n^{-1}\cdot\lambda=  \nu_n$ and $y_n^{-1}\cdot \mu=\nu_n$.

From the finite-dimensional theory and Remark~\ref{rmk:reductive}, we have
\begin{align}
	\big[\Mlie_n\left(\lambda\right):\Llie_n\left(\mu\right)\big]=P^{W_n}_{w_n^0x_n,w_n^0y_n}(1)\,,
\end{align}
where $w_n^0\in W_n$ is the longest element of $W_n$.  Combining this result with Proposition~\ref{prop:MnLn}, we obtain the proposition below.

\begin{prop}
	There exists a positive integer $n_0$ such that, for all $n\geq n_0$,
	\begin{align}
		\big[\Mlie\left(\lambda\right):\Llie\left(\mu\right)\big]=P^{W_n}_{w_n^0x_n,w_n^0y_n}(1)\,.
	\end{align}
\end{prop}

 Fix $x,y\in W$ and set $m(x,y)$ to be the smallest positive integer $m$ such that $x,y\in W_m$.   Theorem~\ref{thm:KLpolystab} dictates that 
 \begin{align}
 	P_{x,y}^{W_{m(x,y)}}(q)=P_{x,y}^{W_n}(q)=P_{x,y}^W(q)\,.
 \end{align}
 This result gives rise to the following proposition.
 
 \begin{prop}
 	For every $x\in W$ and for each regular antidominant weight $\lambda$, 
 	\begin{align}
 		\big[\Mlie(x\cdot \lambda)\big]=\sum_{y\preceq x}\,P^{W_{m(x,y)}}_{w^0_{m(x,y)}x,w^0_{m(x,y)}y}(1)\,\big[\Llie(y\cdot \lambda)\big]\,,
 	\end{align}
 	or equivalently
 	\begin{align}
 		\big[\Llie(x\cdot \lambda)\big]=\sum_{y\preceq x}\,(-1)^{\ell(x)-\ell(y)}\,P^{W_{m(x,y)}}_{x,y}(1)\,\big[\Mlie(y\cdot\lambda)\big]\,.
 	\end{align}
 	(Note that the two equations above are equalities in the Grothendieck group of $\bar\bggO$.)
 \end{prop}
 
 \subsection{$\Hom$ and $\Ext^\bullet$ Functors}
 
 	Unless otherwise specified, $\Ext$ denotes $\Ext_{\bar\bggO}$.  Similarly, $\Hom$ denotes $\Hom_{\bar\bggO}$.
 	
 	\begin{prop}
 		Let $\lambda,\mu\in\hfrak^*$.
 		\label{prop:extintro}
 		
 		\begin{enumerate}[(a)]
 			\item If $M$ is a $\Ulie(\gfrak)$-module such that, for all weights $\upsilon\in \text{supp}(M)$, we have $\lambda \not\prec \upsilon$, then $\Ext^1\big(\Mlie(\lambda),M\big)=0$.  In particular,
 			\begin{align}
 				\Ext^1\big(\Mlie(\lambda),\Llie(\lambda)\big)=0\text{ and }\Ext^1\big(\Mlie(\lambda),\Mlie(\lambda)\big)=0\,.
 			\end{align}
 			\item If $\mu\preceq \lambda$, then $\Ext^1\big(\Mlie(\lambda),\Llie(\mu)\big)=0$.
 			\item If $\mu\prec\lambda$ and $\Nlie(\lambda)$ is the maximal proper submodule of $\Mlie(\lambda)$, then
 			\begin{align}
 				\Ext^1\big(\Llie(\lambda),\Llie(\mu)\big)\cong \Hom\big(\Nlie(\lambda),\Llie(\mu)\big)\,.
 				\label{eq:extandN}
 			\end{align}
 			\item $\Ext^1\big(\Llie(\lambda),\Llie(\lambda)\big)=0$.
 		\end{enumerate}
 		
 		\begin{pf}$\phantom{a}$
 		
 			\begin{enumerate}[(a)]
 				\item
 				Given an extension $0 \to M \overset{i}\longrightarrow E \overset{p}{\longrightarrow} \Mlie(\lambda) \to 0$ in $\bar\bggO$, let $e\in E$ be such that $p(e)$ be a highest-weight vector of $\Mlie(\lambda)$.  Due to the hypothesis, the submodule $V$ of $E$ generated by $e$ is a highest-weight module with highest weight $\lambda$.  Since $V$ is mapped surjectively by $p$ onto $\Mlie(\lambda)$, we conclude that $p$ induces an isomorphism $V \cong \Mlie(\lambda)$, whence the exact sequence splits.
 				\item This is an immediate consequence of (a).
 				\item Starting with the short exact sequence $0\to \Nlie(\lambda) \to \Mlie(\lambda) \to \Llie(\lambda)\to 0$, we get the following long exact sequence of $\Ext$-groups:
 				\begin{align}
 					\ldots \to \Hom\big(\Mlie(\lambda),\Llie(\mu)\big) &\to \Hom\big(\Nlie(\lambda),\Llie(\mu)\big)
 					\nonumber\\
 					&\to \Ext^1\big(\Llie(\lambda),\Llie(\mu)\big) \to \Ext^1\big(\Mlie(\lambda),\Llie(\mu)\big)\to\ldots\,.
 				\end{align}
 				By (b), $\Ext^1\big(\Mlie(\lambda),\Llie(\mu)\big)=0$.  Furthermore, it is clear that $\Hom\big(\Mlie(\lambda),\Llie(\mu)\big)=0$.  Therefore, we have the isomorphism (\ref{eq:extandN}).
 				\item Replace $\mu$ by $\lambda$ in the proof of (c).  We note that $\Hom\big(\Nlie(\lambda),\Llie(\lambda)\big)=0$.  By (b), $\Ext^1\big(\Mlie(\lambda),\Llie(\mu)\big)=0$.  Thus, $\Ext^1\big(\Llie(\lambda),\Llie(\lambda)\big)=0$ as well.
 			\end{enumerate}
 		\end{pf}
 	\end{prop}
 	
 	\begin{prop}
 		Let $\lambda,\mu\in\hfrak^*$.
 		\begin{enumerate}[(a)]
 			\item For every $M,N\in\bar\bggO$ and $k\in\amsbb{Z}_{\geq 0}$, we have $\Ext^k(M,N)\cong \Ext^k\big(N^\vee,M^\vee\big)$.
 			\item The image of any homomorphism $\Mlie(\mu)\to\Vlie(\lambda)$ is a submodule of $\Llie(\lambda)\subseteq \Vlie(\lambda)$.  This means 
			\begin{align}
				\dim_\amsbb{K}\Big(\Hom\big(\Mlie(\mu),\Vlie(\lambda)\big)\Big)=\left\{
					\begin{array}{ll}
						1\,,	&\text{if }\mu=\lambda\,,
						\\
						0\,,	&\text{if }\mu\neq \lambda\,.
					\end{array}
				\right.
			\end{align}
 			\item $\Ext^1\big(\Mlie(\mu),\Vlie(\lambda)\big)=0$ for all $\lambda$ and $\mu$.
 		\end{enumerate}
 		\label{prop:extverma}
 		
 		\begin{pf}$\phantom{a}$
 			
 			\begin{enumerate}[(a)]
 				\item This part is trivial due to the fact that duality is an antiequivalence of the category $\bbggO$ with itself.
 				\item Let $M$ be the image of a nonzero homomorphism $\Mlie(\mu)\to\Vlie(\lambda)$.  Then, $M$ is a highest-weight submodule of $\Vlie(\lambda)$ with highest weight $\mu$.  Since $\Llie(\lambda)$ is contained in every nonzero submodule of $\Vlie(\lambda)$, we see that $\Llie(\lambda)\subseteq M$, so $\mu\preceq \lambda$.  However, the composition factors of $\Vlie(\lambda)$ are the same as those of $\Mlie(\lambda)$, which are simple modules with highest weight less than or equal to $\lambda$.  This means $\mu\succeq \lambda$.  Consequently, $\mu=\lambda$ must hold, whence $M=\Llie(\lambda)$.
 				\item   If $\lambda\not\prec \mu$, then $M:=\Vlie(\mu)$ satisfies the hypothesis of Proposition~\ref{prop:extintro}(a).  Therefore, 
 				\begin{align}\Ext^1\big(\Mlie(\lambda),\Vlie(\mu)\big)=0\,.\end{align}  By (a), we have 
 					\begin{align}
 						\Ext^1\big(\Mlie(\mu),\Vlie(\lambda)\big)\cong\Ext^1\big(\Mlie(\lambda),\Vlie(\mu)\big)=0\,.
					\end{align}
					If $\lambda \preceq \mu$, then $M:=\Vlie(\lambda)$ satisfies the hypothesis of Proposition~\ref{prop:extintro}(a), with $\mu$ replacing $\lambda$ in that proposition.  The same conclusion follows.
 			\end{enumerate}
 		\end{pf}
 	\end{prop}
	
	\begin{rmk}
		As a comparison with Part (c) of the proposition above, we note that $\Ext^1\big(\Vlie(\lambda),\Mlie(\mu)\big)$ can be nontrivial.  For example, consider $\gfrak:=\sllie_\infty$ with $\hfrak:=\hfrak_{\text{A}}$ and $\bfrak:=\bfrak_{\textbf{1A}}$.  If 
		\begin{align}
			\lambda:=\left(0,0,-\frac{1}{2},-\frac{2}{3},-\frac{3}{4},\ldots\right)
		\end{align} and 
		\begin{align}
			\mu:=\left(0,2,\frac{1}{2},\frac{1}{3},\frac{1}{4},\ldots\right)\,,
		\end{align} 
		then we have the equality 
		\begin{align}
			\mu=s_{\epsilon_1-\epsilon_2}\cdot\lambda={(1\;\;2)}\cdot\lambda\,.
		\end{align}  We shall see later that Theorem~\ref{thm:almostdominant} produces the injective hull $\Ilie(\mu) \in \bbggO$.   We have a nonsplitting short exact sequence
		\begin{align}
			0\to \Vlie(\mu) \to \Ilie(\mu) \to \Vlie(\lambda)\to 0\,.
		\end{align}
		As $\mu$ is an almost antidominant weight, $\Mlie(\mu)$ is simple.  Therefore, 
		\begin{align}
			\Mlie(\mu)\cong \Llie(\mu)\cong \big(\Llie(\mu)\big)^\vee\,,\end{align} so that 
		\begin{align}\Vlie(\mu)=\big(\Mlie(\mu)\big)^\vee\cong \Mlie(\mu)\,.\end{align} 
		This shows that $\Ext^1\big(\Vlie(\lambda),\Mlie(\mu)\big)\neq 0$.
		
		Note that the module $\Ilie(\mu)$ given above is an example of a \emph{tilting module}, that is $\Ilie(\mu)$ has both a filtration by submodules with successive co-Verma factors, and a filtration by submodules with successive Verma factors.  To elaborate, we have the following filtration with successive co-Verma factors from the short exact sequence above:
		\begin{align}
			0\subsetneq \Vlie(\mu)\subsetneq \Ilie(\mu)\,.
		\end{align}
		Now, we take $v$ to be a singular vector of $\Ilie(\mu)$ with weight $\lambda$, then the submodule $\Ulie(\gfrak)\cdot v$ is a Verma submodule of $\Ilie(\mu)$ isomorphic to $\Mlie(\lambda)$.  It can be shown that $\Ilie(\mu)/\big(\Ulie(\gfrak)\cdot v\big)\cong \Mlie(\mu)$.  Therefore, $\Ilie(\mu)$ has a Verma filtration
		\begin{align}
			0\subsetneq \Mlie(\lambda)\subsetneq \Ilie(\mu)\,.
		\end{align}
		In fact, $\Ilie(\mu)$ is \emph{self-dual}, i.e., $\Ilie(\mu)=\big(\Ilie(\mu)\big)^\vee$.  That is, $\Ilie(\mu)$ equals the projective cover $\Plie(\mu)$ of $\Llie(\mu)$.  
	\end{rmk}
	
	It is natural to make the following conjecture.  The finite-dimensional version is true, and it is equivalent to the Kazhdan-Lusztig Conjecture (Theorem~\ref{thm:KLconj}).  See also~\cite{bggo}.
	
	\begin{conj}
		Let $\lambda\in\hfrak^*$ be a regular antidominant integral weight.  For $x,y\in W$, we write $m(x,y)$ for the smallest positive integer $m$ such that $x,y\in W_m$.  Then,
		\begin{align}
			P^W_{x,y}(q)=P^{W_{m(x,y)}}_{x,y}(q) = \sum_{i=0}^{\left\lfloor \frac{\ell(x,y)-1}{2}\right\rfloor}\,q^i \,\dim_\amsbb{K}\Big(\Ext^{\ell(x,y)-2i}\big(\Mlie(x\cdot \lambda),\Llie(y\cdot\lambda)\big)\Big)\,.
		\end{align}
	\end{conj}	
 
 \pagebreak

\section{Truncated Category $\bggO$}
\label{ch:trunc}

	In this section, $\gfrak$ is a root-reductive Lie algebra with a Dynkin Borel subalgebra $\bfrak$.

\subsection{Truncation}
	
	As before, $M^\vee$ and $f^\vee$ denote the duals  in $\bar\bggO$ of an object $M$ and a homomorphism $f$, respectively.   We shall now define a truncation method of the category $\bbggO$ using an idea from~\cite{RCW}.

\begin{define}
	 For $\lambda\in\hfrak^*$, we write $\bar\bggO^{\preceq \lambda}$ for the full subcategory of $\bar\bggO$ consisting of all modules $M\in\bar\bggO$ whose weights are less than or equal to $\lambda$ with respect to the partial order $\preceq$ on $\hfrak^*$.
\end{define}

\begin{prop}
	For each $\lambda\in\hfrak^*$, let $\tcyr_{\preceq \lambda}:\bar\bggO\to\bar\bggO^{\preceq\lambda}$ be defined as 
	\begin{align}
		\tcyr_{\preceq\lambda}M\defeq\sum_{\substack{{N\subseteq M}\\{N\in\bar\bggO^{\preceq\lambda}}}}\,N
	\end{align}
	and
	\begin{align}
		\tcyr_{\preceq\lambda}f \defeq f|_{\tcyr_{\preceq\lambda}M}
	\end{align}
	for all $M\in\bar\bggO$ and for all homomorphisms $f:M\to L$ of objects in $\bar\bggO$.
	Then, $\tcyr_{\preceq\lambda}$ is a left-exact (covariant) functor.  We shall call $\tcyr_{\preceq \lambda}$ the \emph{truncation functor} (with the upper bound $\lambda$).
\end{prop}

\begin{cor}
	For each $\lambda\in\hfrak^*$, let $\tcyr^\vee_{\preceq \lambda}:\bar\bggO\to\bar\bggO^{\preceq\lambda}$ be defined as 
	\begin{align}
		\tcyr^\vee_{\preceq\lambda}M\defeq \Big(\tcyr_{\preceq\lambda}\left(M^\vee\right)\Big)^\vee
	\end{align}
	and
	\begin{align}
		\tcyr^\vee_{\preceq\lambda}f \defeq \Big(\tcyr_{\preceq\lambda}\left(f^\vee\right)\Big)^\vee
	\end{align}
	for all $M\in\bar\bggO$ and for all homomorphisms $f:M\to L$ of objects in $\bar\bggO$.
	Then, $\tcyr^\vee_{\preceq\lambda}$ is a right-exact (covariant) functor.  We shall call $\tcyr^\vee_{\preceq \lambda}$ the \emph{dual truncation functor} (with the upper bound $\lambda$).
	\label{cor:modtrunc}
\end{cor}

\begin{prop}
	Let $\lambda\in\hfrak^*$.  If $I$ is an injective object in $\bar\bggO$, then $\tcyr_{\preceq\lambda}I$ is injective in $\bar\bggO^{\preceq \lambda}$.  If $P$ is a projective object in $\bar\bggO$, then $\tcyr^\vee_{\preceq\lambda}P$ is projective in $\bar\bggO^{\preceq \lambda}$. 
	
	\begin{pf}
		Let $I$ be an injective object in $\bar\bggO$ and $0\to \tcyr_{\preceq\lambda}I \to M\to N\to 0$ an exact sequence of objects in $\bar\bggO^{\preceq\lambda}$.  We have the injection $\tcyr_{\preceq\lambda}I \overset{\subseteq}{\longrightarrow} I$.  Because $I$ is an injective object in $\bar\bggO$ and $0\to \tcyr_{\preceq\lambda}I \to M\to N\to 0$ is also an exact sequence of objects in $\bar\bggO$, we conclude that there exists a homomorphism $\phi:M\to I$ such that the diagram below is commutative:
		
		\begin{equation}
		\begin{tikzcd}
		0\arrow{r}& \tcyr_{\preceq \lambda}I \arrow{r} \arrow{d}{\rotatebox{270}{$\subseteq$} }& M \arrow{r}\arrow{ld}{\phi} &N \arrow{r} & 0\\
		&I\,.&&&
		\end{tikzcd}
		\end{equation}

		Since the image of $M$ under $\phi$ must be in $\bar\bggO^{\preceq\lambda}$, we see that $\im(\phi)\subseteq \tcyr_{\preceq\lambda}I$.  Thence, we indeed have a commutative diagram
		
		\begin{equation}
		\begin{tikzcd}
		0\arrow{r}& \tcyr_{\preceq \lambda}I \arrow{r} \arrow{d}{\rotatebox{270}{$=$} }& M \arrow{r}\arrow{ld}{\phi} &N \arrow{r} & 0\\
		&\tcyr_{\preceq \lambda}I \,.&&&
		\end{tikzcd}
		\end{equation}
		Hence, the exact sequence  $0\to \tcyr_{\preceq\lambda}I \to M\to N\to 0$ splits.  Thus, $\tcyr_{\preceq\lambda}I$ is injective.  
		
		For the second part of the proposition, we employ the duality from Corollary~\ref{cor:modtrunc}.  The proof is now complete.
	\end{pf}
\end{prop}

\begin{rmk}
	Note that $\bbggO^{\preceq \lambda}$ is a full subcategory of $\bbggO$ which is closed under extensions, taking quotient objects, and taking direct sums.  That is, $\bbggO^{\preceq\lambda}$ is a \emph{torsion subcategory} of $\bbggO$ (see~\cite{torsionsubcat}).  If $\pazocal{F}^{\preceq \lambda}$ denote the class of objects $M$ in $\bbggO$ such that $\tcyr_{\preceq}M=0$, then $\left(\bbggO^{\preceq \lambda},\pazocal{F}^{\preceq\lambda}\right)$ forms a \emph{torsion pair}.  In other words, we have the following observations:
	\begin{enumerate}[(a)]
		\item $\bbggO^{\preceq \lambda}\cap\pazocal{F}^{\preceq\lambda}=\{0\}$;
		\item If $T\to A\to 0$ is an exact sequence in $\bbggO$ with $T\in \bbggO^{\preceq\lambda}$, then $A\in\pazocal{F}^{\preceq\lambda}$;
		\item If $0\to A\to F$ is an exact sequence in $\bbggO$ with $F\in\pazocal{F}^{\preceq\lambda}$, then $A\in\pazocal{F}^{\preceq\lambda}$;
		\item For every $M\in\bbggO$, there exists an exact sequence $0\to T\to M\to F\to 0$ with $T\in\bbggO^{\preceq\lambda}$ and $F\in\pazocal{F}^{\preceq\lambda}$;
		\item For all $T\in\bbggO^{\preceq\lambda}$ and $F\in\pazocal{F}^{\preceq}$, we have $\Hom_{\bbggO}\left(T,F\right)=0$.
	\end{enumerate}
\end{rmk}

\subsection{Injective Objects}
\label{sec:injectiveobjects}

As in the proof of Proposition~\ref{prop:blocks}, $\gfrak'_n$ and $\bfrak'_n$ denote $\hfrak+\gfrak_n$ and $\hfrak+\bfrak_n$, respectively.

\begin{prop}
		Fix $n\in\amsbb{Z}_{> 0}$.  Suppose that $I_{n+1}$ is an injective object in $\bggO^{\gfrak'_{n+1}}_{\bfrak'_{n+1}}$.  Then, the restriction $I_n:=\Res^{\gfrak'_{n+1}}_{\gfrak'_n}\left(I_{n+1}\right)$ is an injective object in $\bggO^{\gfrak'_n}_{\bfrak'_n}$.		
		\label{prop:injrestricted}
		
		\begin{pf}
			Let $0\to M_n\overset{\varphi_n}{\longrightarrow} N_n$ be an exact sequence in $\bggO^{\gfrak'_n}_{\bfrak'_n}$ along with a homomorphism $f_n:M_n\to I_{n}$ of $\gfrak'_n$-modules.  Now, let $\pfrak'_{n+1}$ denote the parabolic subalgebra
			\begin{align}
				\pfrak'_{n+1}=\gfrak'_n+\bfrak'_{n+1}\,.
			\end{align}
			Take $\left\{x_{\pm\alpha}\,\boldsymbol{|}\,\alpha\in\Delta^+\right\}\cup\left\{h_\beta\,\boldsymbol{|}\,\beta\in\Sigma^+\right\}$ for a Chevalley basis of $\gfrak$.  We equip each object $L$ in $\bggO^{\gfrak'_n}_{\bfrak'_n}$ with a $\pfrak'_{n+1}$-module structure by requiring that,  for each $\bfrak'_{n+1}$-positive root $\alpha$ of $\gfrak'_{n+1}$ which is not a root of $\gfrak'_n$, 
			\begin{align}
				x_\alpha\cdot L=0\,.
			\end{align}

			Note that $\Ulie\left(\gfrak'_{n+1}\right)$ is a free (whence flat) $\Ulie\left(\pfrak'_{n+1}\right)$-module due to the PBW Theorem.  Hence, the parabolic induction functor $\Ulie\left(\gfrak'_{n+1}\right)\underset{\Ulie\left(\pfrak'_{n+1}\right)}{\otimes}\boldsymbol{\_}$ is exact, that is we have an exact sequence of $\gfrak'_{n+1}$-modules
			\begin{align}
				0\to \Ulie\left(\gfrak'_{n+1}\right)\underset{\Ulie\left(\pfrak'_{n+1}\right)}{\otimes}M_n\overset{\varphi_{n+1}}{\longrightarrow} \Ulie\left(\gfrak'_{n+1}\right)\underset{\Ulie\left(\pfrak'_{n+1}\right)}{\otimes}N_n\,,
			\end{align}
			where the $\Ulie\left(\gfrak_{n+1}'\right)$-module homomorphism $\varphi_{n+1}$ is given by $\varphi_{n+1}:=\text{id}_{\Ulie\left(\gfrak_{n+1}'\right)}\otimes \varphi_n$.  Then, we define the homomorphism of $\gfrak'_{n+1}$-modules $f_{n+1}:\Ulie\left(\gfrak_{n+1}'\right)\underset{\Ulie\left(\pfrak'_{n+1}\right)}{\otimes}M_n\to I_{n+1}$ by setting
			\begin{align}
				f_{n+1}(x\otimes v):=x\cdot f_n(v)
			\end{align}
			for all $x\in \Ulie\left(\gfrak'_{n+1}\right)$ and $v\in M_n$.  Since $\Ulie\left(\gfrak'_{n+1}\right)\underset{\Ulie\left(\pfrak'_{n+1}\right)}{\otimes}M_n$ and $\Ulie\left(\gfrak'_{n+1}\right)\underset{\Ulie\left(\pfrak'_{n+1}\right)}{\otimes}N_n$ are objects in $\bggO^{\gfrak'_{n+1}}_{\bfrak'_{n+1}}$, we conclude by injectivity of $I_{n+1}$ in $\bggO^{\gfrak'_{n+1}}_{\bfrak'_{n+1}}$ that there exists a homomorphism of $\gfrak'_{n+1}$-modules $\psi_{n+1}:\Ulie\left(\gfrak'_{n+1}\right)\underset{\Ulie\left(\pfrak'_{n+1}\right)}{\otimes}N_n\to I_{n+1}$ such that 
			\begin{align}
				f_{n+1}=\psi_{n+1}\circ\varphi_{n+1}\,.
			\end{align}
			We then define $\psi_n:N_n\to I_{n+1}$ by setting 
			\begin{align}
				\psi_n(u):=\psi_{n+1}\left(1_{\Ulie\left(\gfrak'_{n+1}\right)}\otimes u\right)
			\end{align}
			for every $u\in N_n$.  It is easy to see that $f_n=\psi_n\circ \varphi_n$ and that $I_n=\Res^{\gfrak'_{n+1}}_{\gfrak'_n}\left(I_{n+1}\right)$ is an object in $\bggO^{\gfrak'_n}_{\bfrak_n'}$, whence $I_{n}$ is injective in $\bggO^{\gfrak'_n}_{\bfrak'_n}$.
		\end{pf}
	\end{prop}

\begin{define}
Let $\lambda\in\hfrak^*$.  Then, we say that $\lambda$ is \emph{dominant} if $(\lambda+\rho)\left(h_\alpha\right)\notin\amsbb{Z}_{<0}$ for all positive roots $\alpha$ of $\gfrak$.  We say that $\lambda$ is \emph{almost dominant} if $(\lambda+\rho)\left(h_\alpha\right)\in\amsbb{Z}_{<0}$ for at most finitely many positive roots $\alpha$ of $\gfrak$.
\label{def:dominantweights}
\end{define}

\begin{thm}
	Let $I_n\in\bggO^{\gfrak'_n}_{\bfrak'_n}$ for each $n\in\amsbb{Z}_{> 0}$ be an injective object.  Suppose that we have an embedding $I_n\overset{\psi_n}{\longrightarrow}I_{n+1}$ for every $n$.  If the direct limit $I=\lim_{\underset{n}{\boldsymbol{\longrightarrow}}}\,I_n$ is an object of $\bar\bggO$, then $I$ is injective in $\bar\bggO$.
	\label{thm:injdirlim}
	
	\begin{pf}
		Let an injective homomorphism $M\overset{\varphi}{\longrightarrow}N$ and a homomorphism $f:M\to I$ be given.  Without loss of generality, we can assume that all weights of $M$, $N$, and $I$ lie within $\lambda+\Lambda$ for some $\lambda\in\hfrak^*$.  In particular, we can assume that the modules $M$, $N$, and $I$ are generated by countably many weight vectors.
		
		We suppose that $N$ is generated by weight vectors $u_1,u_2,u_3,\ldots$, with the corresponding weights $\mu_1,\mu_2,\mu_3,\ldots$.  Let $N_n$ denote the $\gfrak_n$-submodule
		\begin{align}	
			N_n:=\sum_{i=1}^n\,\Ulie\left(\gfrak'_n\right)\cdot u_i\,.
		\end{align}
		Then, we define $M_n$ as $\varphi^{-1}\left(N_n\right)$.  Note that both $M_n$ and $N_n$ are objects in $\bggO^{\gfrak'_n}_{\bfrak'_n}$.  By Proposition~\ref{prop:injrestricted}, for every $m\geq n$, $\text{Res}^{\gfrak'_m}_{\gfrak'_n}\left(I_m\right)$ is an injective module in $\bggO^{\gfrak'_n}_{\bfrak'_n}$.  Therefore, as $I$ has finite-dimensional weight spaces, we can assume without loss of generality that each $I_n$ satisfies the property that 
		\begin{align}
			\dim_\amsbb{K}\big(\left(I_n\right)^{\mu_i}\big)=\dim_\amsbb{K}\left(I^{\mu_i}\right)
		\end{align}
		for every $i=1,2,\ldots,n$.
		
		Let $\varphi_n:=\varphi|_{M_n}$ and $f_n:=f|_{M_n}$.  From the definitions above, we have the diagram of objects of $\bggO^{\gfrak'_n}_{\bfrak'_n}$
		\begin{equation}
		\begin{tikzcd}
		0\arrow{r}& M_n \arrow{r}{\varphi_n} \arrow{d}{f_n} & N_n\\
		&I_n \,.&
		\end{tikzcd}
		\end{equation}
		Because $I_n$ is injective, there exists a $\gfrak'_n$-module homomorphism $F_n:N_n\to I_n$ such that $F_n\circ \varphi_n=f_n$.
		
		Write $\mathcal{F}^1_n$ for the set of all maps $F_n:N_n\to I_n$ such that 
		\begin{align}
			F_n\circ \varphi_n=f_n\,.
		\end{align}
		Take $V^1_n$ to be the $\amsbb{K}$-span of all  vectors of the form $\big(1,F_n\left(u_1\right)\big) \in \amsbb{K}\times I^{\mu_1}$, where $F_n\in \mathcal{F}^1_n$.    Then, $V_n^1$ is a finite-dimensional vector space for every $n$.  Furthermore, we have
		\begin{align}
			V_1^1\supseteq V_2^1\supseteq V_3^1\supseteq \ldots\,.
		\end{align}
		Hence, the inclusion sequence above stabilizes at some $n_1\in\amsbb{Z}_{>0}$.  That is,
		\begin{align}
			V^1:=V_{n_1}^1=V_{n_1+1}^1=V_{n_1+1}^1=\ldots\,.
		\end{align}
		Since $\mathcal{F}^1_n$ is nonempty for every $n$, we conclude that $V^1$ is nonempty.  As $V^1=V_{n_1}^1$, it must contain $\big(1,F_{n_1}\left(u_1\right)\big)$ for some $F_{n_1}\in\mathcal{F}_{n_1}^1$.  We claim that, for every $n\geq n_1$, there exists $F_n\in\mathcal{F}_n^1$ such that 
		\begin{align}
			F_n\left(u_1\right)=F_{n_1}\left(u_1\right)=:v_1\,.
		\end{align}
		
		To verify the claim above, we observe a trivial fact that $\mathcal{F}_n^1$ is closed under \emph{unit linear combinations}.  That is, if $F_{n,1},F_{n,2},\ldots,F_{n,k}\in \mathcal{F}_n^1$, then 
		\begin{align}
			\sum_{i=1}^k\,t_i\,F_{n,i}\in \mathcal{F}_n^1
		\end{align}
		for every $t_1,t_2,\ldots,t_k\in \amsbb{K}$ with $\sum_{i=1}^k\,t_i=1$.  Therefore, if $\left(1,\tilde{u}_n\right)\in V_n^1$, then 
		\begin{align}
			\tilde{u}_n=\sum_{i=1}^k\,t_i\,F_{n,i}\left(u_1\right)\,,
		\end{align}
		for some $t_1,\ldots,t_k\in\amsbb{K}$ with $\sum_{i=1}^k\,t_i=1$ and for some $F_{n,1},\ldots,F_{n,k}\in\mathcal{F}_n^1$, whence with 
		\begin{align}
			F_n:=\sum_{i=1}^k\,t_i\,F_{n,i}\in\mathcal{F}_n^1\,,
		\end{align} we have $\tilde{u}_n=F_n\left(u_1\right)$.   In particular, for every positive integer $n$, $v_1$ is in the image of $F_n$ for some $F_n\in\mathcal{F}_n^1$.
		
		Now, suppose that $v_1,v_2,\ldots,v_l$ have been obtained so that, for each $n\in\amsbb{Z}_{>0}$, there exists a map $F_n:N_n\to I_n$ such that 
		\begin{align}
			F_n\circ\varphi_n=f_n\text{ and }F_n\left(u_i\right)=v_i\text{ for every }i=1,2,\ldots,l\,.
		\label{eq:injcond}
		\end{align}  
		Let $\mathcal{F}_n^l$ denote the set of all $\Ulie\left(\gfrak'_n\right)$-module homomorphisms $F_n$ that obey (\ref{eq:injcond}).  We proceed as before.  Take $V^{l+1}_n \subseteq \amsbb{K}\times I^{\mu_{l+1}}$ to be the $\amsbb{K}$-span of all vectors of the form $\big(1,F_n\left(u_{l+1}\right)\big)$.  Then, 
		\begin{align}
			V^{l+1}_1\supseteq V^{l+1}_2\supseteq \ldots\,,
		\end{align} so that 
		\begin{align}
			V^{l+1}:=V^{l+1}_{n_{l+1}}=V^{l+1}_{n_{l+1}+1}=\ldots
		\end{align} for some positive integer $n_l$.  Then, $V^{l+1}$ is nonzero and contains $\big(1,F_{n_{l+1}}\left(u_{l+1}\right)\big)$ for some $F_{n_{l+1}}\left(u_{l+1}\right)$.  Then, we set $v_{l+1}$ to be $F_{n_{l+1}}\left(u_{l+1}\right)$.  As before, using the fact that $\mathcal{F}^{l+1}_n$ is closed under unit linear combinations, we conclude that, for every positive integer $n$, there exists $F_n\in \mathcal{F}^{l+1}_n$ for which $F_n\left(u_i\right)=v_i$ for $i=1,2,\ldots,l+1$.

		With known values of $v_1,v_2,\ldots$, we can define $F:N\to I$ via extending the conditions 
		\begin{align}
			F\left(u_i\right)=v_i\text{ for every }i=1,2,\ldots\,.
		\end{align}  This gives a well defined map as $u_1,u_2,\ldots$ generate $N$.  By the construction, $F\circ \varphi=f$, so that $I$ is injective.
	\end{pf}
\end{thm}
	
	\begin{thm}
		Let $\lambda\in\hfrak^*$ be almost dominant.  Then, there exists an injective hull $\Ilie(\lambda)$ of the simple module $\Llie(\lambda)$.  In particular, if $\lambda$ is dominant, then $\Ilie(\lambda)=\Vlie(\lambda)$.
		\label{thm:almostdominant}
		
		\begin{pf}
			For each positive integer $n$, we write $\Llie_n(\lambda)$ for the simple module in $\bggO^{\gfrak'_n}_{\bfrak_n'}$ with highest weight $\lambda\in\hfrak^*$ as in the proof of Proposition~\ref{prop:blocks}, and denote by $\Ilie_n(\lambda)$ its injective hull $\Ilie\left(\lambda;\gfrak'_n,\bfrak'_n,\hfrak\right)$ in $\bggO^{\gfrak'_n}_{\bfrak_n'}$.
			Similarly, $\Mlie_n(\lambda)$ and $\Vlie_n(\lambda)$ are, respectively, the Verma module $\Mlie\left(\lambda;\gfrak'_n,\bfrak'_n,\hfrak\right)$ and the co-Verma module $\Vlie\left(\lambda;\gfrak_n',\bfrak'_n,\hfrak\right)$ in $\bggO^{\gfrak'_n}_{\bfrak'_n}$ with highest weight $\lambda$.  
			
			We have
			\begin{align}
				\ch\big(\Ilie_n(\lambda)\big)=\sum_{\mu_n\succeq \lambda_n}\,\big\{\Ilie_n(\lambda):\Vlie_n\left(\mu\right)\big\}\,\ch\big(\Vlie_n\left(\mu\right)\big)\,.
			\end{align}
			Using BGG Reciprocity, we have
			\begin{equation}
				\big\{\Ilie_n(\lambda):\Vlie_n\left(\mu\right)\big\}=\big[\Mlie_n\left(\mu\right):\Llie_n\left(\lambda\right)\big]\,.
			\end{equation}
			For each $\mu\succeq \lambda$, there exists $n_\mu\in\amsbb{Z}_{> 0}$ (due to Proposition~\ref{prop:MnLn}) such that, for all $n\geq n_\mu$, we have
			\begin{equation}
				\big[\Mlie(\mu),\Llie(\lambda)\big]=\big[\Mlie_n\left(\mu\right):\Llie_n\left(\lambda\right)\big]\,.
			\end{equation}
			Because $\lambda$ is almost dominant, there are finitely many $\mu\succeq\lambda$ with $\mu\in W\cdot \lambda$.  Furthermore, the multiplicity $\big\{\Ilie_n(\lambda):\Vlie_n\left(\mu\right)\big\}$  eventually stabilizes at the value $\big[\Mlie(\mu):\Llie(\lambda)\big]<\infty$.  
			
			We have a sequence of embeddings
			\begin{align}
				\Llie_n\left(\lambda\right)\to\Llie_{n+1}\left(\lambda\right)\to \Ilie_{n+1}(\lambda)\,.
			\end{align}  
			By Proposition~\ref{prop:injrestricted}, $\text{Res}_{\gfrak_n'}^{\gfrak'_{n+1}}\big(\Ilie_{n+1}(\lambda)\big)$ is injective in $\bggO^{\gfrak'_{n}}_{\bfrak'_{n}}$.  Since $\Ilie_n(\lambda)$ is the injective hull of $\Llie_n\left(\lambda\right)$ in $\bggO^{\gfrak'_n}_{\bfrak'_n}$, there exists an embedding $\Ilie_n(\lambda)\to \Ilie_{n+1}(\lambda)$.
			
			From the work above, we conclude that every weight space of $\Ilie(\lambda):=\lim_{\underset{n}{\boldsymbol{\longrightarrow}}}\,\Ilie_n(\lambda)$ is finite-dimensional.  This means $\Ilie(\lambda)\in\bar\bggO$, whence $\Ilie(\lambda)$ injective by the previous proposition.  In particular, if $\lambda$ is already dominant, then $\Ilie_n(\lambda)=\Vlie_n\left(\lambda\right)$ for every $n$.  Since the direct limit of $\Vlie_n\left(\lambda\right)$ is just $\Vlie(\lambda)$, the claim follows.
		\end{pf}
	\end{thm}
	
	\begin{thm}
		For a fixed $\lambda\in\hfrak^*$ and $\mu\preceq\lambda$, $\Llie(\mu)$ has an injective hull in $\bar\bggO^{\preceq \lambda}$.
		
		\begin{pf}
			The proof is similar to that of Theorem~\ref{thm:injdirlim}.  We only need to show that the direct limit 
			\begin{align}
				\Ilie^{\preceq\lambda}(\mu):=\lim_{\underset{n}{\boldsymbol{\longrightarrow}}}\,\tcyr_{\preceq\lambda}\Ilie_n\left(\mu\right)
				\end{align}
			 is in $\bbggO^{\preceq\lambda}$, where $\tcyr_{\preceq\lambda}$ also denotes the truncation functor in $\bggO^{\gfrak'_n}_{\bfrak_n'}$ with upper bound $\lambda\in\hfrak^*$.  To this end, we need to verify that $\Ilie^{\preceq\lambda}(\mu)$ has finite-dimensional weight spaces.  
			
			We say that two formal characters $\xi$ and $\zeta$ satisfies $\xi\leq \zeta$ if all coefficients of $e^{\lambda}$ in $\zeta-\xi$ are nonnegative integers.  By studying the formal character of $\tcyr_{\preceq \lambda}\Ilie_n\left(\mu\right)$, it is easy to see that 
			\begin{align}
				\ch\Big(\tcyr_{\preceq\lambda}\Ilie_n\left(\mu\right)\Big) &\leq \sum_{\nu\in\left[\mu,\lambda\right]}\,\big\{\Ilie_n(\mu):\Vlie_n\left(\nu\right)\big\}\,\ch\big(\Vlie_n\left(\nu\right)\big)
				\nonumber\\
				&=\sum_{\nu\in\left[\mu,\lambda\right]}\,\big[\Vlie_n(\nu):\Llie_n\left(\mu\right)\big\}\,\ch\big(\Vlie_n\left(\nu\right)\big)\,.
				\label{eq:truncchar}
			\end{align}
			The right-hand side of (\ref{eq:truncchar}) is bounded as $n\to\infty$.  Therefore, the direct limit $\Ilie^{\preceq\lambda}\left(\mu\right)$ is indeed an object in $\bar\bggO^{\preceq\lambda}$. 
			
			 Since each $\Ilie_n\left(\mu\right)$ is an essential extension of $\Llie_n\left(\mu\right)$ in $\bggO^{\gfrak'_n}_{\bfrak'_n}$, the truncation $\tcyr_{\preceq\lambda}\Ilie_n\left(\mu\right)$ is also an essential extension of $\Llie_n\left(\mu\right)$ in $\left(\bggO^{\gfrak'_n}_{\bfrak'_n}\right)^{\preceq\lambda}$.  Thus, $\Ilie^{\preceq\lambda}\left(\mu\right)$ is indeed the injective hull of $\Llie(\mu)$ in $\bar{\bggO}^{\preceq\lambda}$.
		\end{pf}
	\end{thm}

Since we have introduced the notion of dominant weights, we can also discuss integrable modules in $\bbggO$.  All integrable modules of $\bbggO$ are classified in Theorem~\ref{thm:integrable}, which is easy to prove.

\begin{define}
	Let $R$ be an associative $\amsbb{K}$-algebra.  A left $R$-module $M$ is \emph{integrable} if, for every finitely generated subalgebra $S$ of $R$ and for any $m\in M$, $S\cdot m$ is a finite-dimensional vector space.
\end{define}

\begin{thm}
	A module $M\in\bbggO$ is integrable if and only if it is a direct sum of simple modules $\Llie(\lambda)$, where each $\lambda$ is an integral dominant weight.
	\label{thm:integrable}
\end{thm}	

\subsection{BGG Reciprocity}

The theorem below is stated and proven in \cite[Theorem 3.5(a)]{hwcat}.  Corollaries~\ref{cor:covermafilt} and~\ref{cor:obarhwcat} are immediate consequences of this theorem and the fact that both $\bggO^{\gfrak_n}_{\bfrak_n}$ and $\bbggO^{\gfrak_n}_{\bfrak_n}$ are highest-weight categories.  Corollary~\ref{cor:bhreciprocity} follows from Corollary~\ref{cor:covermafilt} and Theorem~\ref{thm:brauer-humphreys}.

\begin{thm}
	Let $\pazocal{C}$ be a highest-weight category with notations as in Definition~\ref{def:hwcat}.    For $\lambda\in\mathcal{P}$, let $\pazocal{C}^{\preceq \lambda}$ be the full subcategory of $\pazocal{C}$ consisting of objects each whose composition factors $S(\mu)$ satisfies $\mu\preceq\lambda$.  We define the \emph{truncation functor}
	\begin{align}	
		\textrm{T}_{\preceq\lambda}X:=\sum_{\substack{{Y\subseteq X}\\{Y\in \pazocal{C}^{\preceq\lambda}}}}\,Y\,.
	\end{align}
	Then, $\pazocal{C}^{\preceq\lambda}$ is also a highest-weight category with respect to the locally finite partially ordered set $\mathcal{P}^{\preceq\lambda}:=\left\{\mu\in\mathcal{P}\suchthat{}\mu\preceq\lambda\right\}$, the family $\big\{S(\mu)\big\}_{\mu\preceq \lambda}$ of simple objects , the family $\big\{A(\mu)\big\}_{\mu\preceq \lambda}$ of  co-standard objects, and the family $\big\{\textrm{T}_{\preceq\lambda}I(\mu)\big\}_{\mu\preceq\lambda}$.  The family $\big\{V(\mu)\big\}_{\mu\preceq\lambda}$ consists of all standard objects of $\pazocal{C}^{\preceq\lambda}$.

\end{thm}

\begin{cor}
	Let $\lambda_n\in\hfrak_n^*$.   The category $\left(\bggO^{\gfrak_n}_{\bfrak_n}\right)^{\preceq\lambda_n}$ is a highest-weight category with respect to the partially ordered set $\left\{\mu_n\in\hfrak_n^*\suchthat{}\mu_n\preceq\lambda_n\right\}$, the family $\big\{\Llie\left(\mu_n\right)\big\}_{\mu_n\preceq \lambda_n}$ of simple objects, the family $\big\{\Vlie\left(\mu_n\right)\big\}_{\mu_n\preceq \lambda_n}$ of co-standard objects, and the family $\big\{\tcyr_{\preceq\lambda_n}\Ilie\left(\mu_n\right)\big\}_{\mu_n\preceq \lambda_n}$ of injective objects.  In addition, $\big\{\Mlie\left(\mu_n\right)\big\}_{\mu_n\preceq \lambda_n}$ is the family of standard objects.
	\label{cor:covermafilt}
\end{cor}

\begin{cor}
	Let $\lambda_n\in\hfrak_n^*$.   The category $\left(\bbggO^{\gfrak_n}_{\bfrak_n}\right)^{\preceq\lambda_n}$ is a highest-weight category, also with respect to the partially ordered set $\left\{\mu_n\in\hfrak_n^*\suchthat{}\mu_n\preceq\lambda_n\right\}$, the family $\big\{\Llie\left(\mu_n\right)\big\}_{\mu_n\preceq \lambda_n}$ of simple objects, the family $\big\{\Vlie\left(\mu_n\right)\big\}_{\mu_n\preceq \lambda_n}$ of co-standard objects, and the family $\big\{\tcyr_{\preceq\lambda_n}\Ilie\left(\mu_n\right)\big\}_{\mu_n\preceq \lambda_n}$ of injective objects.  Furthermore, $\big\{\Mlie\left(\mu_n\right)\big\}_{\mu_n\preceq \lambda_n}$ is the family of standard objects.
	\label{cor:obarhwcat} 
\end{cor}

\begin{cor}
	For every $\lambda_n,\mu_n,\nu_n\in\hfrak_n^*$ with $\mu_n,\nu_n\preceq \lambda_n$, we have BGG reciprocity
	\begin{align}
		\Big\{\tcyr^\vee_{\preceq \lambda_n}\Plie\left(\mu_n\right):\Mlie\left(\nu_n\right)\Big\}&=\Big\{\tcyr_{\preceq \lambda_n}\Ilie\left(\mu_n\right):\Vlie\left(\nu_n\right)\Big\}\nonumber
		\\
		&=\Big[\Mlie\left(\nu_n\right):\Llie\left(\mu_n\right)\Big]=\Big[\Vlie\left(\nu_n\right):\Llie\left(\mu_n\right)\Big]\,.
	\end{align}
	\label{cor:bhreciprocity}
\end{cor}


\begin{rmk}
	As in the previous subsection, we write $\gfrak'_n$ and $\bfrak_n'$ for  $\hfrak+\gfrak_n$ and $\hfrak+\bfrak_n$, respectively.  It is important to note that Corollary~\ref{cor:covermafilt}, Corollary~\ref{cor:bhreciprocity}, and Corollary~\ref{cor:obarhwcat} hold for $\gfrak'_n$ and $\bfrak'_n$ as well (with proper modifications such as replacing $\Llie\left(\mu_n\right)$ by $\Llie_n(\mu)$, where $\Llie_n(\mu)$ is defined in the proof of Theorem~\ref{thm:almostdominant}).   
	\label{rmk:gprime}
\end{rmk}

Next, we study the category $\left(\bbggO^{\gfrak}_{\bfrak}\right)^{\preceq\lambda}$ for a fixed $\lambda\in\hfrak^*$.  To do so, we first note that every object $M$ of $\left(\bbggO^{\gfrak}_{\bfrak}\right)^{\preceq\lambda}$ is countable-dimensional.  Therefore, $M$ can be generated by countably many weight vectors $v_1,v_2,\ldots$.   We set 
\begin{align}
	M_n:=\sum_{i=1}^k\,\Ulie\left(\gfrak'_n\right)\cdot v_i\,.
\end{align}    
The injective hull of $M_n$ in $\left(\bggO^{\gfrak'_n}_{\bfrak'_n}\right)^{\preceq\lambda}$ is denoted by $I_n$.

Since we have a $\gfrak'_n$-module embedding $M_n\to M_{n+1}$ and $I_{n+1}$ is injective as an object of $\left(\bggO^{\gfrak'_n}_{\bfrak'_n}\right)^{\preceq\lambda}$, there exists a $\gfrak'_n$-module embedding $I_n\to I_{n+1}$.  The question is now whether the direct limit $I:=\lim_{\underset{n}{\boldsymbol{\longrightarrow}}}\,I_n$ is in $\left(\bbggO^{\gfrak}_{\bfrak}\right)^{\preceq\lambda}$; that is, we need to check whether the weight spaces of $I$ are finite dimensional.

Let $\mu\in\hfrak^*$ be such that $\mu\preceq \lambda$.  We want to find $\dim_\amsbb{K}\left(I^\mu\right)$.  To do this, we find a bound on $\dim_\amsbb{K}\left(I_n^{\mu}\right)$.  There are at most $\dim_\amsbb{K}\left(M^\mu\right)$ indecomposable direct summands of $I_n$ having $\mu$ as a weight.  We focus on one of such indecomposable direct summands.  It is of the form $\tcyr_{\preceq \lambda}\Ilie_n\left(\xi\right)$ for some $\xi\preceq \mu$ (here, the $\gfrak'_n$-module $\Ilie_n(\xi)$, as well as $\Vlie_n(\nu)$, is as defined in the proof of Theorem~\ref{thm:almostdominant}).  

The contribution to the weight space with weight $\mu$ of $\tcyr_{\preceq \lambda}\Ilie_n\left(\xi\right)$ can only come from its co-Verma subquotients $\Vlie_n\left(\nu\right)$ with $\mu\preceq \nu\preceq\lambda$.  Thus, we have an upper bound
\begin{align}
	\dim_\amsbb{K}\left(I_n^{\mu}\right)\leq \sum_{\xi}\,\sum_{\nu\in\left[\mu,\lambda\right]}\,\Big\{\tcyr_{\preceq \lambda}\Ilie_n\left(\xi\right):\Vlie_n\left(\nu\right)\Big\}\,\dim_\amsbb{K}\left(\Vlie_n\left(\nu\right)^{\mu}\right)\,,
\end{align}
where $\xi$ runs over possible weights such that $\tcyr_{\preceq \lambda}\Ilie\left(\xi\right)$ is an indecomposable direct summand of $I_n$ with $\mu$ as a weight, and $[\mu,\lambda]$ denotes the set $\big\{\nu\in\hfrak^*\suchthat{}\mu\preceq\nu\preceq\lambda\big\}$.
By the BGG reciprocity, we have 
\begin{align}
	\dim_\amsbb{K}\left(I_n^{\mu}\right)&\leq \sum_{\xi}\,\sum_{\nu\in\left[\mu,\lambda\right]}\,\Big[\Mlie_n\left(\nu\right):\Llie_n\left(\xi\right)\Big]\,\dim_\amsbb{K}\left(\Mlie_n\left(\nu\right)^{\mu}\right)
	\nonumber \\
	&\leq \dim_\amsbb{K}\left(M^\mu\right)\,\sum_{\nu\in\left[\mu,\lambda\right]}\,A_n(\nu)\,\dim_\amsbb{K}\left(\Mlie\left(\nu\right)^{\mu}\right)\,,
	\label{eq:injbound}
\end{align}
where $A_n(\nu)$ is the maximum possible value of $\Big[\Mlie_n\left(\nu\right):\Llie_n\left(\xi\right)\Big]$ with $\xi\preceq \mu$.

We are now ready to prove the proposition below.

\begin{prop}
	If $\lambda$ is an almost antidominant weight, then the truncated category $\left(\bbggO^{\gfrak}_{\bfrak}\right)^{\preceq\lambda}$ has enough injectives (and so, $\left(\bbggO^{\gfrak}_{\bfrak}\right)^{\preceq\lambda}$ has enough projectives as well).
	
	\begin{pf}
		If $\lambda$ is almost antidominant, then $\mu$ is also almost antidominant.  Therefore, there are finitely many weights $\xi\in\hfrak^*$ such that $\xi\preceq \mu$.  Thus, if $A$ denotes the maximum of $\Big[\Mlie(\nu):\Llie(\xi)\Big]$ with $\xi\preceq \mu$ and $\nu\in[\mu,\lambda]$, we have from (\ref{eq:injbound}) that
\begin{align}
	\dim_\amsbb{K}\left(I_n^{\mu}\right) &\leq \dim_\amsbb{K}\left(M^\mu\right)\,\sum_{\nu\in\left[\mu,\lambda\right]}\,A_n(\nu)\,\dim_\amsbb{K}\left(\Mlie\left(\nu\right)^{\mu}\right)
	\nonumber\\
	&\leq A\,\dim_\amsbb{K}\left(M^\mu\right)\,\sum_{\nu\in\left[\mu,\lambda\right]}\,\dim_\amsbb{K}\left(\Mlie\left(\nu\right)^{\mu}\right)<\infty\,,
\end{align}
whenever $n$ is large enough.  Ergo, there exists a universal bound for the dimension of the weight space $I_n^{\mu}$ for all (sufficiently large) $n$.  That is, $\dim_\amsbb{K}(I^\mu)<\infty$ and the claims follows immediately.
	\end{pf}
\end{prop}

Now, we want to show that, for any $\lambda\in\hfrak^*$ and $\mu\preceq\lambda$, the injective hull $\Ilie^{\preceq\lambda}(\mu)$ has a co-standard filtration.  We recall from the finite-dimensional theory and Remark~\ref{rmk:gprime} that there exists a co-Verma filtration
\begin{align}
	0=F_n^0\subsetneq F_n^1\subsetneq \ldots \subsetneq F^{t_n-1}_n\subsetneq F^{t_n}_n=\Ilie_n^{\preceq \lambda}\left(\mu\right)\,,
\end{align}
where $\Ilie_n^{\preceq\lambda}\left(\mu\right)$ denotes the $\gfrak'_n$-module $\tcyr_{\preceq\lambda}\,\Ilie_n(\mu)$.  Since the highest weights of $F_n^{i}/F^{i-1}_n$ are in the interval $\left[\mu,\lambda\right]$, we have by BGG reciprocity that
\begin{align}	
	t_n \leq \sum_{\nu\in\left[\mu,\lambda\right]}\,\Big[\Mlie_n\left(\nu\right):\Llie_n\left(\mu\right)\Big] =  \sum_{\nu\in\left[\mu,\lambda\right]}\,\Big[\Mlie\left(\nu\right):\Llie\left(\mu\right)\Big]<\infty
\end{align}
for sufficiently large $n$.  Thus, there exists a sequence $\left(n_k\right)_{k=1}^\infty$ of positive integers such that 
\begin{align}
	n_1<n_2<n_3\ldots
\end{align} and 
\begin{align}
	t:=t_{n_1}=t_{n_2} =t_{n_3}=\ldots\,.
\end{align}   
Furthermore, as there are only finitely many weights $\upsilon \in W[\lambda]\cdot \lambda$ with $\mu\preceq\upsilon\preceq \lambda$, we may assume without loss of generality (due to the Pigeonhole Principle) that the highest weight $\xi_{n_k}[i]$ of $F^i_{n_{k+1}}/F^{i-1}_{n_{k+1}}$ is the same as the highest weight of $F^{i}_{n_k}/F^{i-1}_{n_k}$ for every $i$ and $k$.  Denote by $\xi[i]\in\hfrak^*$ the common weight $\xi_{n_1}[i],\xi_{n_2}[i],\xi_{n_3}[i],\ldots$.  We need the following lemma.

\begin{lem}
	Let $m,n\in\amsbb{Z}_{>0}$ be such that $m\leq n$.  The Verma module $\Mlie_m\left(\lambda\right)$ is a direct summand of the Verma module $\Mlie_n\left(\lambda\right)$, viewed as a $\gfrak'_m$-module.  Consequently, the co-Verma module $\Vlie_m\left(\lambda\right)$ is also a direct summand of the co-Verma module $\Vlie_n\left(\lambda\right)$, viewed as a $\gfrak_m'$-module.
	\label{lem:covermaemb}
	
	\begin{pf}
		Let $v$ be a highest-weight vector of $\Mlie_n\left(\lambda\right)$.  Let 
		\begin{align}
			\left\{x_{\pm\alpha}\,\boldsymbol{|}\,\alpha\in\Delta^+\right\}\cup\left\{h_\beta\,\boldsymbol{|}\,\beta\in\Sigma^+\right\}
		\end{align}
		be a Chevalley basis of $\gfrak$.  Then, 
		\begin{align}
			\text{Res}^{\gfrak_n'}_{\gfrak_m'}\Mlie_n\left(\lambda\right)=\Ulie\left(\gfrak_n'\right)\cdot v= \Big(\Ulie\left(\gfrak'_m\right)\cdot v \Big)\oplus \left(\sum_{\alpha}\,\Ulie\left(\gfrak'_n\right)\cdot \left(x_{-\alpha}\cdot v\right)\right)\,,
		\end{align}
		where $\alpha$ runs over $\bfrak_n'$-positive roots of $\gfrak'_n$ which are not roots of $\gfrak'_m$, is a direct sum decomposition of $\Mlie_n\left(\lambda\right)$ as a $\gfrak'_m$-module with a direct summand 
		\begin{align}
			\Ulie\left(\gfrak'_m\right)\cdot v\cong \Mlie_m\left(\lambda\right)\,.
		\end{align}  To prove the co-Verma version, we only need to apply the duality functor.
	\end{pf}
\end{lem}

Clearly, $F^1_{n_k}$ is a co-Verma submodule of $I_k:=\Ilie^{\preceq \lambda}_{n_k}\left(\mu\right)$ containing the socle of $I_k$.  By the lemma below, each $F^1_{n_k}$ is unique as it contains the simple module 
\begin{align}
\Ulie\left(\gfrak_{n_k}\right)\cdot v[1]\cong\Llie\big(\xi[1]\big)\,,
\end{align} where $v[1]$ is a singular vector of the socle of $\Ilie^{\preceq\lambda}(\mu)$.   Then, using Lemma~\ref{lem:covermaemb}, we have a embeddings $F^1_{n_k}\to F^1_{n_{k+1}}$, whose direct limit is simply the co-Verma module \begin{align}
	F^1\cong\Vlie\big(\xi[1]\big)\,,
\end{align}
where $\xi[1]$ is clearly equal to $\mu$. 

\begin{lem}
	Let $n\in\amsbb{Z}_{>0}$ and $M_n\in\bggO^{\gfrak'_n}_{\bfrak'_n}$.  Suppose that a simple module $L_n\in\bggO^{\gfrak'_n}_{\bfrak'_n}$ is a submodule of $M_n$.  Then, $M_n$ has at most one co-Verma submodule $V_n$ such that $L_n\subseteq V_n\subseteq M_n$.
	\label{lem:uniquecoverma}
	
	\begin{pf}
		Suppose that $M$ has two co-Verma submodules $V_n$ and $V'_n$ with $L_n\subseteq V_n$ and $L_n\subseteq V'_n$.  Take $N_n:=V_n+V_n'$.  Then, $N_n$ is indecomposable (as $V_n$ and $V'_n$ are both indecomposable with $V_n\cap V'_n\supseteq L_n\supsetneq 0$).  Hence, we have an exact sequence
		\begin{align}
			0 \to V_n \to N_n \to N_n/V_n \to 0\,.
		\end{align}
		By dualizing the exact sequence above, we have
		\begin{align}
			0\to \left(N_n/V_n\right)^\vee \to N_n^\vee \to V_n^\vee \to 0\,.
		\end{align}
		By Proposition~\ref{prop:extintro}, we see that this exact sequence must split.  As $N_n^\vee$ is indecomposable, we conclude that $V_n^\vee=0$ or $\left(N_n/V_n\right)^\vee=0$.  Since $V_n\neq 0$, we must have $N_n/V_n=0$, which leads to $V'_n=V_n$.
	\end{pf}
\end{lem}

Suppose now that, for some positive integer $l<t$, the submodules $0=F^0$, $F^1$, $F^2$, $\ldots$, $F^l$ of $\Ilie^{\preceq \lambda}(\mu)$ have been determined with the property that $F^i$ is the direct limit $\lim\limits_{\longrightarrow}\,F_{n_k}^i$, where the $\gfrak'_{n_k}$-modules $F_{n_k}^i$ are submodules of $I_k$ satisfying the following properties:
\begin{enumerate}[(i)]
\item $0=F^0_{n_k}\subsetneq F^1_{n_k} \subsetneq \ldots \subsetneq F^l_{n_k}\subsetneq I_k$,
	\item $F^{i}_{n_k}/F^{i-1}_{n_k}\cong\Vlie\big(\xi[i]\big)$ for every $i=1,2,\ldots,l$.
\end{enumerate}Then, we proceed by looking at the quotient $I_k/F^{l}_{n_k}$.  Identify each $u+F^l_{n_k}\in I_k/F^l_{n_k}$ as an element of $I/F^l$ via
\begin{align}
	u+F^l_{n_k}\mapsto u+F^l\in I/F^l
	\label{eq:identify}
\end{align} 
(making $I_k/F_{n_k}^l$ a $\gfrak'_{n_k}$-submodule of $I/F^l$).  We have an embedding 
\begin{align}
	\Vlie\big(\xi[l+1]\big) \to I_k/F^l_{n_k}
\end{align} 
for each $k$.  Let $V^{l+1}_k$ be the $\amsbb{K}$-span of all vectors $v\in I_k/F^{l}_{n_k} \subseteq I/F^{l}$ such that $v$ is the image of a singular vector under an embedding $\Vlie\big(\xi[l+1]\big)\to I_k/F^{l+1}_{n_k}$.  Hence, $V^{l+1}_k$ is a nonzero subspace of $\left(I/F^l\right)^{\xi[l+1]}$ and $V^{l+1}_k\supseteq V^{l+1}_{k+1}$ for every $k$.  Because 
\begin{align}
	\dim_\amsbb{K}\left(V_k^{l+1}\right)\leq \dim_\amsbb{K}\left(\left(I/F^l\right)^{\xi[l+1]}\right)<\infty\,,
\end{align} 
there exists $v[l+1]+F^l\in \bigcap_{k\in\amsbb{Z}_{>0}}\,V^{l+1}_k$ which is a nonzero element of $\left(I/F^l\right)^{\xi[l+1]}$. 

Now, by Lemma~\ref{lem:uniquecoverma}, we can show that there is a unique co-Verma submodule $\bar{F}^{l+1}_{n_k}$ of $I_k/F_{n_k}^l$ containing the simple submodule $\Ulie\left(\gfrak'_{n_k}\right)\cdot\left(v[l+1]+F^l\right)$.  Then, the direct limit $\bar{F}^{l+1}$ of $\bar{F}^{l+1}_{n_k}$ must be a co-Verma module of highest weight $\xi[l+1]$.  Let $F^{l+1}$ be the preimage of $\bar{F}^{l+1}$ under the quotient map $I\to I/F^l$.   Then, by induction, we have found a filtration
\begin{align}
	0=F^0\subsetneq F^1\subsetneq F^2\subsetneq \ldots\subsetneq F^{t-1}\subsetneq F^t=\Ilie^{\preceq\lambda}(\mu)
	\label{eq:covermafiltind}
\end{align}
of $\Ilie^{\preceq \lambda}(\mu)$ such that each successive quotient $F^l/F^{l-1}$ is isomorphic to the co-Verma module $\Vlie\big(\xi[l]\big)$.   It can be easily seen that the number of times a co-Verma module $\Vlie(\nu)$ appears as a successive quotient $F^l/F^{l-1}$ in (\ref{eq:covermafiltind}) is independent on the choice of the co-Verma filtration.   We use the notation $\Big\{\Ilie^{\preceq\lambda}(\mu):\Vlie(\nu)\Big\}$ for the number of times that $\Vlie(\nu)$ appears as a successive quotient in (\ref{eq:covermafiltind}). 

Let $\Plie^{\preceq \lambda}(\mu)$ denote $\tcyr^\vee_{\preceq \lambda}\Plie(\mu)=\Big(\tcyr_{\preceq\lambda}\Ilie(\mu)\Big)^\vee$.  Then, by applying duality on the co-Verma filtration (\ref{eq:covermafiltind}), $\Plie^{\preceq \lambda}(\mu)$ has a Verma filtration
\begin{align}
	0=T^0\subsetneq T^1\subsetneq T^2\subsetneq \ldots\subsetneq T^{t-1}\subsetneq T^t=\Plie^{\preceq\lambda}(\mu)\,,
	\label{eq:vermafiltind}
\end{align}
where each successive quotient $T^l/T^{l-1}$ is isomorphic to the Verma module $\Mlie\big(\xi[t+1-l]\big)$.  In particular, $T^t/T^{t-1}\cong \Mlie(\mu)$.  The number of times that $\Mlie(\nu)$ appears as a successive quotient $T^l/T^{l-1}$ in (\ref{eq:vermafiltind}) is also well defined, and is denoted by $\Big\{\Plie^{\preceq\lambda}(\mu):\Mlie(\nu)\Big\}$.

\begin{prop}
	For every $\lambda,\mu\in\hfrak^*$ with $\mu\preceq \lambda$, the injective object $\Ilie^{\preceq\lambda}(\mu)$ has a finite co-standard filtration as in the definition of highest-weight categories (Definition~\ref{def:hwcat}).  Furthermore, we have BGG reciprocity:
	\begin{align}
		\Big\{\Plie^{\preceq\lambda}(\mu):\Mlie(\nu)\Big\}=\Big\{\Ilie^{\preceq\lambda}(\mu):\Vlie(\nu)\Big\}=\Big[\Mlie(\nu):\Llie(\mu)\Big]=\Big[\Vlie(\nu):\Llie(\mu)\Big]\,,
	\end{align}
	for all $\nu\in[\mu,\lambda]$.  
	\label{prop:covermafiltration}
\end{prop}

Finally, we note that, if $\lambda$ is not almost antidominant, then $\Mlie(\lambda)$ is of infinite length and cannot be written as a union of subobjects of finite length.  This is because every submodule $M$ of $\Mlie(\lambda)$ has a singular vector $v\neq 0$.  The submodule $N$ of $M$ generated by $v$ is then a Verma module with highest weight $\mu\preceq\lambda$, which is not almost antidominant.  Ergo, $N$ is of infinite length, and so is $M$.  Thus, $\Mlie(\lambda)$ has no submodules of finite length.  In particular, this implies that $\Mlie(\lambda)$ has trivial socle.

The argument above shows that $\left(\bbggO^{\gfrak}_{\bfrak}\right)^{\preceq \lambda}$ is not locally artinian, whence this category does not satisfy the condition (i) of Definition~\ref{def:hwcat}).  That is, $\left(\bbggO^{\gfrak}_{\bfrak}\right)^{\preceq \lambda}$ is not a highest-weight category.  Combining this observation with the fact that $\left(\bbggO^{\gfrak}_{\bfrak}\right)^{\preceq \lambda}$ has enough injectives when $\lambda$ is almost antidominant, we conclude the following theorem.

\begin{thm}
	The category $\left(\bbggO^\gfrak_\bfrak\right)^{\preceq \lambda}$ is a highest-weight category with respect to the partially ordered set $\left\{\mu\in\hfrak^*\suchthat{}\mu\preceq\lambda\right\}$, the family $\big\{\Llie\left(\mu\right)\big\}_{\mu\preceq \lambda}$ of simple objects, the family $\big\{\Vlie\left(\mu\right)\big\}_{\mu\preceq \lambda}$ of co-standard objects, and the family $\big\{\Ilie^{\preceq\lambda}\left(\mu\right)\big\}_{\mu\preceq \lambda}$ of injective objects if and only if $\lambda$ is an almost antidominant weight.  In the case where $\lambda$ is almost antidominant, $\big\{\Mlie\left(\mu\right)\big\}_{\mu\preceq \lambda}$ is the family of standard objects of the highest-weight category $\left(\bbggO^\gfrak_\bfrak\right)^{\preceq \lambda}$.

\end{thm}

\begin{openq}
For a weight $\lambda\in\hfrak^*$ which is not almost antidominant, does the category $\left(\bbggO^{\gfrak}_{\bfrak}\right)^{\preceq\lambda}$ have enough injectives?
\end{openq}

\begin{define}
	Let $\mathcal{C}$ be an abelian category with an abelian subcategory $\tilde{\mathcal{C}}$.  An object $I\in\mathcal{C}$ is \emph{injective relative to $\tilde{\mathcal{C}}$} if, for any two objects $X,Y\in\tilde{\mathcal{C}}$ and any monomorphism $f\in \Hom_{\tilde{\mathcal{C}}}(X,Y)$, every morphism $g\in \Hom_{\mathcal{C}}(X,I)$ factors through $f$, i.e., there exists $\varphi\in\Hom_{\mathcal{C}}(Y,I)$ such that 
	\begin{align}
		g=\varphi\circ f\,.
	\end{align}
\end{define}

\begin{thm}
	Let $R$ be a ring.  Suppose that $\mathcal{C}$ and $\tilde{\mathcal{C}}$ are abelian subcategories of the category of left $R$-modules with $\tilde{\mathcal{C}}$ being a subcategory of $\mathcal{C}$.  If $M\in\tilde{\mathcal{C}}$ has an injective hull $I$ in $\tilde{\mathcal{C}}$, then for each object $J\in\mathcal{C}$ which is injective relative to $\tilde{\mathcal{C}}$, any embedding $\iota\in\Hom_{\mathcal{C}}(M,J)$ induces an embedding $\varphi\in\Hom_\mathcal{C}(I,J)$.
	\label{thm:relinj}
	
	\begin{pf}
		We have an exact sequence $0\to M\to I$ of objects and morphisms in $\tilde{\mathcal{C}}$ and a homomorphism $\iota\in\Hom_\mathcal{C}(M,J)$.  As $J$ is injective relative to $\tilde{\mathcal{C}}$, there exists a map $\varphi\in\Hom_{\mathcal{C}}(I,J)$ such that the diagram below commutes:
		\begin{equation}
		\begin{tikzcd}
		0\arrow{r}& M \arrow{r}{\subseteq} \arrow{d}{\iota} & I\arrow{dl}{\varphi}\\
		&J \,.&
		\end{tikzcd}
		\label{eq:relinj}
		\end{equation}
		We claim that $\varphi:I\to J$ is an embedding. 
		
		Let $K:=\ker(\varphi)$.  Since $\varphi|_M=\iota$ due to commutativity of (\ref{eq:relinj}) and $\iota$ is an embedding, we must have 
		\begin{align}
			K\cap M=\ker\left(\varphi|_M\right)=\ker(\iota)=0\,.
		\end{align}  Because $I$ is an essential extension of $M$, the condition $K\cap M=0$ implies that $K=0$.  Therefore, $\varphi$ is injective.
	\end{pf}
\end{thm}

	\begin{thm}
		For $\lambda\in\hfrak^*$, the simple module $\Llie(\lambda)$ has an injective hull and a projective cover in $\bbggO$ if and only if $\lambda$ is almost dominant.   In particular, this implies that $\bbggO$ does  not have enough injectives, and therefore, $\bbggO$ is not a highest-weight category.
		
		\begin{pf}
			If $\lambda$ is almost dominant, then Theorem~\ref{thm:almostdominant} shows that $\Llie(\lambda)$ has an injective hull in $\bbggO$, and by duality, it has also a projective cover.  To prove the converse, we suppose on the contrary that $\lambda$ is not almost dominant but $\Llie(\lambda)$ has an injective hull $I$ in $\bbggO$.  
			
			As $\lambda$ is not almost dominant, there exists a sequence of weights $\left(\lambda_i\right)_{i=0}^\infty$ with $\lambda_i\in W[\lambda]\cdot \lambda$ and 
			\begin{align}
				\lambda =\lambda_0 \prec \lambda_1 \prec \lambda_2 \prec \ldots\,.
			\end{align}
			For simplicity, let $I_i$ denote $\Ilie^{\preceq \lambda_i}\left(\lambda\right)$ for $i=0,1,2,\ldots$.  It is clear that $I$ is injective relative to $\bbggO^{\preceq\lambda_i}$ for each $i$.  By Theorem~\ref{thm:relinj}, there exists an embedding of $I_i$ into $I$.
			
			Now, using Proposition~\ref{prop:covermafiltration}, we know that each $I_i$ has a co-Verma filtration 
			\begin{align}
				\Vlie(\lambda)=F_i[0]\subsetneq F_i[1] \subsetneq \ldots \subsetneq F_i\left[k_i\right]=I_i\,.
				\label{eq:covermafilt}
			\end{align}
			Furthermore, as $I_i \subsetneq I_{i+1}$, we have $k_0< k_1< k_2< \ldots$.  For every $j=1,2,\ldots,k_i$, the successive quotient $F_i[j]/F_i[j-1]$ is isomorphic to the co-Verma module $\Vlie\left(\mu_i[j]\right)$ for some $\mu_i[j]\in W[\lambda]\cdot\lambda$ with $\mu_j\succeq \lambda$.  This implies $\dim_\amsbb{K}\left(\big(F_i[j]/F_i[j-1]\big)^\lambda \right)\geq 1$.  Ergo,
			\begin{align}
				\dim_\amsbb{K}\left(I^\lambda\right) \geq \dim_\amsbb{K}\left(I_i^\lambda\right) \geq \sum_{j=1}^{k_i}\,\dim_\amsbb{K}\left(\big(F_i[j]/F_i[j-1]\big)^\lambda \right)\geq k_i
			\end{align}
			for every $i=0,1,2,\ldots$.  As $\lim\limits_{i\tendsto\infty}\,k_i=\infty$, we conclude that $\dim\left(I^\lambda\right)=\infty$, which is absurd.  Hence, $\Llie(\lambda)$ does not have an injective hull in $\bbggO$.  

			Using duality, we also conclude that $\Llie(\lambda)$ does not have a projective cover in $\bbggO$.  The theorem follows.
		\end{pf}
	\end{thm}

\begin{rmk}
	The theorem above also shows that, if $\lambda$ is not almost dominant, then $\Llie(\lambda)$ cannot be embedded into any injective object of $\bbggO$.  In this case, there does not exist a projective object $P\in\bbggO$ along with an surjective $\gfrak$-module homomorphism $P\to \Llie(\lambda)$.
\end{rmk}

\begin{exm}
	When $\gfrak:=\sllie_\infty$, $\hfrak:=\hfrak_{\text{A}}$, and $\bfrak:=\bfrak_{\text{1A}}$, the weight
	\begin{align}
		\lambda:=-2\rho=(0,2,4,6,\ldots)
	\end{align} is not almost dominant.  Thus, $\Llie(\lambda)$ does not have an injective hull, nor does it have a projective cover.
\end{exm}

The theorem above is the reason why we need to truncate the category $\bbggO$.  With truncation, every simple object has an injective hull and a projective cover, and because of that, a version of BGG reciprocity holds.  

\pagebreak

\newpage




\pdfbookmark[0]{Back}{Back}
\renewcommand\refname{References}

\newpage

\printnomenclature%

\newpage

\printindex

\newpage

\blankpage
\end{document}